%% file: SIM_SH.tex
\documentclass[3p,times,english]{elsarticle}
\usepackage{graphicx}
\usepackage{amssymb,amsmath,amsthm,amsbsy}
\usepackage{calligra,calrsfs,mathrsfs}
\usepackage{multirow}
\usepackage{cases}
\usepackage{pgfplots}
\usepackage{tikz}
\usepackage{tkz-euclide}
\usetikzlibrary[patterns]
\usepackage{hyperref}
\usepackage{algorithm}
\usepackage{algpseudocode}
\usepackage{dsfont}
\usepackage{caption, comment}
\usepackage{subcaption}
\usepackage{epstopdf}
\usepackage{epsfig}
\usepackage{array}
\newcolumntype{P}[1]{>{\centering\arraybackslash}p{#1}}
\usepackage{float}
\usepackage[T1]{fontenc}
\usepackage[latin9]{inputenc}
\usepackage{geometry}
\usepackage{color}
\usepackage{stmaryrd}
\usepackage{setspace}
\usepackage{amsthm}
\usepackage{textcomp}
\usepackage{hyperref}
\usepackage{soul,xcolor}
\usepgfplotslibrary{fillbetween}
\usepackage{ulem}
\usepackage{textcomp}
\usepackage[title]{appendix}
\setstcolor{red}
\usepackage{subcaption}
\allowdisplaybreaks
\makeatletter

\newcolumntype{M}[1]{>{\centering\arraybackslash}m{#1}}

\input{format_definitions}

\theoremstyle{plain}

\theoremstyle{remark}
\newtheorem{remark}{Remark}

\biboptions{sort&compress}
\allowdisplaybreaks

\makeatother

\usepackage{babel}

\journal{Journal of Computational Physics}

\begin{document}

\setlength{\unitlength}{1cm}

\begin{frontmatter}

\author[duke,llnl]{Nabil M. Atallah}
\ead{atallah1@llnl.gov}
\author[llnl]{Ketan Mittal}
\ead{mittal3@llnl.gov}
\author[duke]{Guglielmo Scovazzi}
\ead{guglielmo.scovazzi@duke.edu}
\author[llnl]{Vladimir Z. Tomov\corref{ca}}
\ead{tomov2@llnl.gov}
\address[duke]{Department of Civil and Environmental Engineering,
               Duke University, Durham, North Carolina 27708, USA}
\address[llnl]{Lawrence Livermore National Laboratory, Livermore (CA), USA}
\cortext[ca]{Corresponding author: Vladimir Tomov}

\title{A High-Order Shifted Interface Method for \\
       Lagrangian Shock Hydrodynamics}

\begin{abstract}

We present a new method for two-material Lagrangian hydrodynamics, which
combines the Shifted Interface Method (SIM) with a high-order Finite Element Method.
Our approach relies on an exact (or sharp) material interface representation,
that is, it uses the precise location of the material interface.
The interface is represented by the zero level-set of a continuous high-order
finite element function that moves with the material velocity.
This strategy allows to evolve curved material interfaces inside curved elements.
By reformulating the original interface problem over a surrogate (approximate)
interface, located in proximity of the true interface,
the SIM avoids cut cells and the associated problematic issues regarding
implementation, numerical stability, and matrix conditioning.
Accuracy is maintained by modifying the original interface conditions using
Taylor expansions.
We demonstrate the performance of the proposed algorithms on established
numerical benchmarks in one, two and three dimensions.
\end{abstract}


\end{frontmatter}

\input{intro}    
\input{prelims}
\input{SIM_general}

\input{SIM_laghos}
\input{results}


\section{Conclusions}
\label{sec_concl}

We proposed the WSIM, a high-order, two-material Lagrangian
shock hydrodynamics algorithm.
The proposed approach utilizes the exact location of the material interface to
evolve the system, while maintaining a continuous curved interface
compatible with curvilinear grids.
The interface evolution does not involve any geometric operations and is
general with respect to dimension, mesh curvature, and FE discretization order.
We demonstrated that WSIM is capable of
(i) representing evolving interfaces to the order of accuracy of the
underlying high-order discretization spaces,
(ii) equilibrating pressure values for challenging 1D benchmark problems, and
(iii) maintaining sharp material interfaces inside curved elements
for complex two- and three-dimensinoal shock problems. 
For these reasons, WSIM can represent a new pathway in the
(sharp) representation of moving interfaces.

Our future efforts will concentrate on expanding into the
Arbitrary Lagrangian-Eulerian framework.
This extension holds the most promise, particularly in eliminating diffusion
within the material interface regions by remapping the level set function.
Additionally, our plans include extending the methodology to
accommodate an arbitrary number of materials.
Furthermore, we intend to combine the WSIM with the Shifted Boundary Method
to weakly enforce wall boundary conditions in complex geometries.

\section*{Acknowledgments}
The authors of Duke University are gratefully thanking the generous support of Lawrence Livermore National Laboratories, through a Laboratory Directed Research \& Development (LDRD) Agreement.
This work performed under the auspices of the U.S. Department of Energy
by Lawrence Livermore National Laboratory under
Contract DE-AC52-07NA27344, LLNL-JRNL-854557.
Guglielmo Scovazzi has also been partially supported by the National Science Foundation, Division of Mathematical Sciences (DMS), under Grant 2207164.

\appendix
\input{appendix}

\bibliographystyle{plain}
\bibliography{./SIM_SH}

\end{document}

%% file: format_definitions.tex
\newcommand{\ifcomment}{\iffalse}

\newcommand{\RN}[1]{\textup{\uppercase\expandafter{\romannumeral#1}}}

\newcommand{\bs}[1]{\boldsymbol{#1}}

\newcommand{\jump}[1] {\llbracket #1 \rrbracket}

\newcommand{\avg  }[1]{\langle #1 \rangle}

\newcommand{\avgb}[1]{\{\!\{ #1 \}\!\}}    

\newcommand{\ti}[1]{\tilde{#1}}
\newcommand{\G}{\Gamma}

\newcommand{\Om}{\Omega}

\newcommand{\cT}{{\mathcal T}}
\newcommand{\tO}{\tilde{\Omega}^h}

\newcommand{\tG}{\tilde{\Gamma}^{h}}

\newcommand{\tx}{{\tilde{\bs{x}}}}

\newcommand{\tS}{S_{h}}
\DeclareRobustCommand{\tspsb}[2]{{%
		\m@th\ensuremath{%
			^{#1}%
			_{#2}%
		}%
}}
\DeclareRobustCommand{\tsp}[1]{{%
		\m@th\ensuremath{%
			^{#1}%
		}%
}}

\DeclareRobustCommand{\tsb}[1]{{%
		\m@th\ensuremath{%
			_{#1}%
		}%
}}

\newdefinition{rem}{Remark}

%% file: intro.tex
\section{Introduction}
\label{sec:intro}

Lagrangian methods for compressible multi-material shock hydrodynamics
\cite{VonNeumann1950, Benson1992, Loubere2004, Scovazzi2007, Maire2009,
scovazzi2008multi, GScovazzi:2012a, Dobrev2012, Morgan2019, Tokareva2020}
are characterized by a computational mesh that moves with the material velocity.
A key advantage of these methods is that they do not produce any numerical
dissipation around material interfaces, when these are aligned with the
element boundaries of the moving mesh.
The deforming mesh, however, inevitably deteriorates in quality, which can lead
to small time steps or simulation breakdowns.
This motivates the use of Arbitrary Lagrangian/Eulerian (ALE) methods
\cite{Hirt1974, Maire2007, ALOrtega:2011a, Waltz2014, Barlow2016,
Boscheri2016, Tomov2018, Gaburro2020},
where the solution fields are periodically projected on a better mesh.
The new mesh generally loses the alignment with the material interfaces,
producing cut elements that contain multiple materials.
The focus of this paper is how to evolve such mixed elements while maintaining
the key advantage of Lagrangian methods, that is, no dissipation around
material interfaces.

Methods that process mixed elements, or \textit{closure models}
\cite{Kamm2011}, can be classified into two major categories.
The first treats a mixed element as a whole, without considering an explicit
sub-element interface location \cite{Miller2009}.
Material concentrations are represented by volume fractions, and all
material-specific quantites are evolved by introducing an equilibration
mechanism for certain quantity, e.g., pressure
\cite{Yanilkin13, Waltz2021, Waltz2023},
pressure and viscosity as in Section 3.1 of \cite{Kamm2011},
velocity increments as in Section 4.4 of \cite{Yanilkin13},
pressure and heat changes \cite{Despres2007}.
The second approach makes use of material interface
reconstruction within the mixed elements
\cite{Youngs1982, Shashkov2023a, Shashkov2023b}.
Most methods in this category employ of acoustic Riemann solvers in order
to predict quantities like interface velocity, material volume changes,
and pressure \cite{Kamm2010, Barlow2014}.

High-order Finite Element Methods (FEMs) for Lagrangian hydrodynamics on curved
grids have recently gained popularity as they allow both high physics
resolution \cite{Dobrev2012,Guermond2016,Tomov2018}, and efficient computational
performance on the latest computer architectures \cite{Kolev2021, Vargas2022}.
The initial approach to evolve mixed elements in these method was developed
in \cite{Tomov2016}, utilizing the first of the above closure approaches,
mainly because reconstructing curved interfaces within curved elements was
not an existing technology at the time.
While results have been good in most problems, the diffuse-interface nature
of the method tends to propagate small volume fractions away from the
contact regions, especially in simulations involving many ALE steps
over large material displacements \cite{Tomov2016, Tomov2018}.
Diffuse-interface methods that reduce the size of the diffused regions  have
been developed \cite{Xiao2011, Waltz2021, Waltz2023}, but the interface is
still approximate and some amount of diffusion is still present.

Recently, the Shifted Interface Method (SIM)~\cite{Scovazzi2020} was proposed, as an alternative FEM for interface problems.
By reformulating the original interface problem over a surrogate (approximate)
interface, located in proximity of the true interface, the SIM avoids cut cells and the associated
problematic issues regarding implementation, numerical stability, and matrix conditioning.
Accuracy is maintained by modifying the original interface conditions using Taylor expansions.

The SIM was derived from the Shifted Boundary Method (SBM)~\cite{main_shifted_2018-1,main_shifted_2018},
which uses the same strategy for the approximation of boundary conditions.
The SBM was applied to simple hyperbolic systems of conservations laws, such as the equations of acoustics and shallow-water flows~\cite{song_shifted_2018}.
Some preliminary developments in combining the SBM with high-order discretizations were also made in~\cite{atallah2022high}, for the Poisson and Stokes problems.

The SIM approach preserves all favorable properties of high-order FE
discretizations (optimal convergence, generality with respect to space dimension,
element type, and FE polynomial order), while allowing a notion
of exact interface representation without any explicit geometric operations.
The generality of the SIM makes it an obvious approach for evolving exact
curved interfaces in the framework of \cite{Dobrev2012}, which is the goal of
this work.

In this paper we transform the single-material high-order FEM of
\cite{Dobrev2012} into a two-material Lagrangian method that evolves mixed
elements through the SIM.
The material interface is represented by the zero level set (LS) of a continuous
high-order finite element function $\eta$ that moves with the material velocity,
i.e., $\eta$ is constant in time in the Lagrangian frame.
The zero LS of $\eta$, which from now on will also be called the {\it true interface},
 is used to define a {\it surrogate interface} where the interface conditions
will be applied; these surrogates are a subset of the mesh faces.
Distance vectors between the surrogate and true interface
(the zero LS of $\eta$) are computed by a FE distance solver.
The interface conditions are imposed weakly by defining
face integrals over the surrogate faces.
Taylor expansions along the distance vectors are used to compensate for
the difference in location between the true and surrogate interface conditions.
The above procedure allows the method to have a notion of exact (or sharp)
interface location, that is, it uses a precise interface location in the
material-related calculations.
The interface evolution does not involve any geometric operations and is
general with respect to dimension, mesh curvature, and FE discretization order.

The rest of the paper is organized as follows.
Section \ref{sec_prelim} reviews the equations of Lagrangian hydrodynamics,
their FE discretization in \cite{Dobrev2012}, and the targeted interface
conditions.
Section~\ref{sec:gen} introduces the general concepts of the SIM,
starting in a broader context, and then moving to the specifics of
of two-material shock hydrodynamics.
Then Section \ref{sec_sim} gives the FE discretization details in the case of
interest, derives the SIM formulation, and discusses
its conservation properties.
Numerical tests are presented in Section \ref{sec_results}.
Finally Section \ref{sec_concl} provides some conclusions and future work.

%% file: prelims.tex

\section{General equations and finite element discretization}
\label{sec_prelim}

This work builds on an existing single-material Lagrangian method~\cite{Dobrev2012}.
In this section, we briefly go over the underlying equations and their discretization in the finite element
framework; the complete description can be found in~\cite{Dobrev2012}.

\subsection{Governing equations}
\label{sec_equations_SH}

The classical equations of Lagrangian shock hydrodynamics govern the rate of
change in position $\bs{x}$, velocity $\bs{v}$, and specific internal energy $e$ of a
compressible body of fluid as it deforms.
Let $\Omega_0$ and $\Omega$ be open sets in $\mathbb{R}^{d}$ where $d$ is
the number of spatial dimensions.
The motion map
\begin{equation}
\label{eq:motion}
\bs{\varphi}: \Omega_0 \rightarrow \Omega, \quad
\bs{\varphi}(\bs{x}_0, t) = \bs{x} \; , \quad \forall \bs{x}_0 \in \Omega_0 \; , \forall t \geq 0 \; ,
\end{equation}
maps the material coordinate $\bs{x}_0$, representing the initial position of an
infinitesimal material particle of the body, to $\bs{x}$, the position
of the same particle in the current configuration, see Figure \ref{fig_lagmap}.
Using the deformation gradient $\bs{F}(\bs{x}, t) = \nabla_{\bs{x}_0} \bs{\varphi}(\bs{x}, t)$ and
its corresponding Jacobian $J(\bs{x},t) = \det(\bs{F})$, the material density $\rho$ is
evolved pointwise by:
\begin{equation}
\label{eq:mass}
  \rho(\bs{x}, t) = \frac{\rho_0(\bs{x}_0)}{J(\bs{x},t)} \; .
\end{equation}
The material velocity $\bs{v}$ and position $\bs{x}$ are represented by finite
element functions in the space $\mathcal{V} \subset [H^1(\Omega)]^d$.
We will indicate by $\bs{w}$ the kinematic (vector) test function in $\mathcal{V}$. 
The specific internal energy $e$ is discretized in
$\mathcal{E} \subset L_2(\Omega)$, with a corresponding (scalar) test function $\phi$.
Throughout this manuscript, we refer to the pairs of spaces
$(Q_m)^d - \hat{Q}_{m-1}$, by which $(Q_m)^d = \mathcal{V}$
is the Cartesian product
of the space of continuous finite elements on quadrilateral or hexahedral
meshes of degree $m$, and $\hat{Q}_{m-1} = \mathcal{E}$ is the companion space
of discontinuous finite elements of order one less than the kinematic space.
The semi-discrete system is:
\begin{subequations}
\begin{align}
\label{eq_weak_form}
\frac{d \boldsymbol{x}}{dt} 
&=\;
\boldsymbol{v} \; , 
\\
\int_{\Omega(t)}  \bs{w} \cdot \left( \rho \frac{d \bs{v}}{dt} \right)
&=\;
\int_{\Omega(t)} \nabla \bs{w} : \left( p \bs{I} - \bs{\sigma}_a \right) \; , 
\\  
\int_{\Omega(t)} \phi \, \rho \frac{d e}{dt} 
&=\;
\int_{\Omega(t)}  \phi \, \left( -p \nabla \cdot \bs{v} + \bs{\sigma}_a : \nabla \bs{v} \right) \; ,
\end{align}
\end{subequations}
where $p$ is the thermodynamic pressure, computed from an equation of state
of the type $p = \mathcal{P}(\rho,{e})$, and $\bs{\sigma}_a$ is the artificial viscosity stress tensor, 
which provides numerical stability in the presence of shocks, see Section 6 in \cite{Dobrev2012}.
Note also that $d/dt$ indicates the total material derivative (or Lagrangian derivative), i.e., for some field $\vartheta$,
we have $d \vartheta / dt = \partial \vartheta / \partial t + \bs{v} \cdot \nabla \vartheta$.

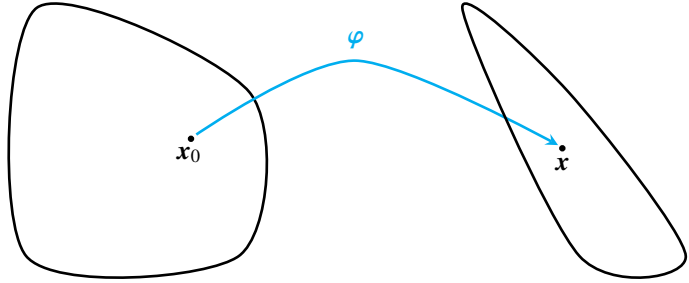
\begin{figure}[h!]
\centering
  \begin{tikzpicture}[line width=1.0, scale = 1.4]
    \draw [black] plot [smooth cycle] coordinates{(-2.5, -3.7) (-.5, -4.5) (-0.6, -6.05)  (-2.6, -6.05)};
    \draw [cyan,-stealth] plot [smooth] coordinates{(-1.0, -4.9) (0.5, -4.2) (2.4, -5.0)};
    \draw [black] plot [smooth cycle] coordinates{(1.5, -3.7) (2.5, -4.5) (3.6, -6.05)  (2.6, -6.05)};
    \filldraw [black] (-1.05, -4.94) circle (0.5pt);
    \draw [black] (-1.07, -5.09) node{$\bs{x}_0$};
    \draw [cyan] (0.5, -4.0) node{$\bs{\varphi}$};
    \filldraw [black] (2.44, -5.03) circle (0.5pt);
    \draw [black] (2.44, -5.18) node{$\bs{x}$};
	\end{tikzpicture}
	\captionsetup{font=normalsize}
  \caption{Sketch of the Lagrangian map $\bs{\varphi}$.}
\label{fig_lagmap}
\end{figure}
	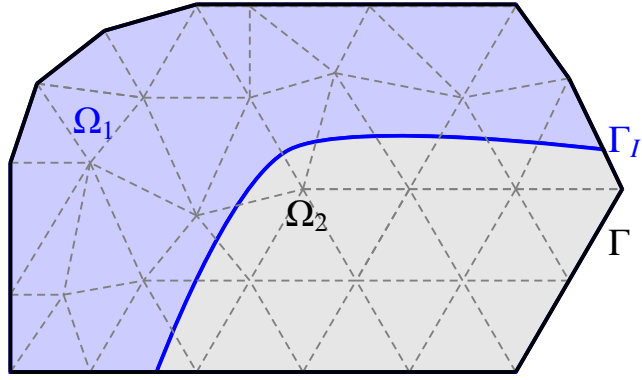
\begin{figure}[h!]
		\centering
	\begin{tikzpicture}[scale=0.7]
	\draw [black, draw=none,name path=bdr] plot coordinates {  (-0.75,-3.4641) (-3.5,-3.4641)  (-3.5,-2)  (-3.5,0.5) (-3,2) (-1.73205,3) (0,3.5) (6,3.5) (7,2.1) (7.7,0.75) };
	\draw [line width = 0.5mm,blue, fill=blue!20!, name path=true] plot[smooth] coordinates {(-0.75,-3.4641)  (1.75,0.75) (7.65,0.75) } -- (7,2.1)-- (6,3.5)  -- (0,3.5) -- (-1.73205,3)  -- (-3,2) -- (-3.5,0.5) -- (-3.5,-2) -- (-3.5,-3.4641)  --cycle    ;
	\draw [line width = 0.5mm,blue, fill=gray!20!, name path=true] plot[smooth] coordinates {(-0.75,-3.4641)  (1.75,0.75) (7.65,0.75) } -- (8,0)-- (6,-3.4641)  -- (-0.75,-3.4641)   --cycle;
	\draw[line width = 0.25mm,densely dashed,gray] (-1,1.73205) -- (0,3.4641);
	\draw[line width = 0.25mm,densely dashed,gray] (0,3.4641) -- (1,1.73205);
	\draw[line width = 0.25mm,densely dashed,gray] (1,1.73205) -- (1,3.4641);
	\draw[line width = 0.25mm,densely dashed,gray] (1,3.4641) -- (0,3.4641);
	\draw[line width = 0.25mm,densely dashed,gray] (1,3.4641) -- (2.6,2.2);
	\draw[line width = 0.25mm,densely dashed,gray] (1,3.4641) -- (3.25,3.4641);
	\draw[line width = 0.25mm,densely dashed,gray] (3.25,3.4641) -- (2.6,2.2);
	\draw[line width = 0.25mm,densely dashed,gray] (3.25,3.4641) -- (5,1.73205);
	\draw[line width = 0.25mm,densely dashed,gray] (3.25,3.4641) -- (6,3.4641);
	\draw[line width = 0.25mm,densely dashed,gray] (6,3.4641) -- (5,1.73205);
	\draw[line width = 0.25mm,densely dashed,gray] (6,3.4641) -- (7,2.1);
	\draw[line width = 0.25mm,densely dashed,gray] (0,-0.5) -- (-2,0.5);
	\draw[line width = 0.25mm,densely dashed,gray] (-2,0.5) -- (-1,1.73205);
	\draw[line width = 0.25mm,densely dashed,gray] (-2,0.5) -- (-1,-1.73205);
	\draw[line width = 0.25mm,densely dashed,gray] (-2,0.5) -- (-2.5,-2);
	\draw[line width = 0.25mm,densely dashed,gray] (-2.5,-2) -- (-1,-1.73205);
	\draw[line width = 0.25mm,densely dashed,gray] (-2.5,-2) -- (-2,-3.4641);
	\draw[line width = 0.25mm,densely dashed,gray] (0,-0.5) -- (-1,1.73205);
	\draw[line width = 0.25mm,densely dashed,gray] (-1,1.73205) -- (1,1.73205);
	\draw[line width = 0.25mm,densely dashed,gray] (0,-0.5) -- (2,0);
	\draw[line width = 0.25mm,densely dashed,gray] (2,0) -- (1,1.73205);
	\draw[line width = 0.25mm,densely dashed,gray] (1,1.73205) -- (0,-0.5);
	\draw[line width = 0.25mm,densely dashed,gray] (2,0) -- (2.6,2.2);
	\draw[line width = 0.25mm,densely dashed,gray] (2.6,2.2) -- (1,1.73205);
	\draw[line width = 0.25mm,densely dashed,gray] (2,0) -- (4,0);
	\draw[line width = 0.25mm,densely dashed,gray] (4,0) -- (2.6,2.2);
	\draw[line width = 0.25mm,densely dashed,gray] (2.6,2.2) -- (5,1.73205);
	\draw[line width = 0.25mm,densely dashed,gray] (5,1.73205) -- (4,0);
	\draw[line width = 0.25mm,densely dashed,gray] (4,0) -- (6,0);
	\draw[line width = 0.25mm,densely dashed,gray] (6,0) -- (5,1.73205);
	\draw[line width = 0.25mm,densely dashed,gray] (6,0) -- (7,2.1);
	\draw[line width = 0.25mm,densely dashed,gray] (7,2.1) -- (5,1.73205);
	\draw[line width = 0.25mm,densely dashed,gray] (6,0) -- (8,0);
	\draw[line width = 0.25mm,densely dashed,gray] (8,0) -- (7,2.1);
	\draw[line width = 0.25mm,densely dashed,gray] (0,-0.5) -- (-1,-1.73205);
	\draw[line width = 0.25mm,densely dashed, gray] (-1,-1.73205) -- (1,-1.73205);
	\draw[line width = 0.25mm,densely dashed,gray] (2,0) -- (1,-1.73205);
	\draw[line width = 0.25mm,densely dashed,gray] (1,-1.73205) -- (0,-0.5);
	\draw[line width = 0.25mm,densely dashed,gray] (2,0) -- (3,-1.73205);
	\draw[line width = 0.25mm,densely dashed,gray] (3,-1.73205) -- (1,-1.73205);
	\draw[line width = 0.25mm,densely dashed,gray] (4,0) -- (3,-1.73205);
	\draw[line width = 0.25mm,densely dashed,gray] (2,0) -- (4,0);
	\draw[line width = 0.25mm,densely dashed,gray] (4,0) -- (3,-1.73205);
	\draw[line width = 0.25mm,densely dashed,gray] (3,-1.73205) -- (5,-1.73205);
	\draw[line width = 0.25mm,densely dashed,gray] (5,-1.73205) -- (4,0);
	\draw[line width = 0.25mm,densely dashed,gray] (4,0) -- (6,0);
	\draw[line width = 0.25mm,densely dashed,gray] (6,0) -- (5,-1.73205);
	\draw[line width = 0.25mm,densely dashed,gray] (6,0) -- (7,-1.73205);
	\draw[line width = 0.25mm,densely dashed,gray] (7,-1.73205) -- (5,-1.73205);
	\draw[line width = 0.25mm,densely dashed,gray] (6,0) -- (8,0);
	\draw[line width = 0.25mm,densely dashed,gray] (8,0) -- (7,-1.73205);
	\draw[line width = 0.25mm,densely dashed,gray] (-1.73205,3) -- (0,3.5);
	\draw[line width = 0.25mm,densely dashed,gray] (-3,2) -- (-1.73205,3);
	\draw[line width = 0.25mm,densely dashed,gray] (-1,1.73205) -- (-1.73205,3);
	\draw[line width = 0.25mm,densely dashed,gray] (-3,2) -- (-1,1.73205);
	\draw[line width = 0.25mm,densely dashed,gray] (-3,2) -- (-1,1.73205);
	\draw[line width = 0.25mm,densely dashed,gray] (-3,2) -- (-3.5,0.5);
	\draw[line width = 0.25mm,densely dashed,gray] (-3,2) -- (-2,0.5);
	\draw[line width = 0.25mm,densely dashed,gray] (-3.5,0.5) -- (-2,0.5);
	\draw[line width = 0.25mm,densely dashed,gray] (-3.5,-2) -- (-3.5,0.5);
	\draw[line width = 0.25mm,densely dashed,gray] (-3.5,-2) -- (-2,0.5);
	\draw[line width = 0.25mm,densely dashed,gray] (-3.5,-2) -- (-3.5,-3.4641);
	\draw[line width = 0.25mm,densely dashed,gray] (-3.5,-2) -- (-2.5,-2);
	\draw[line width = 0.25mm,densely dashed,gray] (-3.5,-3.4641) -- (-2.5,-2);
	\draw[line width = 0.25mm,densely dashed,gray] (-2,-3.4641) -- (-3.5,-3.4641);
	\draw[line width = 0.25mm,densely dashed,gray] (0,-3.4641) -- (-2,-3.4641);
	\draw[line width = 0.25mm,densely dashed,gray] (-2,-3.4641) -- (-1,-1.73205);
	\draw[line width = 0.25mm,densely dashed,gray]  (-1,-1.73205) -- (0,-3.4641);
	\draw[line width = 0.25mm,densely dashed,gray] (0,-3.4641) -- (1,-1.73205);
	\draw[line width = 0.25mm,densely dashed,gray] (0,-3.4641) -- (2,-3.4641);
	\draw[line width = 0.25mm,densely dashed,gray] (2,-3.4641) -- (1,-1.73205);
	\draw[line width = 0.25mm,densely dashed,gray] (2,-3.4641) -- (3,-1.73205);
	\draw[line width = 0.25mm,densely dashed,gray] (2,-3.4641) -- (4,-3.4641);
	\draw[line width = 0.25mm,densely dashed,gray] (4,-3.4641) -- (3,-1.73205);
	\draw[line width = 0.25mm,densely dashed,gray] (4,-3.4641) -- (5,-1.73205);
	\draw[line width = 0.25mm,densely dashed,gray] (4,-3.4641) -- (6,-3.4641);
	\draw[line width = 0.25mm,densely dashed,gray] (6,-3.4641) -- (5,-1.73205);
	\draw[line width = 0.25mm,densely dashed,gray] (6,-3.4641) -- (7,-1.73205);
	\draw[line width = 0.5mm,black]  (-3.5,-3.4641) -- (6,-3.4641);
	\draw[line width = 0.5mm,black]  (-3.5,-3.4641) -- (-3.5,0.5);
	\draw[line width = 0.5mm,black]  (-3.5,0.5) -- (-3,2);
	\draw[line width = 0.5mm,black]  (-3,2) -- (-1.73205,3);
	\draw[line width = 0.5mm,black]  (-1.73205,3) -- (0,3.5);
	\draw[line width = 0.5mm,black]  (0,3.5) -- (6,3.5);
	\draw[line width = 0.5mm,black]  (6,3.5) -- (7,2.1) ;
	\draw[line width = 0.5mm,black]  (7,2.1) -- (8,0);

	\draw[line width = 0.5mm,black]  (6,-3.4641) --  (8,0) ;
	\node[text width=5cm] at (1.25,1.25) {\Large${\color{blue}\Omega_{1}}$};
	\node[text width=0.5cm] at (8.1,0.9) {\Large${\color{blue}\G_{I}}$};
	\node[text width=5cm] at (5.25,-0.5) {\Large${\color{black}\Omega_{2}}$};
	\node[text width=0.5cm] at (8.1,-1.0) {\Large${\color{black}\G}$};
	\end{tikzpicture}
	\caption{The true domain $\Om$, with boundary $\G=\partial \Om$, and the internal interface $\G_I$, which decomposes $\Om$ as $\Om = \Omega_{1} \cup \Omega_{2} \cup \G_{I}$. 
	}
	\label{domain}
\end{figure}
\paragraph{Interface conditions}
When multiple materials are present in the domain $\Om$, the principles of
mass and momentum conservation imply specific conditions at
the interfaces separating the materials.
For the sake of simplicity and without loss of generality, consider the case of two materials and an interface separating them.
Let $\Om$ be the physical domain of the problem, a bounded and connected open region in $\mathbb{R}^{n}$ ($n=2$ or $3$).
Let $\Gamma :=\partial\Om$ be the Lipschitz boundary of $\Om$, composed as the union of $C^2$ curves (in two dimensions)
or surfaces (in three dimensions) that intersect at a finite number of vertices or edges. 
We assume that $\Om$ is partitioned by an interface $\G_I$ into two non-overlapping
domains satisfying $\Omega =  \Omega_{1} \cup \Omega_{2} \cup \G_{I}$, as shown in Figure \ref{domain}.
Here $\Omega_{1} \cup \Omega_{2} \neq \emptyset$ and
$\G_{I} = \text{clos}(\Omega_{1}) \cap \text{clos}(\Omega_{2})$
defines the internal interface (also assumed Lipschitz)
between $ \Omega_{1}$ and  $\Omega_{2}$.
Then the system of equations~\eqref{eq_weak_form} would hold on
$\Omega = \Omega_{1} \cup \Omega_{2}$ with the following
interface conditions on $\G_{I}$:
\begin{subequations}
\label{eq_interface_cond}
\begin{align}
\jump{p} &= 0 \; , 
\label{eq:SH_StrongInterfaceD}
\\
\jump{\bs{v}} \cdot \bs{n}_{1}&= 0 
\; .
\label{eq:SH_StrongInterfaceN}
\end{align}
\end{subequations}
Note that $\jump{p} = (p_{1}-p_{2}) \bs{n}_{1}$,
$\jump{\bs{v}} = \bs{v}^{+}-\bs{v}^{-}$
indicate jumps in pressure and velocity, with $\bs{n}_{1}$ and $\bs{n}_{2}$
the outward-pointing normals to $\partial \Om_{1}$ and $\partial \Om_{2}$,
respectively.
It is also useful to define the average
$\avgb{\zeta}_{\gamma} = \gamma \, \zeta_{1} + (1-\gamma)\, \zeta_{2}$
for $\gamma \, \in \, \{0,1\}$.

\begin{remark}
\label{rem_p_cond}
We choose to enforce the pressure condition \eqref{eq:SH_StrongInterfaceD}
as this is common in most closure model approaches
\cite{Yanilkin13, Barlow2014, Tomov2016, Waltz2021}.
Applications often require a more sophisticated condition on the normal component of the full stresss tensor, including the artificial viscosity,
namely,
\[
\jump{(-p \bs{I} + \bs{\sigma}_a) \bs{n}_1} = 0 
\; .
\]
In principle the shifted approach can be used to enforce such condition,
and we will consider this in the future.
The present work focuses on the simpler pressure condition
\eqref{eq:SH_StrongInterfaceD}.
\end{remark}

\subsection{Integration and infinite dimensional function spaces}
In the sequel, $(v,w)_\omega = \int_\omega v \,w$ denotes the $L^2$ inner product on a subset $\omega \subset \Omega$ and $\avg{v,w}_{\zeta} = \int_{\zeta} v \, w$ denotes the $L^2$ inner product on a subset $\zeta$ of one of the boundaries or interfaces present in $\Omega$.
Throughout the paper, we will use the Sobolev spaces $H^m(\Om)=W^{m,2}(\Om)$ of index of regularity $m \geq 0$ and index of summability 2, equipped with the (scaled) norm
\begin{equation}
\|v \|_{H^{m}(\Om)} 
= \left( \| v \|^2_{L^2(\Om)} + \sum_{k = 1}^{m} \| l(\Om)^k  \bs{D}^k v \|^2_{L^2(\Om)} \right)^{1/2} \; ,
\end{equation}
where $\bs{D}^{k}$ is the $k$th-order spatial derivative operator and $l(A)= (\mathrm{meas}(A))^{1/{n_d}}$ is a characteristic length of the domain $A$. Note that $H^0(\Om)=L^{2}(\Om)$.  As usual, we use a simplified notation for norms and semi-norms, i.e., we set $\| v \|_{m;\Om}=\| v \|_{H^m(\Om)}$ and $| v |_{k;\Om}= 
\| \bs{D}^k v \|_{0;\Om}= \| \bs{D}^k v \|_{L^2(\Om)}$.

%% file: SIM_general.tex
\section{A general introduction to the Shifted Interface Method}
\label{sec:gen}

\begin{figure}
	\begin{subfigure}[h]{.45\textwidth}\centering
		\begin{tikzpicture}[scale=0.7]
		\draw [black, draw=none,name path=surr] plot coordinates { (-2,-3.4641) (-1,-1.73205) (0,-0.5) (1,1.73205) (2.6,2.2) (5,1.73205) (7,2.1) (7.62,0.8) };
		\draw [black, draw=none,name path=bdr] plot coordinates {  (-3.5,-3.4641)  (-3.5,-2)  (-3.5,0.5) (-3,2) (-1.73205,3) (0,3.5) (6,3.5) (7,2.1)};
		\draw [blue, draw=none, name path=inter] plot coordinates {(-2,-3.4641) (0,-3.4641) (1,-1.73205) (2,0) (4,0) (6,0) (8,0) };
		\draw [black, draw=none,name path=outersurr] plot coordinates {  (8,0) (2,0) (0,-3.4641) };
		\draw [black, draw=none,name path=outerbdr] plot coordinates {  (8,0) (6,-3.4641) (0,-3.4641) };
		\tikzfillbetween[of=bdr and surr,split]{red!5!};
		\tikzfillbetween[of=surr and inter, split]{olive!5!};
		\tikzfillbetween[of=outersurr and outerbdr, split]{gray!10!};
		\draw[line width = 0.25mm,densely dashed,gray] (-1,1.73205) -- (0,3.4641);
		\draw[line width = 0.25mm,densely dashed,gray] (0,3.4641) -- (1,1.73205);
		\draw[line width = 0.25mm,densely dashed,gray] (1,1.73205) -- (1,3.4641);
		\draw[line width = 0.25mm,densely dashed,gray] (1,3.4641) -- (0,3.4641);
		\draw[line width = 0.25mm,densely dashed,gray] (1,3.4641) -- (2.6,2.2);
		\draw[line width = 0.25mm,densely dashed,gray] (1,3.4641) -- (3.25,3.4641);
		\draw[line width = 0.25mm,densely dashed,gray] (3.25,3.4641) -- (2.6,2.2);
		\draw[line width = 0.25mm,densely dashed,gray] (3.25,3.4641) -- (5,1.73205);
		\draw[line width = 0.25mm,densely dashed,gray] (3.25,3.4641) -- (6,3.4641);
		\draw[line width = 0.25mm,densely dashed,gray] (6,3.4641) -- (5,1.73205);
		\draw[line width = 0.25mm,densely dashed,gray] (6,3.4641) -- (7,2.1);
		\draw[line width = 0.25mm,densely dashed,gray] (0,-0.5) -- (-2,0.5);
		\draw[line width = 0.25mm,densely dashed,gray] (-2,0.5) -- (-1,1.73205);
		\draw[line width = 0.25mm,densely dashed,gray] (-2,0.5) -- (-1,-1.73205);
		\draw[line width = 0.25mm,densely dashed,gray] (-2,0.5) -- (-2.5,-2);
		\draw[line width = 0.25mm,densely dashed,gray] (-2.5,-2) -- (-1,-1.73205);
		\draw[line width = 0.25mm,densely dashed,gray] (-2.5,-2) -- (-2,-3.4641);
		\draw[line width = 0.25mm,densely dashed,gray] (0,-0.5) -- (-1,1.73205);
		\draw[line width = 0.25mm,densely dashed,gray] (-1,1.73205) -- (1,1.73205);
		\draw[line width = 0.25mm,orange] (0,-0.5) -- (2,0);
		\draw[line width = 0.25mm,orange] (2,0) -- (1,1.73205);
		\draw[line width = 0.25mm,orange] (1,1.73205) -- (0,-0.5);
		\draw[line width = 0.25mm,orange] (2,0) -- (2.6,2.2);
		\draw[line width = 0.25mm,densely dashed,gray] (2.6,2.2) -- (1,1.73205);
		\draw[line width = 0.25mm,densely dashed,gray] (2,0) -- (4,0);
		\draw[line width = 0.25mm,orange] (4,0) -- (2.6,2.2);
		\draw[line width = 0.25mm,densely dashed,gray] (2.6,2.2) -- (5,1.73205);
		\draw[line width = 0.25mm,orange] (5,1.73205) -- (4,0);
		\draw[line width = 0.25mm,densely dashed,gray] (4,0) -- (6,0);
		\draw[line width = 0.25mm,orange] (6,0) -- (5,1.73205);
		\draw[line width = 0.25mm,orange] (6,0) -- (7,2.1);
		\draw[line width = 0.25mm,densely dashed,gray] (7,2.1) -- (5,1.73205);
		\draw[line width = 0.25mm,densely dashed,gray] (6,0) -- (8,0);
		\draw[line width = 0.5mm,cyan] (8,0) -- (7,2.1);
		\draw[line width = 0.25mm,densely dashed,gray] (0,-0.5) -- (-1,-1.73205);
		\draw[line width = 0.25mm,densely dashed, orange] (-1,-1.73205) -- (1,-1.73205);
		\draw[line width = 0.25mm,densely dashed,gray] (2,0) -- (1,-1.73205);
		\draw[line width = 0.25mm,orange] (1,-1.73205) -- (0,-0.5);
		\draw[line width = 0.25mm,densely dashed,gray] (2,0) -- (3,-1.73205);
		\draw[line width = 0.25mm,densely dashed,gray] (3,-1.73205) -- (1,-1.73205);
		\draw[line width = 0.25mm,densely dashed,gray] (4,0) -- (3,-1.73205);
		\draw[line width = 0.25mm,densely dashed,gray] (2,0) -- (4,0);
		\draw[line width = 0.25mm,densely dashed,gray] (4,0) -- (3,-1.73205);
		\draw[line width = 0.25mm,densely dashed,gray] (3,-1.73205) -- (5,-1.73205);
		\draw[line width = 0.25mm,densely dashed,gray] (5,-1.73205) -- (4,0);
		\draw[line width = 0.25mm,densely dashed,gray] (4,0) -- (6,0);
		\draw[line width = 0.25mm,densely dashed,gray] (6,0) -- (5,-1.73205);
		\draw[line width = 0.25mm,densely dashed,gray] (6,0) -- (7,-1.73205);
		\draw[line width = 0.25mm,densely dashed,gray] (7,-1.73205) -- (5,-1.73205);
		\draw[line width = 0.25mm,densely dashed,gray] (6,0) -- (8,0);
		\draw[line width = 0.25mm,densely dashed,gray] (8,0) -- (7,-1.73205);
		\draw[line width = 0.25mm,densely dashed,gray] (-1.73205,3) -- (0,3.5);
		\draw[line width = 0.25mm,densely dashed,gray] (-3,2) -- (-1.73205,3);
		\draw[line width = 0.25mm,densely dashed,gray] (-1,1.73205) -- (-1.73205,3);
		\draw[line width = 0.25mm,densely dashed,gray] (-3,2) -- (-1,1.73205);
		\draw[line width = 0.25mm,densely dashed,gray] (-3,2) -- (-1,1.73205);
		\draw[line width = 0.25mm,densely dashed,gray] (-3,2) -- (-3.5,0.5);
		\draw[line width = 0.25mm,densely dashed,gray] (-3,2) -- (-2,0.5);
		\draw[line width = 0.25mm,densely dashed,gray] (-3.5,0.5) -- (-2,0.5);
		\draw[line width = 0.25mm,densely dashed,gray] (-3.5,-2) -- (-3.5,0.5);
		\draw[line width = 0.25mm,densely dashed,gray] (-3.5,-2) -- (-2,0.5);
		\draw[line width = 0.25mm,densely dashed,gray] (-3.5,-2) -- (-3.5,-3.4641);
		\draw[line width = 0.25mm,densely dashed,gray] (-3.5,-2) -- (-2.5,-2);
		\draw[line width = 0.25mm,densely dashed,gray] (-3.5,-3.4641) -- (-2.5,-2);
		\draw[line width = 0.25mm,densely dashed,gray] (-2,-3.4641) -- (-3.5,-3.4641);
		\draw[line width = 0.5mm,cyan] (0,-3.4641) -- (-2,-3.4641);
		\draw[line width = 0.25mm,orange] (-2,-3.4641) -- (-1,-1.73205);
		\draw[line width = 0.25mm,orange]  (-1,-1.73205) -- (0,-3.4641);
		\draw[line width = 0.25mm,densely dashed,gray] (0,-3.4641) -- (1,-1.73205);
		\draw[line width = 0.25mm,densely dashed,gray] (0,-3.4641) -- (2,-3.4641);
		\draw[line width = 0.25mm,densely dashed,gray] (2,-3.4641) -- (1,-1.73205);
		\draw[line width = 0.25mm,densely dashed,gray] (2,-3.4641) -- (3,-1.73205);
		\draw[line width = 0.25mm,densely dashed,gray] (2,-3.4641) -- (4,-3.4641);
		\draw[line width = 0.25mm,densely dashed,gray] (4,-3.4641) -- (3,-1.73205);
		\draw[line width = 0.25mm,densely dashed,gray] (4,-3.4641) -- (5,-1.73205);
		\draw[line width = 0.25mm,densely dashed,gray] (4,-3.4641) -- (6,-3.4641);
		\draw[line width = 0.25mm,densely dashed,gray] (6,-3.4641) -- (5,-1.73205);	
		\draw[line width = 0.25mm,densely dashed,gray] (6,-3.4641) -- (7,-1.73205);	
		\draw [line width = 0.5mm,blue, name path=true] plot[smooth] coordinates {(-0.75,-3.5) (1.75,0.75) (7.7,0.75)};
		\draw[line width = 0.5mm,red] (1,1.73205) -- (2.6,2.2);
		\draw[line width = 0.5mm,red] (2.6,2.2) -- (5,1.73205);
		\draw[line width = 0.5mm,red] (5,1.73205) --  (7,2.1);
		\draw[line width = 0.5mm,red] (1,1.73205) -- (0,-0.5);
		\draw[line width = 0.5mm,red] (0,-0.5) -- (-1,-1.73205);
		\draw[line width = 0.5mm,red] (-1,-1.73205) -- (-2,-3.4641);
		\draw[line width = 0.5mm,red] (-2,-3.4641) -- (-3.5,-3.4641);
		\draw[line width = 0.5mm,red] (-3.5,-3.4641) -- (-3.5,-2);
		\draw[line width = 0.5mm,red] (-3.5,-2) --  (-3.5,0.5) ;
		\draw[line width = 0.5mm,red] (-3.5,0.5)-- (-3,2);
		\draw[line width = 0.5mm,red] (-1.73205,3)-- (-3,2);
		\draw[line width = 0.5mm,red] (-1.73205,3)-- (0,3.5);
		\draw[line width = 0.5mm,red] (6,3.5)-- (0,3.5);
		\draw[line width = 0.5mm,red] (6,3.5)-- (7,2.1);
		\draw[line width = 0.5mm,olive]  (0,-3.4641) -- (1,-1.73205);
		\draw[line width = 0.5mm,olive]  (1,-1.73205) --  (2,0) ;
		\draw[line width = 0.5mm,olive] (2,0)-- (8,0);
		\draw[line width = 0.5mm,olive]  (0,-3.4641) -- (6,-3.4641);
		\draw[line width = 0.5mm,olive]  (6,-3.4641) --  (8,0) ;
		\node[text width=5cm] at (7.15,2.5) {\large${\color{red}\ti{\G}_{1}}$};
		\node[text width=5cm] at (7.15,3.85) {\large${\color{red}\G_{1}}$};
		\node[text width=3cm] at (1.85,1.25) {\large${\color{red}\ti{\Om}_{1}}$};
		\node[text width=0.5cm] at (1.25,0.8) {\large${\color{blue}\G_{I}}$};
		\node[text width=0.5cm] at (3,-0.5) {\large${\color{olive}{\ti{\G}}_{2}}$};
		\node[text width=0.5cm] at (3,-3.95) {\large${\color{olive}{\G}_{2}}$};
		\node[text width=0.5cm] at (4,1.5) {\large${\color{orange}{\Om}_{cut}}$};
		\node[text width=0.5cm] at (6.7,0.4) {\large${\color{orange}\mathcal{E}^o_{c}}$};
		\node[text width=0.5cm] at (8,1.4) {\large${\color{cyan}\G_{c}}$};
		\node[text width=0.5cm] at (-1,-3.95) {\large${\color{cyan}\G_{c}}$};
		\node[text width=5cm] at (5.25,-2.25) {\large${\color{olive}\ti{\Om}_{2}}$};
		\end{tikzpicture}
		\caption{The surrogate domains $\ti{\Om}_{1}$ and $\ti{\Om}_{2}$, the surrogate interfaces $\ti{\G}_{1}$ and $\ti{\G}_{2}$, the set of cut elements $\Om_{cut}$, and cut boundary edges $\G_{c}$.}
		\label{SBM}
	\end{subfigure}
	\hspace{1.5cm}
	\begin{subfigure}[h]{.45\textwidth}
		\centering
		\begin{tikzpicture}[scale=0.85]
		\draw[line width = 0.25mm,densely dashed,gray] (0,0.5) -- (-1.5,3);
		\draw[line width = 0.25mm,densely dashed,gray] (-1.5,3) -- (0.5,5);
		\draw[line width = 0.25mm,densely dashed,gray] (0,0.5) -- (2.5,2);
		\draw[line width = 0.25mm,densely dashed,gray] (2.5,2) -- (0.5,5);
		\draw[line width = 0.25mm,densely dashed,gray] (0.5,5) -- (0,0.5);
		\draw [line width = 0.5mm,blue, name path=true] plot[smooth] coordinates {(1,-0.5) (2.25,2.5) (0.75,6)};
		\draw[line width = 0.5mm,red] (0,0.5) -- (0.5,5);
		\node[text width=0.5cm] at (0.5,5.5) {\large${\color{red}\tG}$};
		\node[text width=0.5cm] at (1.75,5.5) {\large${\color{blue}\G}$};
		\node[text width=0.5cm] at (1.25,3.25) {\large$\bs{d}$};
		\node[text width=0.5cm] at (3,3.5) {\large$\bs{n}$};
		\node[text width=0.5cm] at (2.7,2.25) {\large$\bs{\tau}$};
		\draw[->,line width = 0.25mm,-latex] (0.25,2.75) -- (2.12,3.1);
		\draw[->,line width = 0.25mm,-latex] (2.12,3.1) -- (2.40,2.1);
		\draw[->,line width = 0.25mm,-latex] (2.12,3.1) -- (2.95,3.25);
		\end{tikzpicture}
		\caption{The distance vector $\bs{d}$, the true normal $\bs{n}$ and the true tangent $\bs{\tau}$.}
		\label{ntd}
	\end{subfigure}
	\caption{The surrogate domain, its boundary, and the distance vector $\bs{d}$.}
\end{figure}
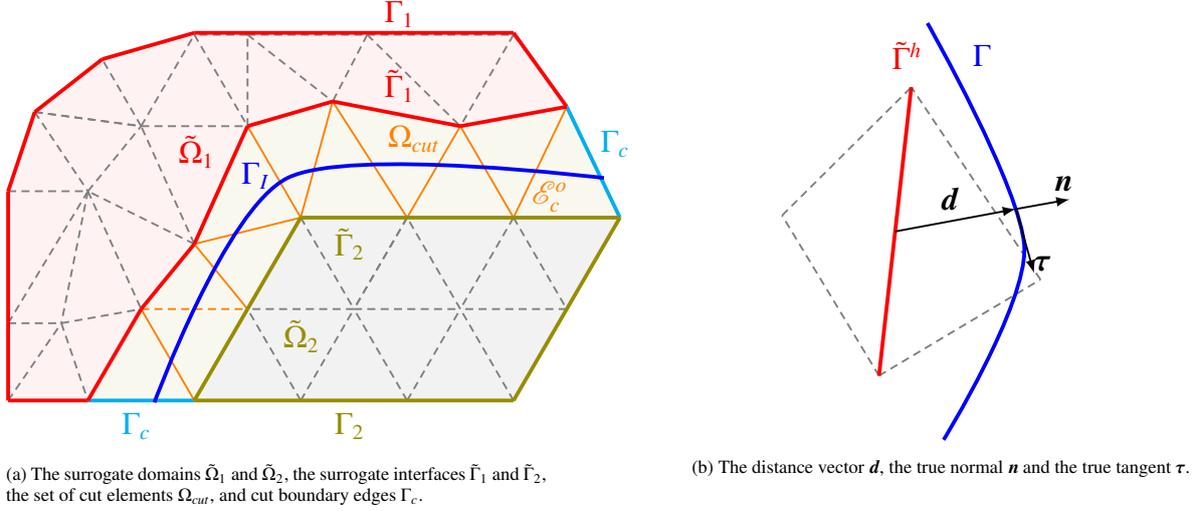
\subsection{The triangulation $\mathcal{T}^h$ of $\Om$}
Let $\mathcal{T}^h$ be a family of admissible and shape-regular triangulations of $\Om$. 
In what follows, we will make the specific assumption that the triangulation $\mathcal{T}^h$ is {\it fitted} to the boundary $\G$ of $\Om$.
The set of all edges (in two dimensions) or faces (in three dimensions) of the triangulation $\mathcal{T}^h$ is indicated as $\mathcal{E}$.
We indicate by $h_T$ the diameter of an element $T \in \mathcal{T}^h$ and by $h$ the piecewise constant function in $\Om$ such that $h_{|T}=h_T$ for all $T \in \mathcal{T}^h$.

\subsection{The true domain, the surrogate domain and maps}
As shown in Figure \ref{SBM}, consider a grid, fitted to the boundary $\G = \partial \Om$, but not to $\G_I$.
We can then define the triangulations $\ti{\mathcal{T}}_{1}^h$ and $\ti{\mathcal{T}}_{2}^h$ by selecting those elements that are {\it{strictly}} contained in the closure of $\Om_{1}$ and $\Om_{2}$ respectively, i.e.,
\begin{align}
\ti{\mathcal{T}}_{i}^h := \{ T \in \mathcal{T}^h : T \subset \text{clos}(\Om_i)  \}\,,
\end{align}
where $i = 1$ or $2$.
$\ti{\mathcal{T}}_{1}^h$ and $\ti{\mathcal{T}}_{2}^h$ then identify the  surrogate domains  $\ti{\Om}_{1}$ and $\ti{\Om}_{2}$ (the areas shaded in red and dark green in Figure \ref{SBM}).
Namely, for $i = 1$ or $2$,
\begin{align}
\ti{\Om}_{i} := \text{int} \left( \bigcup_{T \in \ti{\mathcal{T}}_{i}^h}  T \right) \subseteq \Om_i  \,.
\end{align}
Furthermore, the set of all edges (in two dimensions) or faces (in three dimensions) of the triangulation $\ti{\mathcal{T}}_{1}^h$ and $\ti{\mathcal{T}}_{2}^h$ is indicated as $\mathcal{E}_{1}$ and $\mathcal{E}_{2}$, respectively.
We also define the subset of elements that are cut by $\G_{I}$, namely
\begin{align}
\mathcal{T}_c^h := \{ T \in \mathcal{T}^h : T \cap \G_I \neq \emptyset \}\, ,
\end{align}
and the corresponding {\it cut domain}
\begin{align}
\Om_c := \text{int} \left( \bigcup_{T \in \mathcal{T}_c^h}  T \right) \, .
\end{align}
We denote by 
\begin{align}
\G_i:= \ti{\Om}_i \cap \G
\end{align}
the body-fitted faces and by 
\begin{align}
\ti{\G}_i:=  \Om_{c} \cap \ti{\Om}_{i}
\end{align}
the surrogate faces.
The set of all interior edges/faces of the triangulation $\ti{\mathcal{T}}^h$, $\ti{\mathcal{T}}_{1}^h$ and $\ti{\mathcal{T}}_{2}^h$ is indicated by $\mathcal{E}^{o} = \mathcal{E} \setminus \G$, $\mathcal{E}_{1}^{o} = \mathcal{E}_{1} \setminus \G_{1}$ and $\mathcal{E}_{2}^{o} = \mathcal{E}_{2} \setminus \G_{2}$, respectively. $\mathcal{E}_c^o= \{ E \in \mathcal{E} \setminus \G : E \cap \G_I \neq \emptyset \}$ indicates the union of {\it interior} edges/faces (in two/three dimensions, resp.) that are cut by $\G_I$ (marked in orange in Figure \ref{SBM}). Finally, we denote by 
\begin{align}
\G_c:= \partial \Om_c \setminus \left(\ti{\G}_1 \cup \ti{\G}_2 \right)
\end{align}
the subset of edges/faces on the boundary $\G$ of $\Om$ that belong to elements in $\Om_c$ (depicted in light blue, in Figure \ref{SBM}).

\subsection{Maps from the surrogate to the true interfaces}
\label{sec:maps}

We now define the map 
\begin{subequations}
\begin{align}
\ti{\bs{M}}_{i}:&\; \ti{\G}_{i} \to \G \; , \, \, \forall \, {i} \, \,  \in \{1,2\} \\
    &   \ti{\bs{x}}_{i} \mapsto \bs{x}   \; ,
\end{align}
\end{subequations}
which maps a point $\ti{\bs{x}} \in \ti{\G}_{i}$ on the surrogate free surface to a point $\bs{x} \in \G_I$ on the true free surface. More generally, $\ti{\bs{M}}_{i}$ is a map between a surrogate and true boundaries in an immersed/embedded/unfitted domain computation, and can be built using the {\it closest-point projection} of points in $\ti{\G}_{i}$ onto $\G_I$, as shown in Figure~\ref{ntd}. 
In particular, it will become very important to characterize the map $\ti{\bs{M}}_{i}$ through a distance vector function 
\begin{align}
\label{eq:Mmap}
\bs{d}_{\ti{\bs{M}}_{i}} (\ti{\bs{x}})
\, = \, 
\bs{x}-\ti{\bs{x}}_{i}
\, = \, 
[ \, \ti{\bs{M}}_{i}-\bs{I} \, ] (\ti{\bs{x}})
\; .
\end{align}
As shown in Figure \ref{ntd}, if the closest-point projection is used, the vector $\bs{d}$ is aligned with $\bs{n}$. This case is typical when $\G_I$ is a smooth surface, while we refer to~\cite{TheoreticalPoissonAtallahCanutoScovazzi2020} for the case when edges/corners are present or when multiple boundary conditions are enforced on a smooth surface.
Moreover, we will denote by $\ti{\bs{n}}_{i}$ the unit outward-pointing normal to the surrogate faces  $\ti{\G}_i$, to be distinguished from the outward-pointing normal $\bs{n}$ to the true interface $\G_I$.

Equivalently, we can find a map $ \bs{M}_c: E \in \mathcal{E}^o_c \to \G_F $ , defined by means of the distance $ \bs{d}_{\bs{M}_c} = \bs{x}-\bs{x}_{c}$, from a point on  $ \bs{x}_{c}\in \mathcal{E}_c^o$, to its closest point projection on the true interface, $ \bs{x}\in\G_I $.
For the sake of simplicity and whenever there is no chance of confusion, the ``tilde'' and the subscript $c$ will be omitted from the map symbols $ \ti{\bs{M}}$ and $\bs{M}_c$, and will simply write ``$\bs{M}$.'' 
Similarly, we will omit the subscripts of $\bs{d}_{\ti{\bs{M}}_{i}} $ and $ \bs{d}_{\bs{M}_c}$, and write ``$\bs{d}$.'' 
\begin{rem}
In this work, the distance vector $\bs{d}$ is defined by means of the zero level set of a distance function, which is always well defined if the interface $\G_I$ is Lipschitz.
This approach is presented in full details in Section~\ref{sec_distance}. 
We also point the reader to the recent analysis in~\cite{TheoreticalPoissonAtallahCanutoScovazzi2020}, where the SBM and associated boundary maps are implemented for general domains with corners and edges. 
Most importantly, the methods that we are about to introduce do not depend on how $\bs{M}$ is constructed.
\end{rem}

\subsection{General strategy for interface conditions}
\label{sec:genStrat}

While the governing equations are discretized over $\Om$, the interface
conditions ~\eqref{eq:SH_StrongInterfaceD}-~\eqref{eq:SH_StrongInterfaceN} will be imposed on $\ti{\G}_{i} \cup \mathcal{E}_{c}^{o}$ rather than $\G_{I}$.
As such, the challenge is to appropriately impose interface conditions on $\ti{\G}_{i} \cup \mathcal{E}_{c}^{o}$ that would mimic~\eqref{eq:SH_StrongInterfaceD}-~\eqref{eq:SH_StrongInterfaceN}.
To this end, we resort to the $m^{th}$-order Taylor expansion of the variables of interest, centered at $\ti{\G}_{i} \cup \mathcal{E}_{c}^{o}$.
This approach allows to {\it shift} the interface conditions from $\G_{I}$ to $\ti{\G}_{i} \cup \mathcal{E}_{c}^{o}$.

Starting with the interface condition involving the pressure, let
us assume that $p$ is sufficiently smooth so as to admit a $m$th-order Taylor expansion pointwise, and let us denote by $\mathcal{D}^{k}_{\bs{d}}$ the
$k$th-order directional derivative in the direction of $\bs{d}$, defined  as $\mathcal{D}^{k}_{\bs{d}} p = \displaystyle{\sum_{\bs{\alpha} \in \mathbb{N}^n, |\bs{\alpha}|=k} \frac{k!}{\bs{\alpha}!}   \frac{\partial^k p}{\partial \bs{x}^{\bs{\alpha}}} \bs{d}^{\bs{\alpha}}  }$.
Then, we can write
\begin{align}
\jump{p(\bs{x}) } = \jump{p(\ti{\bs{x}}+\bs{d}(\ti{\bs{x}}))}
=
\jump{p(\tx) +  \sum_{k = 1}^{m-1}  \frac{\mathcal{D}^{k}_{\bs{d}} \, p(\ti{\bs{x}})}{k!} + (R^{m}(p,\bs{d}))(\tx) } &= 0  \, ,
\label{eq:shifted_SH_InterfaceD}
\end{align}
where the remainder $R^{m}(p,\bs{d}) = o(\Vert \bs{d}\Vert^{m})$ as
$\Vert \bs{d}\Vert \to 0$.
Now we introduce the operator
\begin{equation}
\label{eq:def-bndS}
\tS^{m-1} p := p+ \sum_{k = 1}^{m-1} \frac{\mathcal{D}^{k}_{\bs{d}} \, p}{k!},
\end{equation}
Then the Taylor expansion can be used to enforce \eqref{eq:SH_StrongInterfaceD}
on $\ti{\G}_{i} \cup \mathcal{E}_{c}^{o}$ rather than $\G_{I}$, as 
\begin{equation}
\label{eq:trace-u}
\jump{\tS^{m-1} p  + R^{m-1}(p,\bs{d})} 
= 0,
\end{equation}
which can be equivalently written as
	\begin{equation}\label{eq:trace-u-rewrite}
	\jump{p} 
	= 
	- \jump{ \sum_{k = 1}^{m-1}  \frac{\mathcal{D}^{k}_{\bs{d}} \, p }{k!} + R^{m}(p,\bs{d})}.
	\end{equation}
Neglecting the remainder $R^{m}(p,\bs{d})$, we obtain the final expression
of the {\it shifted} interface condition on the pressure
\begin{equation}
\label{eq:Finalu-g}
  \jump{\tS^{m-1} p} 
  = 
  \left( p_{1}+ \sum_{k = 1}^{m-1}  \frac{\mathcal{D}^{k}_{\bs{d}} \, p_{1} }{k!} 
  - p_{2} - \sum_{k = 1}^{m-1}  \frac{\mathcal{D}^{k}_{\bs{d}} \, p_{2} }{k!} \right)   \bs{n}_{1}(\bs{M} (\tx))
   = 0  \, , \quad
  \mbox{on } \ti{\G}_{i} \cup \mathcal{E}_{c}^{o} \;,
\end{equation}
which will be weakly enforced on the discretization $p_h$ of $p$ that
will be introduced later.
Similarly, for a vector field $\bs{v}$, we enforce 
\begin{equation}\label{eq:bsu-g}
\jump{\bs{S}_{h}^{m} \bs{v} + \bs{R}^{m+1}(\bs{v},\bs{d}) } \cdot \bs{n}_{1}(\bs{M} (\tx)) = 0 \, , \quad  \mbox{on } \ti{\G}_{i} \cup \mathcal{E}_{c}^{o}
\end{equation}
where $\bs{S}_{h}^{m} \bs{v} := \displaystyle{\bs{v}+ \sum_{k = 1}^{m}  \frac{\mathcal{D}^{k}_{\bs{d}} \, \bs{v} }{k!}}$
and $\bs{R}^{m+1}(\bs{v},\bs{d})$ is the Taylor expansion remainder of
$\bs{v}$ on $\ti{\G}_{i} \cup \mathcal{E}_{c}^{o}$. Again, neglecting the the remainder $\bs{R}^{m+1}(\bs{v},\bs{d})$, we obtain the {\it shifted} vector boundary condition 
\begin{equation}
\label{eq:Finalbsu-g}
  \jump{\bs{S}_{h}^{m} \bs{v}}  \cdot \bs{n}_{1}(\bs{M} (\tx)) = 0\, , \quad
  \mbox{on } \ti{\G}_{i} \cup \mathcal{E}_{c}^{o} \; .
\end{equation}
which is equivalent to 
\begin{equation}
\label{eq:Finalbsu-g-equiv}
\jump{ \bs{v} } \cdot \bs{n}_{1}(\bs{M} (\tx)) = -\left(\jump{\sum_{k = 1}^{m}  \frac{\mathcal{D}^{k}_{\bs{d}} \, \bs{v} }{k!}} \right) \cdot \bs{n}_{1}(\bs{M} (\tx))  \, , \quad
\mbox{on } \ti{\G}_{i} \cup \mathcal{E}_{c}^{o} \; .
\end{equation}
In what follows, whenever it does not cause confusion, we will simply write $\bs{n}_{1}$ in place of $\bs{n}_{1}(\bs{M} (\tx))$.


\subsection{Interface representation and distance computation}
\label{sec_distance}

As we are interested in time-dependent problems, in which the material
interface moves with the fluid velocity, the interface location must
be represented by a dynamic discrete function.
We utilize the standard Level-Set (LS) approach \cite{Osher1994} where the interface
is described by the zero LS of a finite element function $\eta$.
Since the amount of deformation in time would be limited by the resolution of
the kinematic FE space, it is natural to choose $\eta$ in $Q_m$, i.e., the
scalar version of the continuous FE space $\mathcal{V} = (Q_m)^d$
of the position $\bs{x}$ and velocity $\bs{v}$.
Once initialized, the LS function $\eta$ is simply evolved by
\begin{equation}
\label{eq_eta}
\frac{d \eta}{dt} = 0 \; ,
\end{equation}
which guarantees that the interface does not move into new elements and
stays continuous throughout the Lagrangian simulation, while it follows the
velocity curvature.

All shifted integrals require knowledge of the pointwise
distance vector $\bs{d}(\bs{x})$ to the zero LS of $\eta$.
There are many methods to compute such distances in the literature and we
refer the reader to a recent summary given in \cite{Belyaev2015}.
In this work we make use of two distance computation methods.
The first one is a simple renormalization
procedure given by
\begin{equation}
\label{eq_dist_normal}
d(\bs{x}) = \frac{\eta(\bs{x})}
            {\sqrt{ \eta(\bs{x})^2 + |\nabla \eta(\bs{x})|^2}} \; ,
\end{equation}
where $d(\bs{x}) \in Q_m$ is the signed scalar distance
from $\bs{x}$ to the zero LS of $\eta$ (see Section 6 of \cite{Belyaev2015}).
This is a very fast pointwise computation which gives reasonable accuracy
only in the vicinity of the zero level set.
This outcome usually suffices for our purposes since all shifted computations
primarily focus on the first band of elements surrounding the zero level set.

The second method to compute a FE distance function $\bs{d}(\bs{x})$ utilizes
the p-Laplacian distance computation (see Section 7 in \cite{Belyaev2015}).
That is, the magnitude $d(\bs{x}) \in Q_m$ of the distance vector
$\bs{d} \in Q_m^d$ is the solution of:
\begin{equation}
\label{eq_dist_plap}
\begin{split}
& \nabla \cdot \left( |\nabla d|^{p-2} \nabla d \right) = -1
  \text{ in } \Omega \; , \\
& d = 0 \text{ on } \partial \Omega, \quad 2 \leq p < \infty \; .
\end{split}
\end{equation}
This formulation computes the signed distance with respect to the boundary of
any domain $\Omega$, but the same approach can be customized for an arbitrary
LS of a discrete function, by choosing appropriate finite element
basis functions (see Section 3 in \cite{Rvachev1982}).
This approach is more expensive computationally, but offers improved accuracy.

In principle the distance function $\bs{d}(\bs{x})$ must be
updated after every Lagrangian time step through
\eqref{eq_dist_normal} or \eqref{eq_dist_plap}, as the mesh deformation
at $t^{n+1}$ deteriorates the signed distance property of $\bs{d}^n$.
In the case of the more expensive calculation \eqref{eq_dist_plap}, the update
can be deferred by a few steps, assuming that the deterioration is not rapid.
An alternative approach is to evolve directly the signed distance FE function
$d(\bs{x})$ with an additional redistancing procedure, for example as in
\cite{Kuzmin2013, Kuzmin2019}.

%% file: SIM_laghos.tex
\section{The Weighted Shifted Interface Method}
\label{sec_sim}

The Weighted Shifted Interface Method (WSIM) is inspired by the Weighted Shifted Boundary Method (WSBM), 
which was proposed in~\cite{colomes2021weighted} for boundary value problems associated with the Navier-Stokes equations with free surfaces.
The key idea in both methods is to weight the variational formulation by the volume fraction of the material(s), as a way to enhance the mass conservation properties of the overall algorithm.
In what follows, we adapt the WSBM to the specific case of shock hydrodynamics interfaces, and in particular the high-order computational framework proposed in~\cite{Dobrev2012}.


\subsection{Finite Element spaces}
\label{sec:weak_discr_poisson}
For an element $T \in \cT_h $, let the volume fraction, $\alpha_{i,T}$, for material $i \, \in \{1,2\}$ be defined as
\begin{align}
\label{eq_alpha}
\alpha_{1, T}(t) = \frac{\int_{T(t)} H(\eta(\bs{x}, t))}
                      {{\mathrm{meas}(T(t))}} \, , \quad
\alpha_{2, T}(t) = 1 - \alpha_{1, T}(t) \, ,
\end{align}
where $H$ is the Heaviside step function.
Clearly $\alpha_{i,T} \in [0 \, , 1 ]$ and more precisely: $\alpha_{i,T} = 1$ if $T \subset \ti{\Om}_{i}$; $ 0 < \alpha_{i,T} < 1$ if $T \in \Om_{c}$.
We also construct the global function
$\alpha_{i}(\bs{x}, t) \in \ti{\Om}_{i} \cup \Om_{c}$
such that $\alpha_{i}|_{T} = \alpha_{T}$.
Hence, $\alpha_{i} \in \mathcal{P}^{0}(\ti{\mathcal{T}}_{i}^{h}\cup \mathcal{T}_c^h)$ where $\mathcal{P}^{0}(\ti{\mathcal{T}}_{i}^{h} \cup \mathcal{T}_c^h)$ is the set of piecewise constant functions over $\ti{\mathcal{T}}_{i}^{h} \cup \mathcal{T}_c^h$.
Also note that \eqref{eq_alpha} is time dependent due to the mesh motion.
The numerator of \eqref{eq_alpha} is computed approximately,
by taking a large number of quadrature points in $T$ that
are either inside or outside of $\Om_{i}$;
more sophisticated integration methods can also be used to obtain
better accuracy, for example \cite{Muller2013}.

We now define the following {\it{weighted}} velocity and pressure spaces:
\begin{subequations}
\begin{equation}
\begin{split}
\bs{V}^h_{\alpha_{i}}(\Om, \mathcal{T}^{h} ) := 
  \bigg\{ \bs{v}^h_{\alpha_{i}} \in L_2(\Om)  :
  \bs{v}^h_{\alpha_{i}} = \alpha_{i} \, \bs{v}^{h}, ~
  \bs{v}^{h} \in H^{1}(\Om)^{n}, ~
  \bs{v}^{h} \cdot \bs{n} = 0 \, \text{on} \, \G_{1} \cup \G_{2} \cup \G_{c}
  \, ,
  \\ 
  {\bs{v}^h}_{|T} \in \mathcal{P}^m(T)  \, , \,
  \forall T \in \mathcal{T}^{h} \, \text{and} \, i \in \{1,2\} \bigg\},
\end{split}
\end{equation}
\begin{equation}
  Q^h_{\alpha_{i}}(\bar{\Om}_{i}, \bar{\mathcal{T}}_{i}^{h} ) := \left\{ q^h_{\alpha_{i}} \in L_2(\bar{\Om}_{i})  :
  q^h_{\alpha_{i}} = \alpha_{i} \, q^{h}, q^{h} \in L_{2}(\bar{\Om}_{i}) \, | \, \ {q^h}_{|T} \in \mathcal{P}^{m-1}(T)  \, , \, \forall T \in \bar{\mathcal{T}}_{i}^{h} \, \text{and} \, i \in \{1,2\}\right\},
\end{equation}
\end{subequations}
where $\bar{\Om}_{i} = \ti{\Om}_{i} \cup \Om_{c}$ and $\bar{\mathcal{T}}_{i}^{h} = \ti{\mathcal{T}}_{i}^{h} \cup \mathcal{T}^{h}_{c}$.


\subsection{Discrete material representation}
\label{sec_materials}

In addition to the material-dependent volume fractions \eqref{eq_alpha}, we
introduce material-dependent densities $\rho_i$ and specific internal
energies $e_i$,  discretized as smooth polynomials on each element.
These variables allow to exactly represent a given
discontinuity inside an element, which would be impossible with only a single
density or energy.
These fields are initialized in the following manner:
\begin{equation}
\label{eq_init_rho_e}
  \rho_i(\bs{x}, t = 0) =
  \begin{cases}
    \rho_0(\bs{x})|_i & \alpha_i(\bs{x}, 0) > 0 \,, \\
    0 &                 \alpha_i(\bs{x}, 0) = 0 \,,
  \end{cases}
  \qquad
  e_i(\bs{x}, t=0) =
  \begin{cases}
    e_0(\bs{x})|_i & \alpha_i(\bs{x}, 0) > 0 \,, \\
    0 &              \alpha_i(\bs{x}, 0) = 0 \,.
  \end{cases}
\end{equation}
Here $\rho_0|_i$ and $e_0|_i$ provide additional values for material $i$ in
mixed elements, at points where the material's values are not defined by the initial conditions, i.e., $\rho_0|_i$ and $e_0|_i$ extend the initial
conditions for material $i$ wherever $\alpha_i(\bs{x}, t) > 0$.
The time evolution of $\alpha_i$ in an element $T \in \mathcal{T}^h$ is
computed by \eqref{eq_alpha};
this formula and $d \eta / dt = 0$ imply that changes in the volume
fractions will occur due to local deformations of the mixed elements.
While approaches based on closure models \cite{Waltz2021, Tomov2016} perform
more complex operations to evolve volume fractions, in this work their evolution is based
exclusively on the fluid velocity and we rely on the methods presented in
Section \ref{sec_weak_form} for the enforcement of interface conditions in the 
momentum and internal energy equations.

The evolution of material-specific densities $\rho_i$ is governed by
the principle of exact local mass conservation:
\begin{equation}
\label{eq_rho}
  \rho_i(\bs{x}, t) =
    \frac{\rho_i(\bs{x}_0, 0) \alpha_i(\bs{x}_0, 0)}
         {J(\bs{x}, t) \alpha_i(\bs{x}, t)} \, ,
\end{equation}
where $\bs{\varphi}(\bs{x}_0, t) = \bs{x}$ and 
$\alpha_i(\bs{x}, t) \rho_i(\bs{x}, t)$ is the local mass of material $i$.
Finally, the evolution of the specific internal energies $e_i$ is governed
by the weak forms of the internal energy equations as explained in
Section \ref{sec_weak_form}.


\subsection{A derivation of the weighted-shifted variational formulation}
\label{sec_weak_form}

\begin{figure}
	\begin{subfigure}[!htb]{.4\textwidth}
		\begin{tikzpicture}[scale=0.5]
		\path [draw=orange,fill=gray!30!,line width=1.5pt,name] plot coordinates {  (2,0)  (7,0)  (7,5) (2,5) (2,0)};
		\path [draw=black,fill=gray!30!,line width=1.5pt,name] plot coordinates {  (-3,0)  (2,0)  (2,5) (-3,5) (-3,0)};
		\path [draw=black,line width=1.5pt,name] plot coordinates {  (7,0)  (12,0)  (12,5) (7,5) (7,0)};
		\draw [line width = 0.5mm,blue, name path=true] plot[smooth] coordinates {(4,-0.5) (3,2.5) (4.5,5.5)};
		\draw [-stealth](7,3) -- (3.5,3.5);
		\draw [-stealth](7,2) -- (9,2);
		\node[text width=5cm] at (6.1,0.75) {\large${\color{black}\ti{\G}_{1}}$};
		\node[text width=5cm] at (6.2,4) {\large${\color{black}+}$};
		\node[text width=5cm] at (7.2,4) {\large${\color{black}-}$};
		\node[text width=5cm] at (11.2,4) {\large${\color{black}+}$};
		\node[text width=5cm] at (12.2,4) {\large${\color{black}-}$};
		\node[text width=5cm] at (12.25,0.75) {\large${\color{black}\ti{\G}_{2}}$};
		\node[text width=5cm] at (10.5,5.6) {\large${\color{orange}{\mathcal{E}}^{o}_{c}}$};
		\node[text width=5cm] at (9,5.95) {\large${\color{blue}{\G}_{c}}$};
		\node[text width=5cm] at (3.75,1.75) {\large${\color{black}\ti{\Om}_{1}}$};
		\node[text width=5cm] at (7.35,0.75) {\large${\color{black}\alpha_{1}}$};
		\node[text width=5cm] at (8.8,1.5) {\large${\color{black}\alpha_{2}}$};
	\node[text width=5cm] at (10.5,2.75) {\large${\color{black}\bs{d}}$};
		\node[text width=5cm] at (14.5,1.75)
		 {\large${\color{black}\ti{\Om}_{2}}$};
		 \node[text width=5cm] at (12.5,2.5) {\large${\color{black}\ti{\bs{n}}^{+}}$};
		\end{tikzpicture}
		\caption{Case 1: Treating the cut elements as part of material $\#1$
             weighted by $\alpha_{1}$.}
	\label{fig:Case1}
	\end{subfigure}
	\hspace{2cm}
	\begin{subfigure}[!htb]{.4\textwidth}
		\begin{tikzpicture}[scale=0.5]
		\path [draw=orange,fill=gray!30!,line width=1.5pt,name] plot coordinates {  (2,0)  (7,0)  (7,5) (2,5) (2,0)};
		\path [draw=black,line width=1.5pt,name] plot coordinates {  (-3,0)  (2,0)  (2,5) (-3,5) (-3,0)};
		\path [draw=black,fill=gray!30!,line width=1.5pt,name] plot coordinates {  (7,0)  (12,0)  (12,5) (7,5) (7,0)};
		\draw [line width = 0.5mm,blue, name path=true] plot[smooth] coordinates {(4,-0.5) (3,2.5) (4.5,5.5)};
		\draw [-stealth](2,3) -- (3.1,2.75);
			\draw [-stealth](7,2) -- (9,2);
		\node[text width=5cm] at (6.1,0.75) {\large${\color{black}\ti{\G}_{1}}$};
		\node[text width=5cm] at (6.2,4) {\large${\color{black}+}$};
		\node[text width=5cm] at (7.2,4) {\large${\color{black}-}$};
		\node[text width=5cm] at (11.2,4) {\large${\color{black}+}$};
		\node[text width=5cm] at (12.2,4) {\large${\color{black}-}$};
		\node[text width=5cm] at (12.25,0.75) {\large${\color{black}\ti{\G}_{2}}$};
		\node[text width=5cm] at (10.5,5.6) {\large${\color{orange}{\mathcal{E}}^{o}_{c}}$};
		\node[text width=5cm] at (9,5.95) {\large${\color{blue}{\G}_{c}}$};
		\node[text width=5cm] at (3.75,1.75) {\large${\color{black}\ti{\Om}_{1}}$};
		\node[text width=5cm] at (7.35,0.75) {\large${\color{black}\alpha_{1}}$};
		\node[text width=5cm] at (8.8,1.5) {\large${\color{black}\alpha_{2}}$};
		\node[text width=5cm] at (7.25,2.25) {\large${\color{black}\bs{d}}$};
		\node[text width=5cm] at (14.5,1.75) {\large${\color{black}\ti{\Om}_{2}}$};
	  \node[text width=5cm] at (12.5,2.5) {\large${\color{black}\ti{\bs{n}}^{+}}$};
		\end{tikzpicture}
		\caption{Case 2: Treating the cut elements as part of material $\#2$
             weighted by $\alpha_{2}$.}
	\label{fig:Case2}
	\end{subfigure}
	\caption{The two cases for treating the cut elements.}
	\label{fig:Cases}
\end{figure}
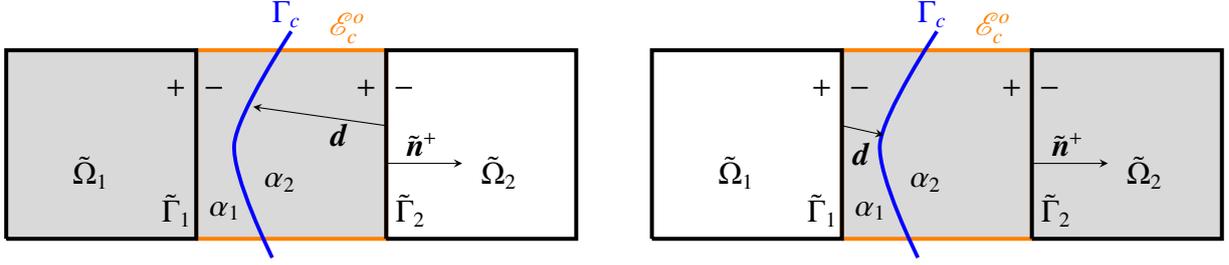

Let
$\bar{\mathcal{E}}_{1} =
\ti{\G}_{1} \cup \mathcal{E}^{o}_{c} \cup \mathcal{E}_{1}^{o}$ and $\bar{\mathcal{E}}_{2} =
\ti{\G}_{2} \cup \mathcal{E}^{o}_{c} \cup \mathcal{E}_{2}^{o}$
be the unions of the edges/faces of the grid in two/three dimensions where the materials $\#1$ and $\#2$ are present, respectively, see Figure \ref{fig:Cases}.
That is to say, if an edge/face is in $\bar{\mathcal{E}}_{i}$, then material $i$ is present at least on one of the sides of that edge/face (possibly both).
For each side, the possibly discontinuous values of a field are indicated with the superscripts $+$ and $-$, as shown in Figure~\ref{fig:Cases}.
Across each side of the grid, the pressure $p_{1}$ of material $\#1$ and $p_{2}$ of material $\#2$ need to {\it separately} satisfy the following continuity conditions:
\begin{subequations}
	\label{eq:continuity_Lagrangian}
	\begin{align}
	[p_{1}] &= \bs{0} , \quad \mbox{ on } \bar{\mathcal{E}}_{1} \; ,
	\label{eq:continuity_p1}
	\\
	[p_{2}] &= \bs{0} , \quad \mbox{ on } \bar{\mathcal{E}}_{2} \; ,
	\label{eq:continuity_p2}
	\end{align}
\end{subequations}
which will be enforced weakly, in what follows.

Next we introduce a definition of the jump of a field across an edge/face which is complementary to~\eqref{eq_interface_cond}.
For $i=1$ or $2$, $[\zeta_i] = \zeta_{i}^{+} \bs{n}^{+} + \zeta_{i}^{-}  \bs{n}^{-} $ for any scalar field $\zeta_i$ and $[\bs{\zeta}_i] = \bs{\zeta}_{i}^{+} \cdot \bs{n}^{+} + \bs{\zeta}_{i}^{-} \cdot \bs{n}^{-} $ for any vector field $\bs{\zeta}_i$.
Observe that while~\eqref{eq_interface_cond} express jump conditions across a
material interface, the previous definitions express a jump condition for a field associated with a single material across an edge/face of the grid. 
See again Figure~\ref{fig:Cases}.

To derive the weak form of the momentum equation, we start by considering material $\#1$, as sketched in Case 1 of Figure \ref{fig:Case1}.
Multiplying the strong form of the momentum equation by the test function $\bs{\psi}^{h}_{\alpha_{1}} \in \bs{V}_{\alpha_{1}}^h$ and integrating over the generic element $T$ of the finite element mesh yields
\begin{align}
	(\rho_{1} \, \dot{\bs{v}} , \alpha_{1} \, \bs{\psi}  )_{T}
	- 	(p_{1} , \alpha_{1} \, \nabla \cdot \bs{\psi}   )_{T}
  +  (\mu_1  \nabla^s \bs{v}, \alpha_1 \nabla \bs{\psi})_{T}  
  +  \langle  p_{1}  , \alpha_{1} \,  \bs{\psi}\cdot \bs{n}  \rangle _{\partial T}
	= 0.
	\end{align}
The term $(\mu_1 \nabla^s \bs{v}, \alpha_1 \nabla \bs{\psi})_{T} $ contains the artificial viscosity tensor $\sigma_a= \mu_1 \nabla^s \bs{v}$:
its complete definition can be found in \cite{Dobrev2012}.
Since the interface conditions \eqref{eq_interface_cond} do not involve the
artificial viscosity tensor, this term does not appear in any of the face integrals.
This is equivalent to weakly enforce the continuity of the normal trace of the artificial viscosity on each of the mesh faces.
Summing over all the elements we have
\begin{equation}
\begin{split}
		\label{eq:baseDGO1d}
	&( \rho_{1} \, \dot{\bs{v}},  \alpha_{1} \, \bs{\psi} )_{\bar{\Om}_{1}}
	- 	( p_{1} ,  \alpha_{1} \, \nabla \cdot \bs{\psi}   )_{\bar{\Om}_{1} }
   +   (\mu_1 \alpha_1 \nabla^s \bs{v}, \nabla \bs{\psi})_{\bar{\Om}_{1} }
\\ &
	 + \langle [p_{1}] ,  \{\alpha_{1}\}_{\alpha^{-}_{1}} \bs{\psi} \rangle _{\bar{\mathcal{E}}_{1} }
	+ \langle \{p_{1}\}_{1-\alpha^{-}_{1}} ,  [\alpha_{1}] \bs{\psi} \rangle _{\bar{\mathcal{E}}_{1}}
	+ \langle p^{+}_{1} , \alpha^{+}_{1} \bs{\psi} \cdot \ti{\bs{n}}^{+}  \rangle _{\ti{\G}_{2}}
	= 0,
\end{split}
\end{equation}
where  $\{\zeta\}_{\gamma} = \gamma \, \zeta^{+} + (1-\gamma) \, \zeta^{-}$
for $\gamma \in [0,1]$ and any scalar or vector field $\zeta$.
Enforcing \eqref{eq:continuity_p1} gives
\begin{equation}
\begin{split}
\label{eq:baseDGO1d2}
& (\rho_{1} \, \dot{\bs{v}},  \alpha_{1}  \bs{\psi} )_{\bar{\Om}_{1} }
	- 	(p_{1} , \alpha_{1} \, \nabla \cdot \bs{\psi}  )_{\bar{\Om}_{1} }
   +   (\mu_1 \alpha_1 \nabla^s \bs{v}, \nabla \bs{\psi})_{\bar{\Om}_{1} }
\\ &
	+ \langle \{p_{1}\}_{1-\alpha^{-}_{1}} ,  [\alpha_{1}] \, \bs{\psi} \rangle _{\ti{\G}_{1} \cup \mathcal{E}^{o}_{c}}
	+ \langle p^{+}_{1}   ,  \alpha^{+}_{1} \, \bs{\psi} \cdot  \ti{\bs{n}}^{+}  \rangle _{\ti{\G}_{2}}
	= 0.
\end{split}
\end{equation}
Considering now material $\#2$, as depicted in Case 2 of Figure \ref{fig:Case2}, we test the momentum equation with $\bs{\psi}^{h}_{\alpha_{2}} \in \bs{V}_{\alpha_{2}}^h$, and we follow analogous steps as above to enforce \eqref{eq:continuity_p2}:
	\begin{align}
		\label{eq:baseDGO1d3}
	(\rho_{2} \, \dot{\bs{v}} , \alpha_{2} \, \bs{\psi}   )_{\bar{\Om}_{2}}
	- 	(  p_{2}, \alpha_{2} \, \nabla \cdot \bs{\psi}  )_{\bar{\Om}_{2}}
   +   (\mu_2 \alpha_2 \nabla^s \bs{v}, \nabla \bs{\psi})_{\bar{\Om}_{2} }
	+ \langle \{p_{2}\}_{1-\alpha^{+}_{2}} ,  [\alpha_{2}] \, \bs{\psi} \rangle _{\ti{\G}_{2} \cup \mathcal{E}^{o}_{c}}
	+ \langle p^{-}_{2} ,  \alpha^{-}_{2}  \, \bs{\psi} \cdot \ti{\bs{n}}^{-}   \rangle _{\ti{\G}_{1}}
	= 0.
	\end{align}
Summing~\eqref{eq:baseDGO1d2} and~\eqref{eq:baseDGO1d3} and adding and subtracting $\langle \alpha^{-}_{2}\,p_{1}^{-} ,  \, \bs{\psi} \cdot \ti{\bs{n}}^{-}   \rangle _{\ti{\G}_{1}} $ and $\langle \alpha^{+}_{1}\,p_{2}^{+} ,  \, \bs{\psi} \cdot \ti{\bs{n}}^{+}   \rangle _{\ti{\G}_{2}} $ gives
\begin{equation}
\begin{split}
\label{eq:baseDGO1d4}
(\rho_{i}   \dot{\bs{v}} ,  \alpha_{i} \bs{\psi} )_{\Om}
- 	( p_{i}  , \alpha_{i} \nabla \cdot \bs{\psi}  )_{\Om}
+   (\mu_i \alpha_i \nabla^s \bs{v}, \nabla \bs{\psi})_{\Om}
&+ \langle \jump{p^{-}} \cdot \bs{n}_{1},  \alpha^{-}_{2}  \bs{\psi} \cdot \ti{\bs{n}}^{+} \rangle _{\ti{\G}_{1}}
+ \langle \jump{p^{+}} \cdot \bs{n}_{2} ,  \alpha^{+}_{1}  \bs{\psi} \cdot \ti{\bs{n}}^{-} \rangle _{\ti{\G}_{2}}
\\ & 
+ \langle \{p_{1}\}_{1-\alpha^{-}_{1}} [\alpha_{1}] + \{p_{2}\}_{1-\alpha^{+}_{2}} [\alpha_{2}]   ,  \bs{\psi} \rangle _{\mathcal{E}^{o}_{c}}
= 0.
\end{split}
\end{equation}
Enforcing the shifted conditions \eqref{eq:Finalu-g} on $\ti{\G}_{1}$ and $\ti{\G}_{2}$ yields
\begin{equation}
\begin{split}
\label{eq:baseDGO1d5}
(\rho_{i}  \, \dot{\bs{v}}, \alpha_{i} \, \bs{\psi}   )_{\Om}
&- 	(p_{i} , \alpha_{i}  \, \nabla \cdot \bs{\psi}  )_{\Om}
+   (\mu_i \alpha_i \nabla^s \bs{v}, \nabla \bs{\psi})_{\Om}
\\ &
- \langle \jump{ \sum_{k = 1}^{m-1}  \frac{\mathcal{D}^{k}_{\bs{d}} \, p^{-} }{k!}  } \cdot \bs{n}_{1} ,  \alpha^{-}_{2} \, \bs{\psi} \cdot \ti{\bs{n}}^{+} \rangle _{\ti{\G}_{1}}
 -\langle \jump{ \sum_{k = 1}^{m-1}  \frac{\mathcal{D}^{k}_{\bs{d}} \, p^{+} }{k!}  } \cdot \bs{n}_{2}  ,  \alpha^{+}_{1} \, \bs{\psi}  \cdot \ti{\bs{n}}^{-} \rangle _{\ti{\G}_{2}}
\\ & + \langle \{p_{1}\}_{1-\alpha^{-}_{1}} [\alpha_{1}]  + \{p_{2}\}_{1-\alpha^{+}_{2}} \, [\alpha_{2}]   , \bs{\psi} \rangle _{\mathcal{E}^{o}_{c}}
= 0.
\end{split}
\end{equation}
Noting that $ [\alpha \, p ] =  \{ \alpha \}_{\gamma} \, [ p ] + [\alpha ] \, \{ p \}_{1-\gamma}  $ and $\alpha_{2} = 1-\alpha_{1}$, and recalling  \eqref{eq:continuity_Lagrangian} and \eqref{eq:Finalu-g}, we have
\begin{align}
\label{eq:expression_cut}
\{p_{1}\}_{1-\alpha^{-}_{1}} [\alpha_{1}]  + \{p_{2}\}_{1-\alpha^{+}_{2}} \, [\alpha_{2}] &= 
[ \alpha_{1} \, p_{1} + \alpha_{2} \, p_{2} ]
\nonumber \\ &
=\left(\alpha^{+}_{1} \, p^{+}_{1} - \alpha^{-}_{1} \, p^{-}_{1} + \alpha^{+}_{2} \, p^{+}_{2} - \alpha^{-}_{2} \, p^{-}_{2} \right) \ti{\bs{n}}^{+}
\nonumber	\\ &
= \left( \frac{(\alpha^{+}_{1} - \alpha^{-}_{1})}{2}\, p^{+}_{1} + \frac{(\alpha^{+}_{1} - \alpha^{-}_{1})}{2} \, p^{-}_{1} - \frac{(\alpha^{+}_{1} - \alpha^{-}_{1})}{2}\, p^{+}_{2} - \frac{(\alpha^{+}_{1} - \alpha^{-}_{1})}{2} \, p^{-}_{2} \right) \ti{\bs{n}}^{+}
\nonumber	\\ &
=  \frac{\alpha^{+}_{1} - \alpha^{-}_{1}}{2}\left(  \jump{p^{-}}\cdot \bs{n}_{1}   + \jump{p^{+}} \cdot \bs{n}_{1} \right)  \ti{\bs{n}}^{+}
\nonumber	\\ &
= -\frac{\alpha^{+}_{1} - \alpha^{-}_{1}}{2}\left( \jump{ \sum_{k = 1}^{m-1}  \frac{\mathcal{D}^{k}_{\bs{d}} \, p^{-} }{k!}  } \cdot \bs{n}_{1}   +\jump{ \sum_{k = 1}^{m-1}  \frac{\mathcal{D}^{k}_{\bs{d}} \, p^{+} }{k!}  } \cdot \bs{n}_{1} \right)  \ti{\bs{n}}^{+}
\; .
\end{align}
Substituting~\eqref{eq:expression_cut} into~\eqref{eq:baseDGO1d5} yields the
final two-material momentum equation with shifted interface conditions
(note that we also added back the ignored volumetric viscosity terms):
\begin{equation}
\begin{split}
\label{eq_weak_v}
(\rho_{i}  \, \dot{\bs{v}}, \alpha_{i} \, \bs{\psi}   )_{\Om}
& - 	(p_{i} , \alpha_{i}  \, \nabla \cdot \bs{\psi}  )_{\Om}
  +   (\mu_i \alpha_i \nabla^s \bs{v}, \nabla \bs{\psi})_{\Om}
\\ &
- \langle \jump{ \sum_{k = 1}^{m-1}  \frac{\mathcal{D}^{k}_{\bs{d}} \, p^{-} }{k!}  } \cdot \bs{n}_{1} ,  \alpha^{-}_{2} \, \bs{\psi} \cdot \ti{\bs{n}}^{+} \rangle _{\ti{\G}_{1}}
-\langle \jump{ \sum_{k = 1}^{m-1}  \frac{\mathcal{D}^{k}_{\bs{d}} \, p^{+} }{k!}  } \cdot \bs{n}_{2}  ,  \alpha^{+}_{1} \, \bs{\psi}  \cdot \ti{\bs{n}}^{-} \rangle _{\ti{\G}_{2}}
\\ & 
-\langle \frac{\alpha^{+}_{1} - \alpha^{-}_{1}}{2}\left( \jump{ \sum_{k = 1}^{m-1}  \frac{\mathcal{D}^{k}_{\bs{d}} \, p^{-} }{k!}  } \cdot \bs{n}_{1}   +\jump{ \sum_{k = 1}^{m-1}  \frac{\mathcal{D}^{k}_{\bs{d}} \, p^{+} }{k!}  } \cdot \bs{n}_{1} \right)   , \bs{\psi}   \cdot \ti{\bs{n}}^{+} \rangle_{\mathcal{E}^{o}_{c}}
= \;
0
\,.
\end{split}
\end{equation}
\\
To derive the weak form of the internal energy equation, we follow the typical steps of a local discontinuous Galerkin (LDG) formulation.
Starting with material \#$1$ and multiplying the strong form of the
internal equation by the test function
$\phi^{h}_{\alpha_{1}} \in Q_{\alpha_{1}}^h$ and
over a generic element $T \in \bar{\Om}_{1}$ yields
	\begin{align}
	(\rho_{1} \dot{e}_{1}  , \alpha_{1} \, \phi  )_{T}
	- ( \bs{v}  ,  \nabla (\alpha_{1} \, \phi \, p_{1}) )_{T}
	+\langle\bs{v}\cdot\bs{n}  ,  \alpha_{1} \, \phi \, p_{1} \rangle_{\partial T}
   - ( \mu_1 \alpha_1 \nabla^s \bs{v}, \nabla \bs{v} )_{T}
	= 0.
	\end{align}
	Replacing the boundary velocity trace, $\bs{v}$, by a \textit{numerical flux},  $\hat{\bs{v}}$, yields
\begin{align}
	(\rho_{1} \dot{e}_{1}  ,  \alpha_{1} \, \phi  )_{T}
		- ( \bs{v} ,  \nabla (\alpha_{1} \, \phi \, p_{1})    )_{T}
		+\langle \hat{\bs{v}}\cdot\bs{n}  ,  \alpha_{1} \, \phi \, p_{1} \rangle_{\partial T}
   - ( \mu_1 \alpha_1 \nabla^s \bs{v}, \nabla \bs{v} )_{T}
= 0.
\end{align}
	The numerical flux $\hat{\bs{v}}$ accomplishes the linking of the numerical solution between neighboring elements, and will be defined momentarily. Summing over all the elements and integrating back by parts we have
\begin{equation}
\label{eq:baseDGO1}
\begin{split}
		(\rho_{1} \dot{e}_{1}  ,  \alpha_{1} \, \phi  )_{\bar{\Om}_{1}} 
	 + (p_{1} \, \nabla \cdot \bs{v}  ,   \alpha_{1} \,  \phi \,  )_{\bar{\Om}_{1}}
  &- ( \mu_1 \alpha_1 \nabla^s \bs{v}, \nabla \bs{v} )_{\bar{\Om}_{1}}
		\\ &
		+  \langle \jump{\hat{\bs{v}}} \cdot \ti{\bs{n}}^{+} ,  \{\alpha_{1} \, \phi \, p_{1}\}_{\gamma_{1}}  \rangle _{\ti{\G}_{1} \cup \ti{\G}_{2} \cup \mathcal{E}^{o}_{c}}
		+  \langle \{\hat{\bs{v}} \} _{1-\gamma_{1}}  ,  [\alpha_{1} \, \phi \, p_{1}] \rangle _{\ti{\G}_{1} \cup \ti{\G}_{2} \cup \mathcal{E}^{o}_{c}}
		\\ & 
			-  \langle \jump{\bs{v}} \cdot \ti{\bs{n}}^{+} ,  \{\alpha_{1} \, \phi \, p_{1}\}_{\gamma_{1}}  \rangle _{\ti{\G}_{1} \cup \ti{\G}_{2} \cup \mathcal{E}^{o}_{c}}
	-  \langle \{\bs{v} \} _{1-\gamma_{1}} ,  [\alpha_{1} \, \phi \, p_{1} ]  \rangle _{\ti{\G}_{1} \cup \ti{\G}_{2} \cup \mathcal{E}^{o}_{c}} = 0.
\end{split}
\end{equation}
We define the normal and tangential components of $\hat{\bs{v}}$ as
\begin{equation}
	\begin{aligned}
	\label{eq:decompflux}
	\jump{\hat{\bs{v}}} \cdot \bs{n}_{1} = - \jump{ \sum_{k = 1}^{m}  \frac{\mathcal{D}^{k}_{\bs{d}} \, \bs{v} }{k!}   } \cdot \bs{n}_{1},  \qquad
	\jump{\hat{\bs{v}}} \cdot \bs{\tau}_{1}^{i}=   \jump{\bs{v}} \cdot \bs{\tau}_{1}^{i}, 
	\end{aligned}
\end{equation}
where $\bs{\tau}_{1}^{i}$ is a unit tangent vector to the true interface (repeated index notation is implied). The normal component of the numerical flux, $\jump{\hat{\bs{v}}} \cdot \bs{n}_{1}$, mimics \eqref{eq:Finalbsu-g-equiv} while  $\jump{\hat{\bs{v}}} \cdot \bs{\tau}_{1}^{i}$
is the same as the body-fitted case.
Noting that the relation
$\jump{\hat{\bs{v}}} = (\jump{\hat{\bs{v}}} \cdot \bs{n}_{1})
\bs{n}_{1}+(\jump{\hat{\bs{v}}}  \cdot \bs{\tau}_{1}^{i})\bs{\tau}_{1}^{i}$,
holds for any vector $\hat{\bs{v}}$,
the interface conditions for $\bs{v}$ can be reflected by setting
\begin{equation}
\label{eq:v_flux}
\begin{aligned}
\jump{\hat{\bs{v}}} = - \left(\jump{ \sum_{k = 1}^{m}  \frac{\mathcal{D}^{k}_{\bs{d}} \, \bs{v} }{k!} } \cdot \bs{n}_{1} \right) \bs{n}_{1} +
	(\jump{\bs{v}} \cdot \bs{\tau}_{1}^{i}) \bs{\tau}_{1}^{i} ,  \qquad
\avgb{\hat{\bs{v}}}_{1-\gamma} =   \avgb{\bs{v}}_{1-\gamma}, 
\end{aligned}
\end{equation}
where the choice for $\avgb{\hat{\bs{v}}}_{1-\gamma}$ is exactly the same as that for the body-fitted.
Once the above substitutions are made, and using the fact that $\jump{\bs{v}} = 0$ on
internal faces as the velocity space $\mathcal{V}$ is continuous,
the discrete variational form reads
	\begin{align}
	\label{eq:dg_bf_intermediate_5O1}
	(\rho_{1} \dot{e}_{1}  , \alpha_{1} \, \phi  )_{\bar{\Om}_{1}} 
+(p_{1} \, \nabla \cdot \bs{v}  ,   \alpha_{1} \,  \phi \,  )_{\bar{\Om}_{1}}
- ( \mu_1 \alpha_1 \nabla^s \bs{v}, \nabla \bs{v} )_{\bar{\Om}_{1}}
	-\avg{  \jump { \bs{S}_{h}^{m} \bs{v}} \cdot \bs{n}_{1} (\ti{\bs{n}}^{+} \cdot \bs{n}_{1}),   \{\alpha_{1} \, \phi \, p_{1}\}_{\gamma_{1}} }_{\ti{\G}_{1} \cup \ti{\G}_{2} \cup \mathcal{E}^{o}_{c}}
	= 0.
\end{align}
Finally, to improve stability, we add a diffusion-type double jump terms at the
surrogate interfaces to obtain the final internal energy variational form
for material $\#1$:
\begin{equation}
\begin{split}
\label{eq_weak_e_1}
&	(\rho_{1} \dot{e}_{1}  , \alpha_{1} \phi  )_{\bar{\Om}_{1}} 
+(p_{1}  \nabla \cdot \bs{v} ,   \alpha_{1}   \phi  )_{\bar{\Om}_{1}}
- (\mu_1 \alpha_1 \nabla^s \bs{v}, \nabla \bs{v} )_{\bar{\Om}_{1}} 
-\avg{ \jump { \bs{S}_{h}^{m} \bs{v}} \cdot \bs{n}_{1} (\ti{\bs{n}}^{+} \cdot \bs{n}_{1})   ,  \{\alpha_{1} \phi  p_{1}\}_{\gamma_{1}} }_{\ti{\G}_{1} \cup \ti{\G}_{2} \cup \mathcal{E}^{o}_{c}}
\\ &
+ \langle (1-\alpha_1)
          \{ d ~ |\nabla \bs{v}| \}_{\gamma}
          \left(p_{1}^{+} - p_{2}^{-} +
                \sum_{k=1}^{m-1} \frac{\mathcal{D}^{k}_{\bs{d}} \,
                \left(p_1^{+} - p_2^{-} \right) }{k!} \right),
                [\tS^{m-1} \phi ] \cdot \ti{\bs{n}}^{+}        
  \rangle_{\ti{\G}_{1}} \\
&
+ \langle \alpha_1
          \{ d ~ |\nabla \bs{v}| \}_{\gamma}
          \left(p_{1}^{+} - p_{2}^{-} +
                \sum_{k=1}^{m-1} \frac{\mathcal{D}^{k}_{\bs{d}} \,
                \left(p_1^{+} - p_2^{-} \right) }{k!} \right),
                [\tS^{m-1} \phi ] \cdot \ti{\bs{n}}^{+}        
  \rangle_{\ti{\G}_{2}} = 0.
\end{split}
\end{equation}
Note that the diffusion terms in \eqref{eq_weak_e_1} are shifted, i.e.,
they represent the pressure jump at the true interface.
Similarly, to derive the weak form of the internal energy equation for
material $\#2$, we follow the same procedure as for material $\#1$,
but this time we multiply the strong form of the internal equation
by the test function $\phi^{h}_{\alpha_{2}} \in Q_{\alpha_{2}}^h$.
Thus, the final internal energy variational form for material $\#2$ reads
\begin{equation}
\begin{split}
\label{eq_weak_e_2}
&	(\rho_{2}  \dot{e}_{2} , \alpha_{2} \phi   )_{\bar{\Om}_{2}} 
	 + (p_{2}  \nabla \cdot \bs{v} , \alpha_{2} \phi  )_{\bar{\Om}_{2}} 
   - ( \mu_2 \alpha_2 \nabla^s \bs{v}, \nabla \bs{v} )_{\bar{\Om}_{2}}
-	\avg{  \jump { \bs{S}_{h}^{m} \bs{v}} \cdot \bs{n}_{1} (\ti{\bs{n}}^{+} \cdot \bs{n}_{1})  ,  \{\alpha_{2}  \phi  p_{2}\}_{\gamma_{2}} }_{\ti{\G}_{1} \cup \ti{\G}_{2} \cup \mathcal{E}^{o}_{c}}
\\ &
+ \langle \alpha_2
          \{ d ~ |\nabla \bs{v}| \}_{\gamma}
          \left(p_{1}^{+} - p_{2}^{-} +
                \sum_{k=1}^{m-1} \frac{\mathcal{D}^{k}_{\bs{d}} \,
                \left(p_1^{+} - p_2^{-} \right) }{k!} \right),
                [\tS^{m-1} \phi ] \cdot \ti{\bs{n}}^{+}          
  \rangle_{\ti{\G}_1} \\
&
+ \langle (1-\alpha_2)
          \{ d ~ |\nabla \bs{v}| \}_{\gamma}
          \left(p_{1}^{+} - p_{2}^{-} +
                \sum_{k=1}^{m-1} \frac{\mathcal{D}^{k}_{\bs{d}} \,
                \left(p_1^{+} - p_2^{-} \right) }{k!} \right),
                [\tS^{m-1} \phi ] \cdot \ti{\bs{n}}^{+}        
  \rangle_{\ti{\G}_{2}} = 0.
\end{split}
\end{equation}
The conservation properties of the above formulation are discussed in \ref{sec_conservation}.

\begin{remark}
Note that the mass matrices in
\eqref{eq_weak_v}, \eqref{eq_weak_e_1}, and \eqref{eq_weak_e_2}
remain constant in time due to \eqref{eq_rho}.
\end{remark}

%% file: results.tex
\section{Numerical Results}
\label{sec_results}

All simulations are performed in a customized version of the open-source
Laghos proxy application \cite{Laghos2019}, which is based on the
MFEM finite element library \cite{MFEM2021}.


\subsection{Smooth two-dimensional Taylor-Green Vortex}
\label{sec_tg}

The purpose of this example is to confirm that the method can retain
high-order convergence for a smooth problem.
The example can be seen as a patch-like test where a smooth single-material
simulation is altered by adding an artificial interface between two
identical materials.

\begin{figure}[t!]
\centerline
{
  \includegraphics[width=0.27\textwidth]{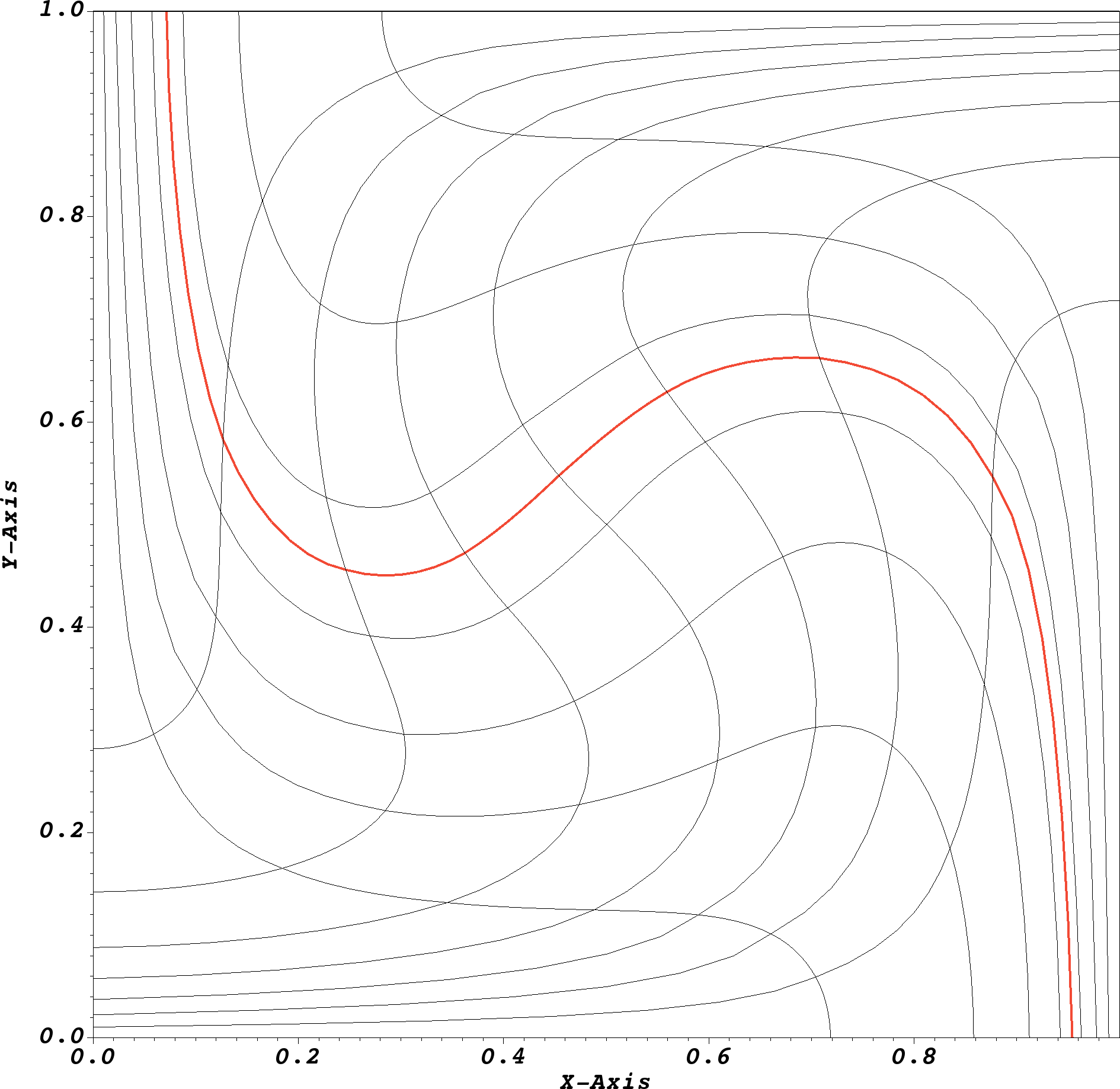}\hfil
  \includegraphics[width=0.27\textwidth]{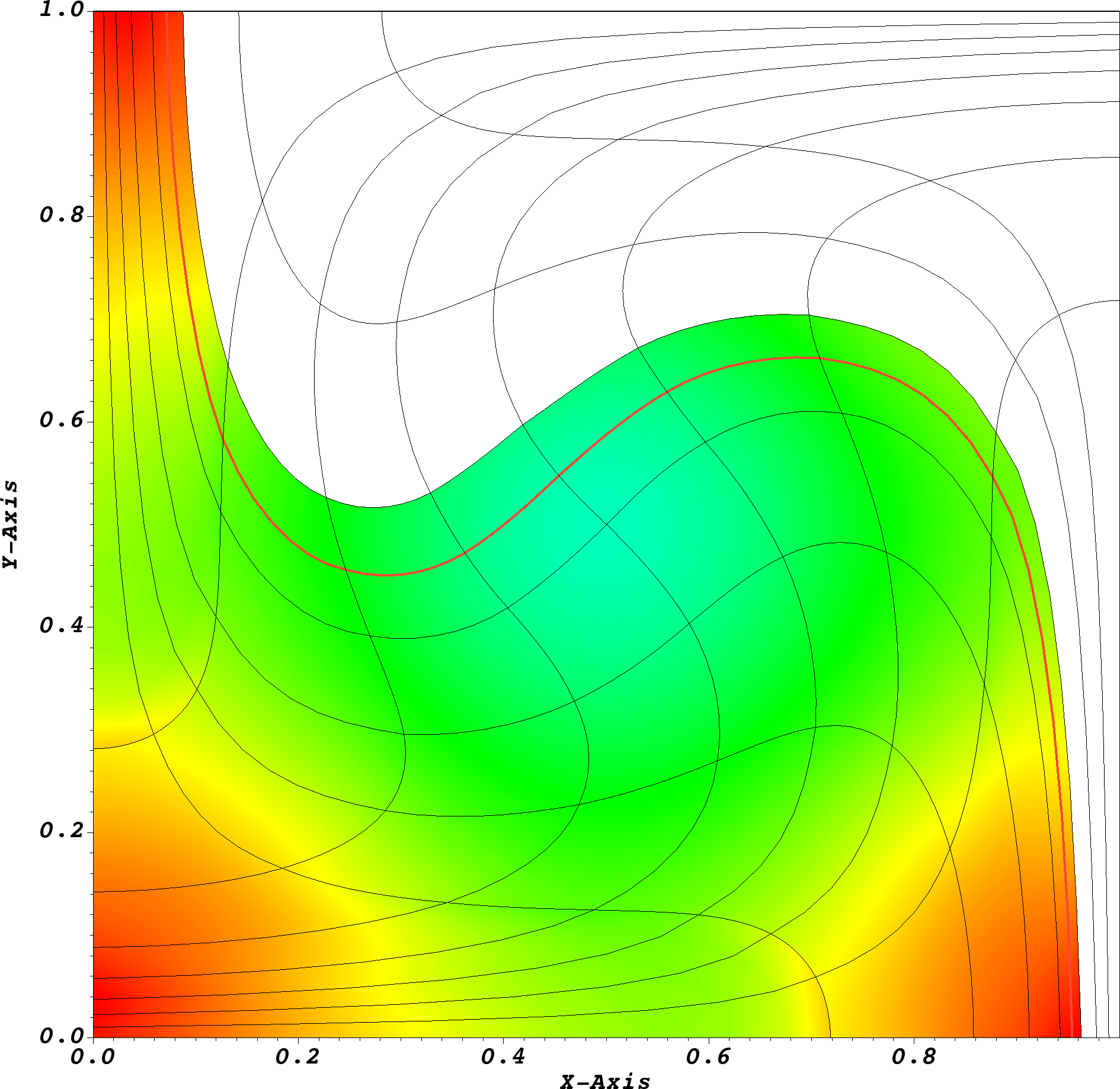} \hfil
  \includegraphics[width=0.27\textwidth]{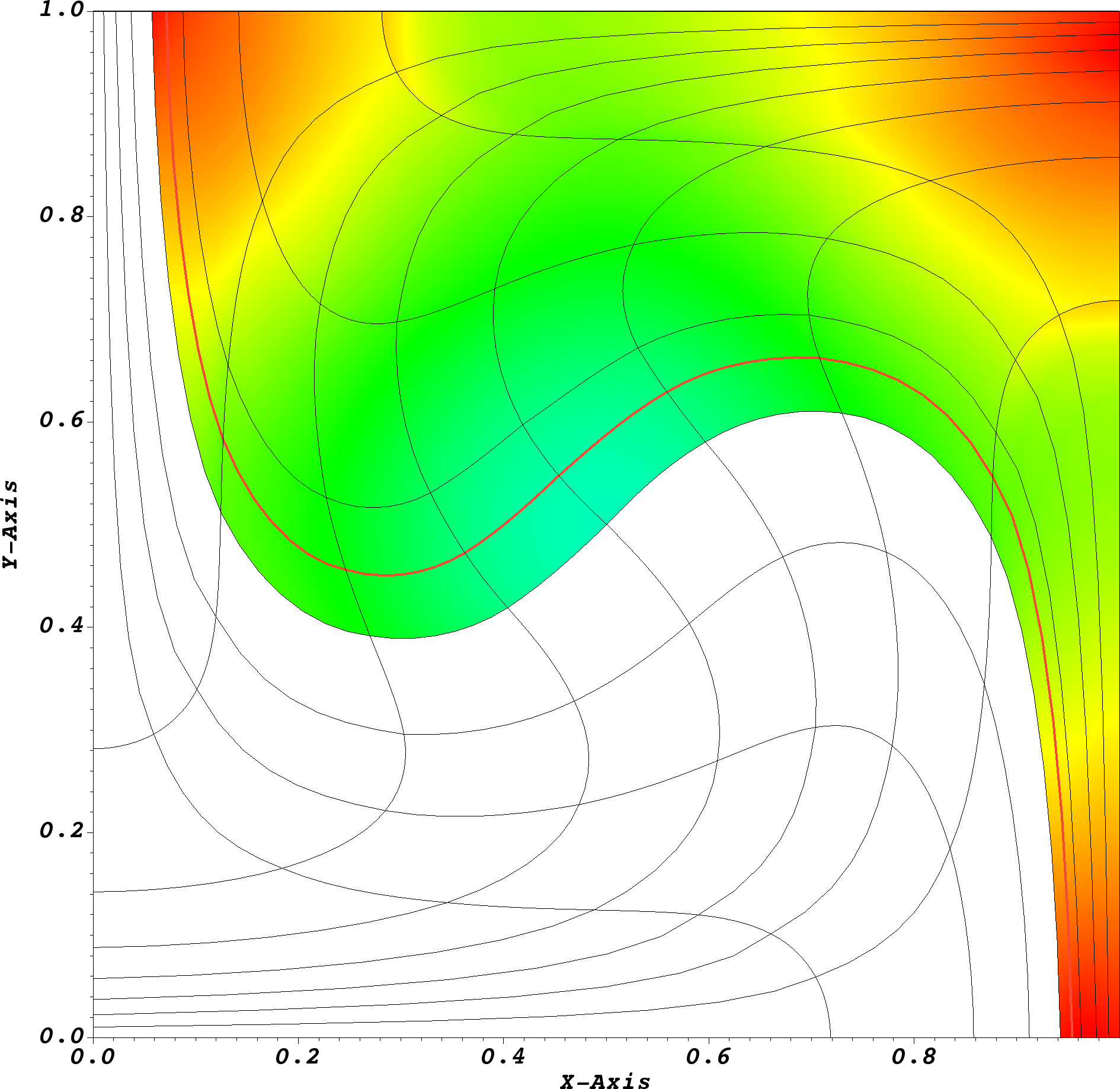}
}
\caption{Final mesh and interface position (left), and
         specific internal energies (middle and right)
         at $t=0.75$ for a $Q_3-Q_2$ discretization of
         the two-dimensional Taylor Green vortex.}
\label{fig_tg}
\end{figure}

The domain is $[0, 1] \times [0, 1]$ with $\bs{v} \cdot \bs{n} = 0$ (slip) boundary conditions.
For all tests, a vertical interface is initially
positioned at $x = 0.5 + \Delta x$.
The initial velocity, densities and pressures are set as:
\[
\bs{v} = \{ \sin(\pi x) \cos(\pi y), -\cos(\pi x) \sin(\pi y) \}, \quad
\rho_i = 1, \quad
p_i = \frac{\rho_i}{4}(\cos(2 \pi x) + \cos(2 \pi y)), \quad
\gamma_i = \frac{5}{3},
\]
where the material-specific quantities are initialized only on their
side of the interface, and in the cut elements.
The right-hand sides of the specific internal energy equations include the
following source terms:
\[
e_{i,src} = \frac{3 \pi}{8} \alpha_i
  \left( \cos(3 \pi x) \cos(\pi y) - \cos(\pi x) \cos(3 \pi y) \right).
\]
This setup guarantees that the velocity, material energies, and densities
do not change in time, i.e. $\partial / \partial t = 0$,
but they do vary along particle trajectories, i.e., $d / d t \neq 0$,
which allows to compare to the exact solution at any time;
see Section 8.1 in \cite{Dobrev2012} for more details.
Figure \ref{fig_tg} shows an example results for a
$Q_3-Q_2$ simulation with 2 mesh refinements ($8 \times 8$ elements).
As our shifted face terms contain jumps in $\nabla p$ at the surrogate
interfaces, the numerical discretization of $p_k$ introduces additional errors.
This test demonstrates that these errors converge with the expected rate as
we increase the polynomial order and take higher Taylor expansions of $p_k$.

\begin{figure}[b!]
\centerline
{
  \includegraphics[width=0.3\textwidth]{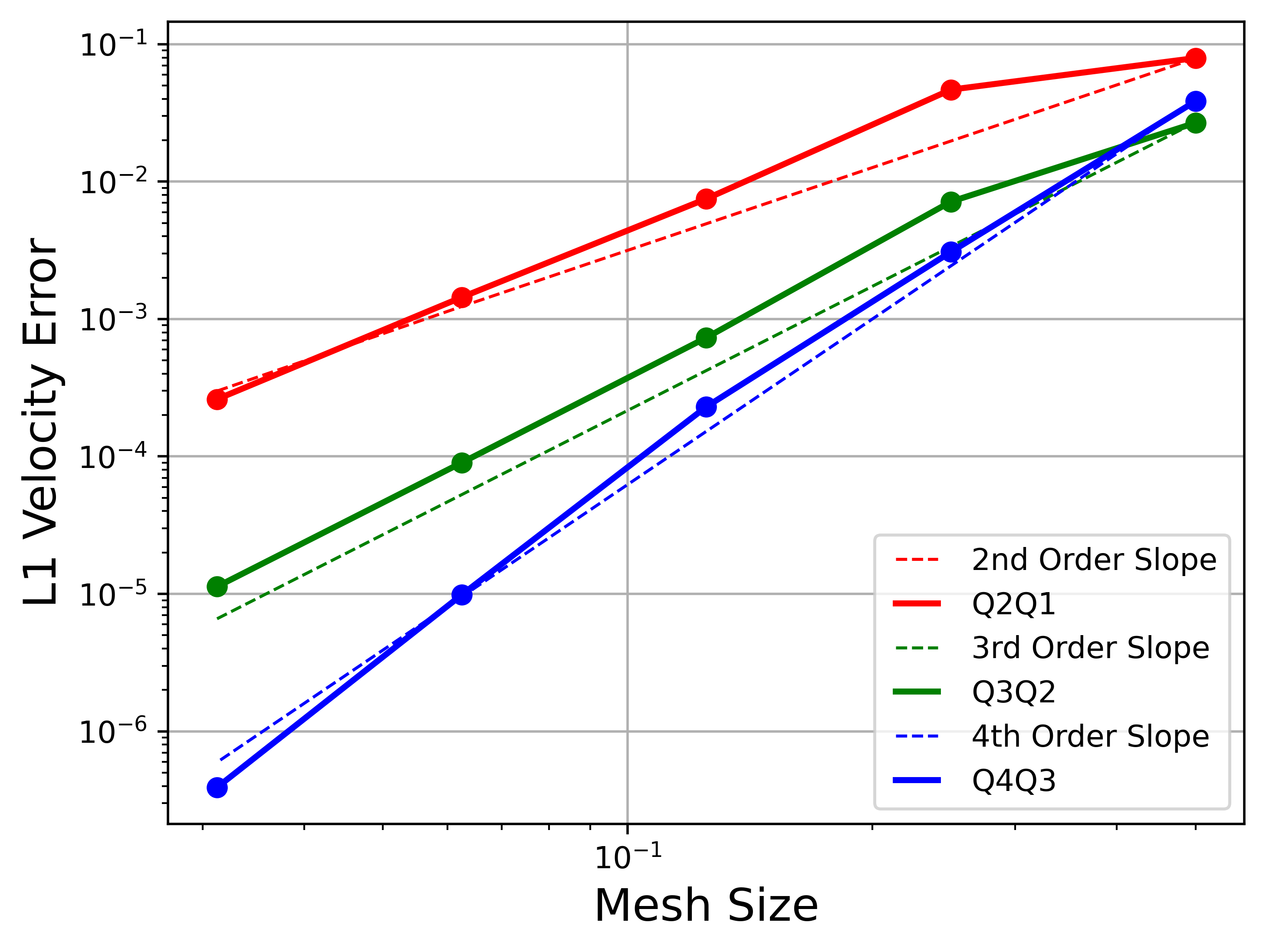}  \hfil
  \includegraphics[width=0.3\textwidth]{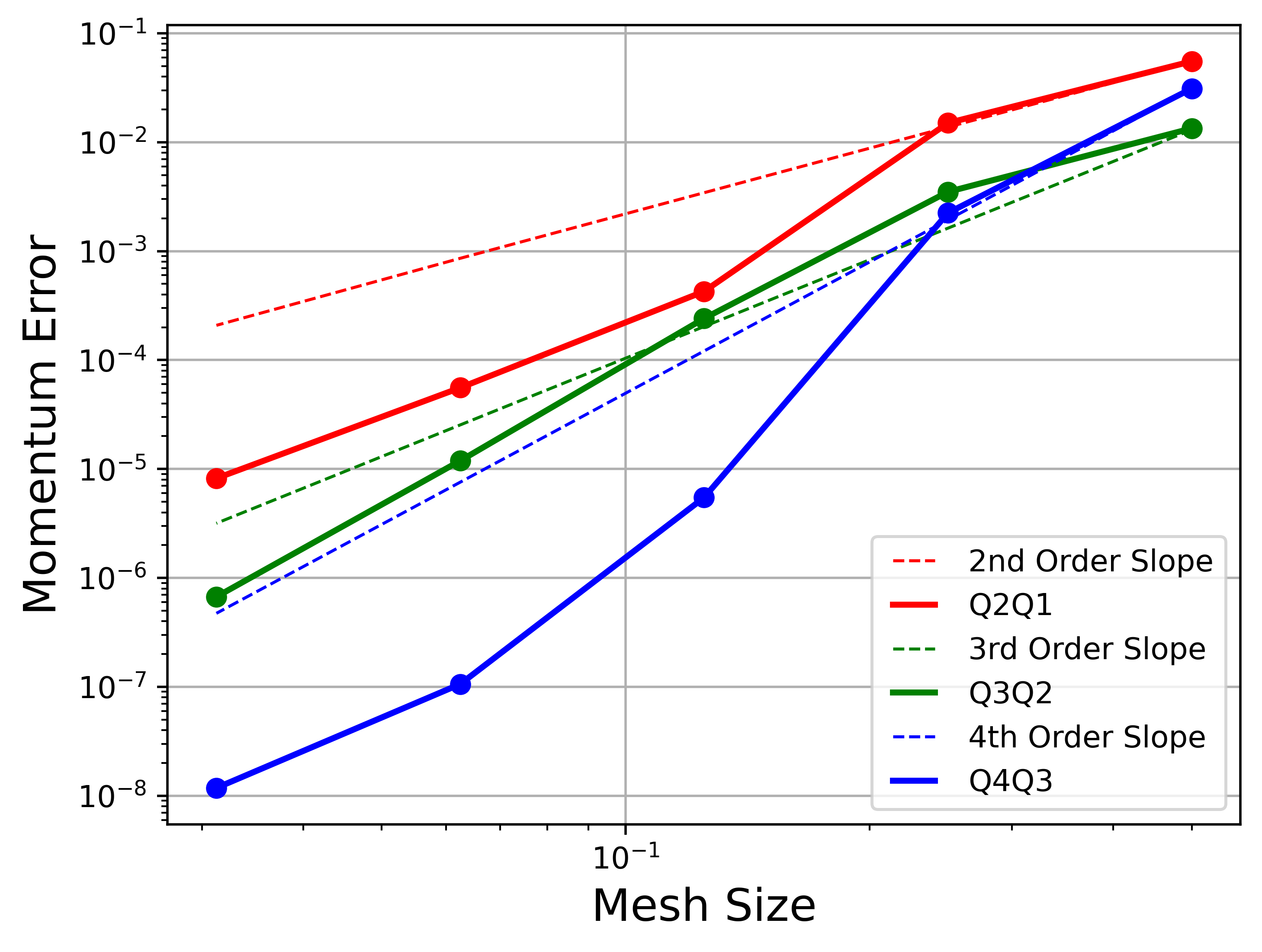} \hfil
  \includegraphics[width=0.3\textwidth]{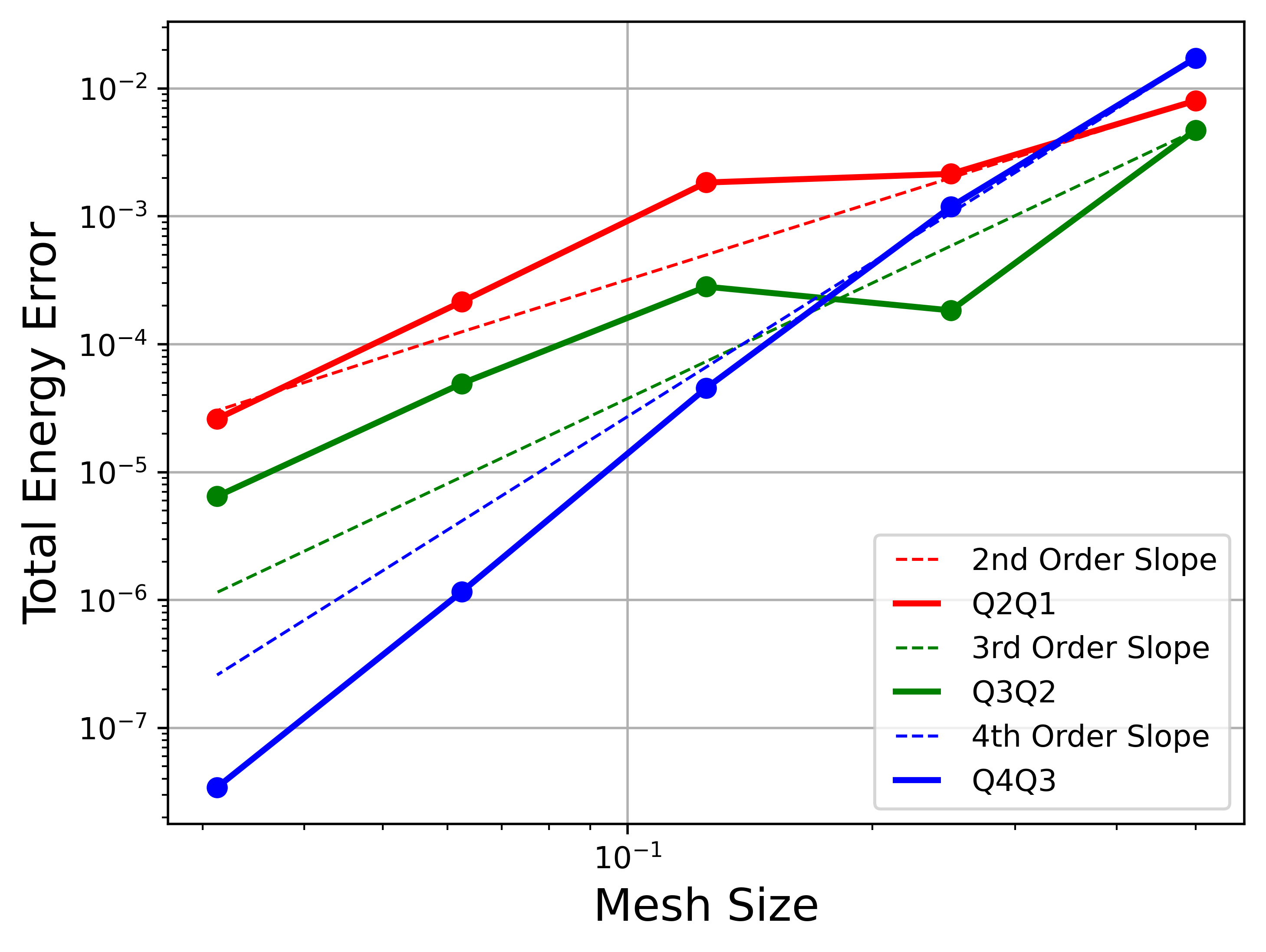}
}
\vspace{-3mm}
\caption{Convergence plots: L1 velocity errors, momentum errors,
         and total energy errors for the two-dimensional Taylor Green problem.}
\label{fig_tg_rates}
\end{figure}

For all simulations we use reconstruct 4th order pressure functions
and take 4 Taylor expansion terms in the shifted integrals.
The distances to the interface are recomputed at every time by the
normalization-based distance solver.
We test convergence for $Q_2-Q_1, Q_3-Q_2$, and $Q_4-Q_3$ discretizations,
which should produce 2nd, 3rd, and 4th order, respectively.
Table \ref{tab_tg} shows the L1 velocity errors there we observe the
expected convergence rates.
Figure \ref{fig_tg_rates} shows the convergence plot corresponding to
Table \ref{tab_tg}, and convergence plots of the total momentum and
total energy errors, which also converge as expected.

\begin{table}[h!]
\begin{center}
\begin{tabular}{c | c c | c c | c c}
\hline
~ref~ & ~Q2Q1 L1 error ~ & ~rate~
      & ~Q3Q2 L1 error ~ & ~rate~
      & ~Q4Q3 L1 error ~ & ~rate~  \\
\hline
0   & 7.93E-2 & -	   & 2.68E-2 & - 	  & 3.85E-2 & -    \\
1   & 4.67E-2 & 0.76 & 7.15E-3 & 1.91 & 3.08E-3 & 3.64 \\
2   & 7.50E-3 & 2.64 & 7.32E-4 & 3.29 & 2.30E-4 & 3.74 \\
3   & 1.44E-3 & 2.38 & 9.01E-5 & 3.02 & 9.84E-6 & 4.55 \\
4   & 2.61E-4 & 2.46 & 1.13E-5 & 2.99 & 3.91E-7 & 4.65 \\
\hline
\end{tabular}
\end{center}
\vspace{-3mm}
\caption{Convergence tests for the two-dimensional Taylor Green problem.}
\label{tab_tg}
\end{table}

\subsection{Two-material one-dimensional Sod shock tube}
\label{sec_sod}

Next we consider a two-material Riemann problem in one dimension.
This problem tests the Shifted Interface Method's ability
to equilibrate a pressure discontinuity, and allows to compare the behavior
of the method to case of a fitted interface.

The domain is $[0,1]$ with $v(0) = v(1) = 0$ boundary conditions.
Both materials are ideal gases with equation of state $p = (\gamma-1) \rho e$.
The interface is at $x_c = 0.5$ for the fitted simulations and is moved
to $x_c = 0.5 + 0.5 \Delta x$ for the shifted interface simulations.
The problem is run to a final time of $t = 0.2$.
The two initial states are:
\[
(v, \rho, e, p, \gamma)=
  \begin{cases}
     (0, 1, 2, 2, 2)         & \text{if } x < x_c ~~ \text{(Left material)}, \\
     (0, 0.125, 2, 0.1, 1.4) & \text{if } x < x_c ~~ \text{(Right material)}. \\
  \end{cases}
\]

\begin{figure}[h!]
\centerline
{
  \includegraphics[width=0.3\textwidth]{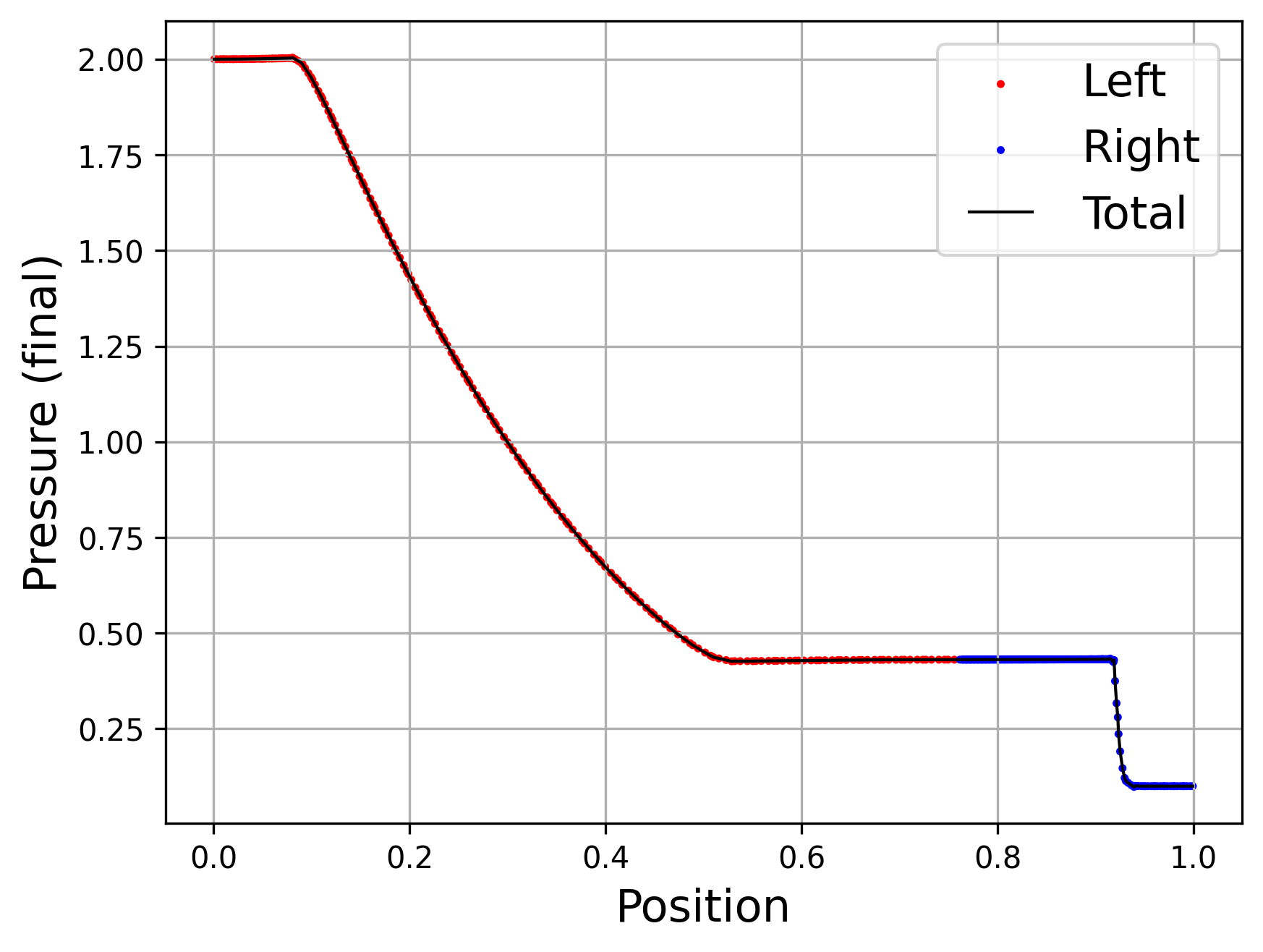} \hfil
  \includegraphics[width=0.3\textwidth]{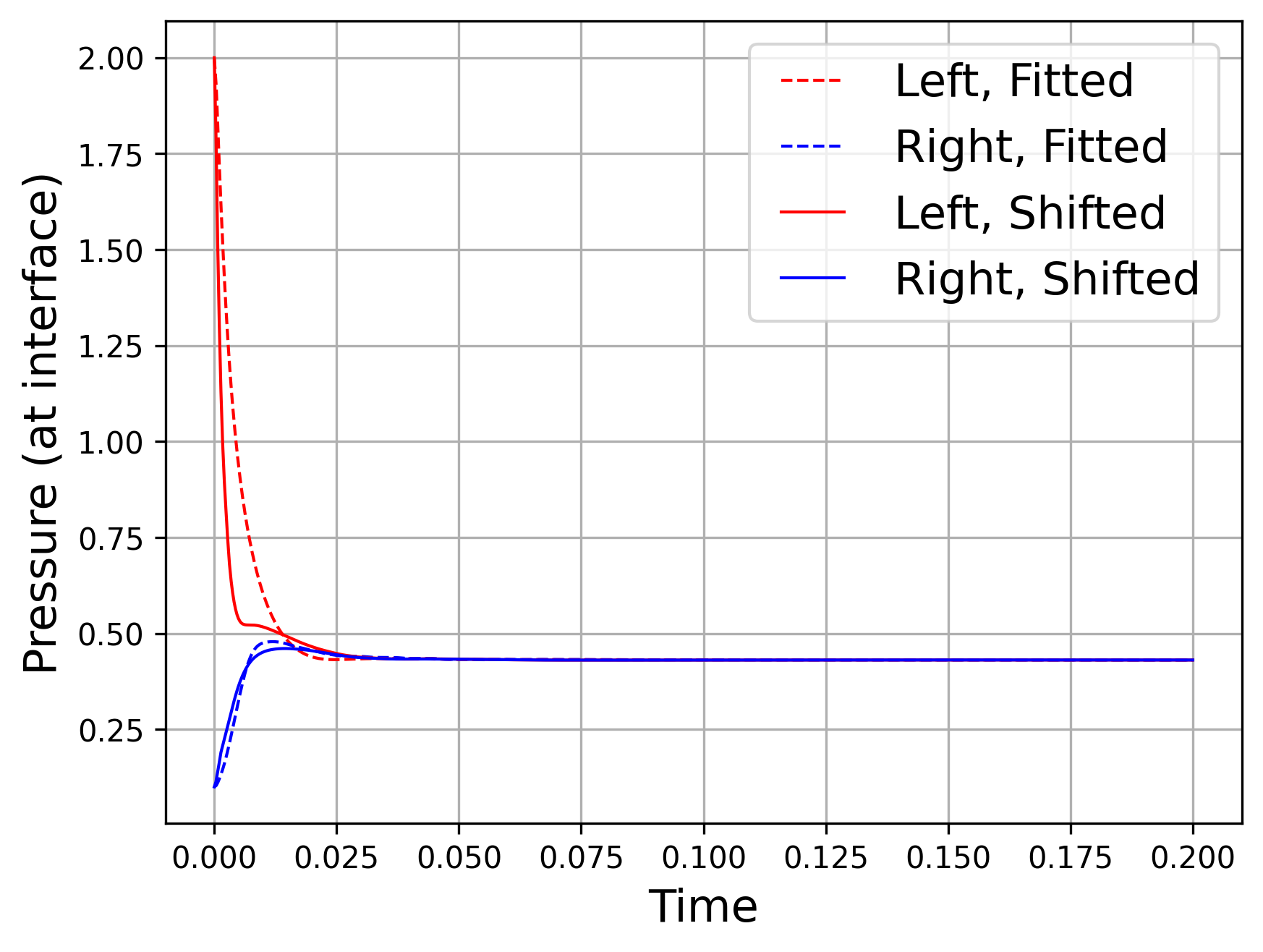}  \hfil
  \includegraphics[width=0.3\textwidth]{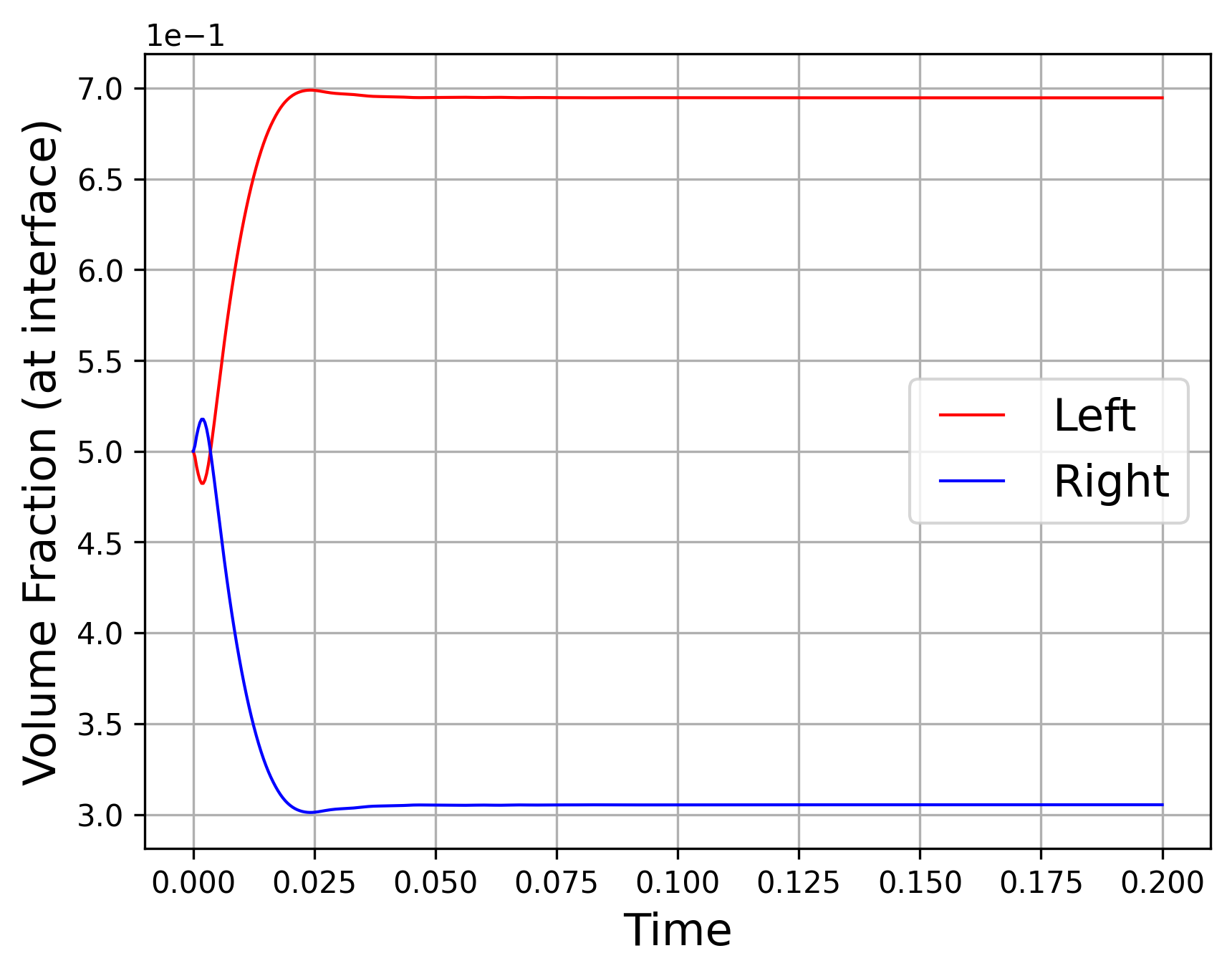}
}
\caption{Simulation with \textit{active} shifted interface terms:
         final material pressures (left),
         pressures' time-history at the interface (middle),
         volume fractions' time-history at the interface (right),
         for the one-dimensional Sod tube test.}
\label{fig_sod_shift}
\end{figure}

\begin{figure}[h!]
\centerline
{
  \includegraphics[width=0.3\textwidth]{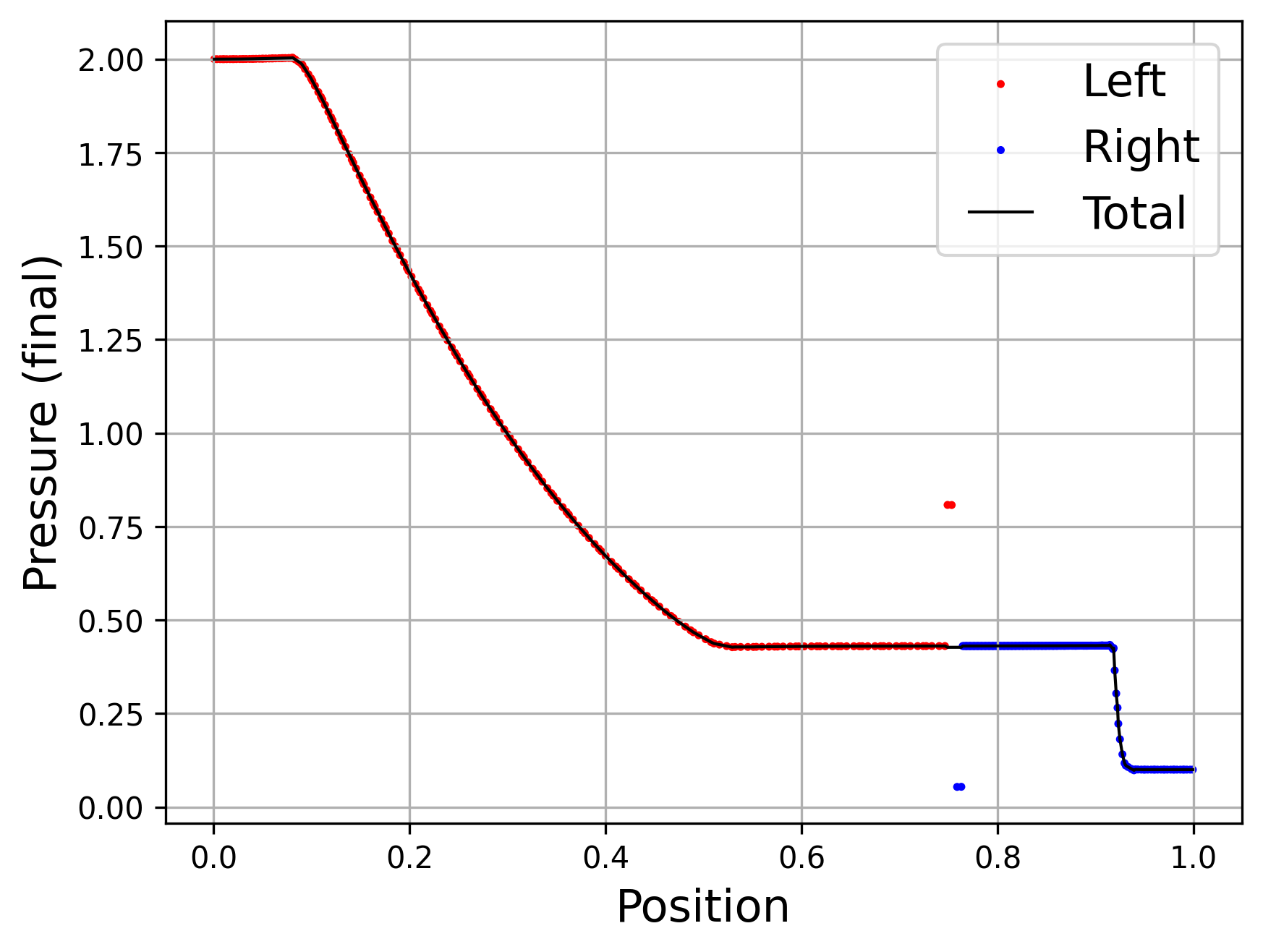} \hfil
  \includegraphics[width=0.3\textwidth]{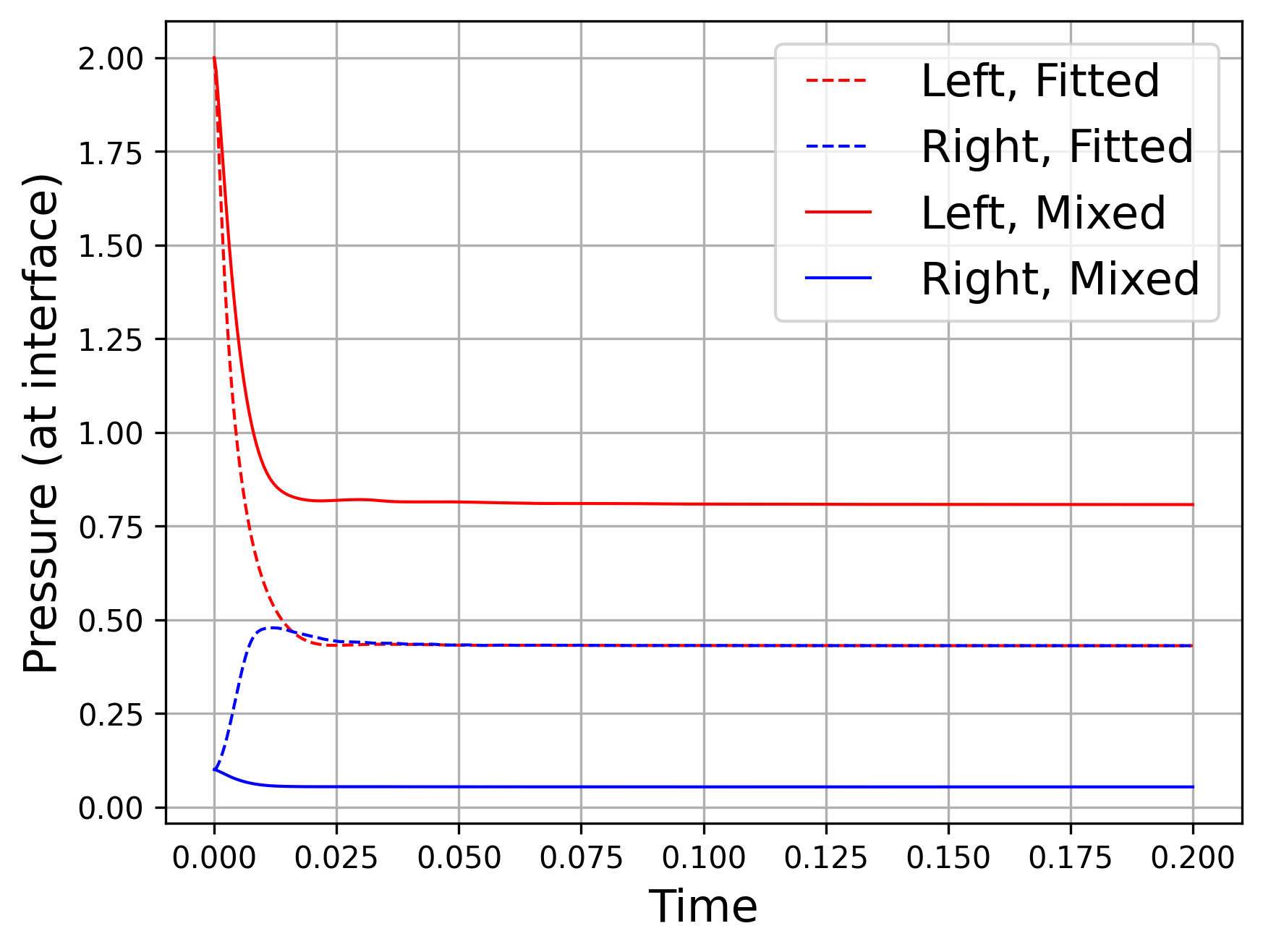}  \hfil
  \includegraphics[width=0.3\textwidth]{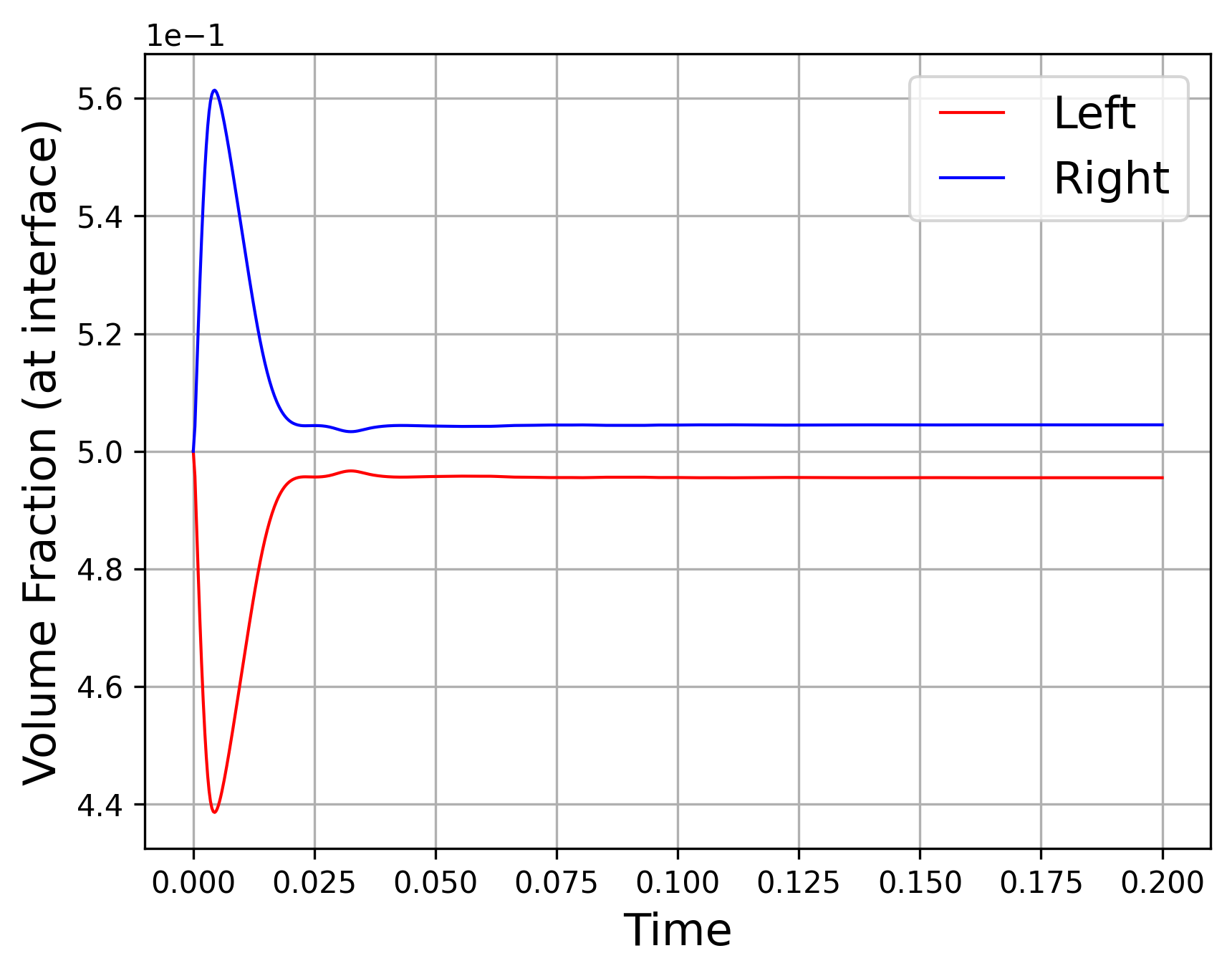}
}
\caption{Simulation with \textit{inactive} shifted interface terms:
         final material pressures (left),
         pressures' time-history at the interface (middle),
         volume fractions' time-history at the interface (right),
         for the one-dimensional Sod tube test.}
\label{fig_sod_mix}
\end{figure}

\begin{figure}[h!]
\centerline
{
  \includegraphics[width=0.3\textwidth]{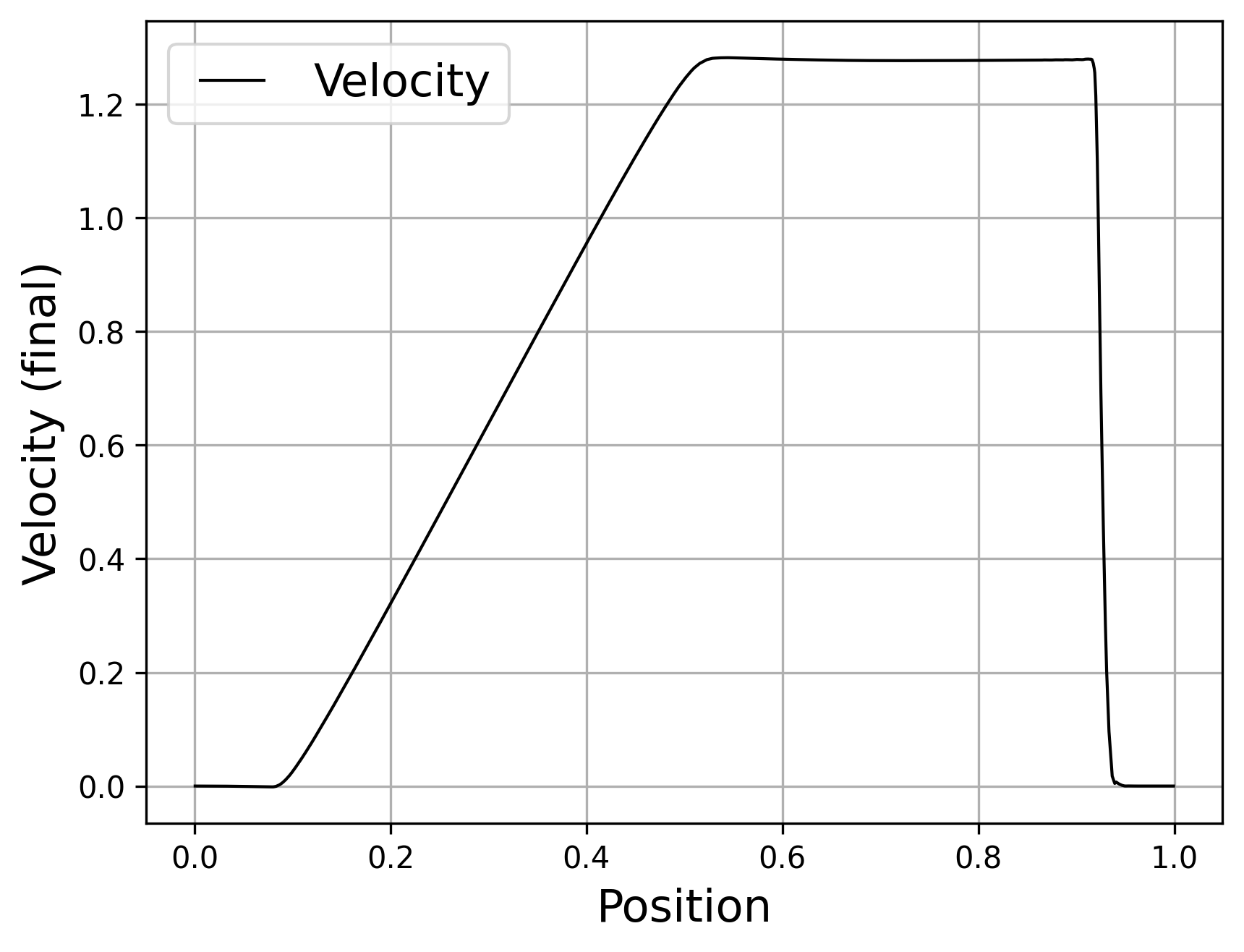}   \hfil
  \includegraphics[width=0.3\textwidth]{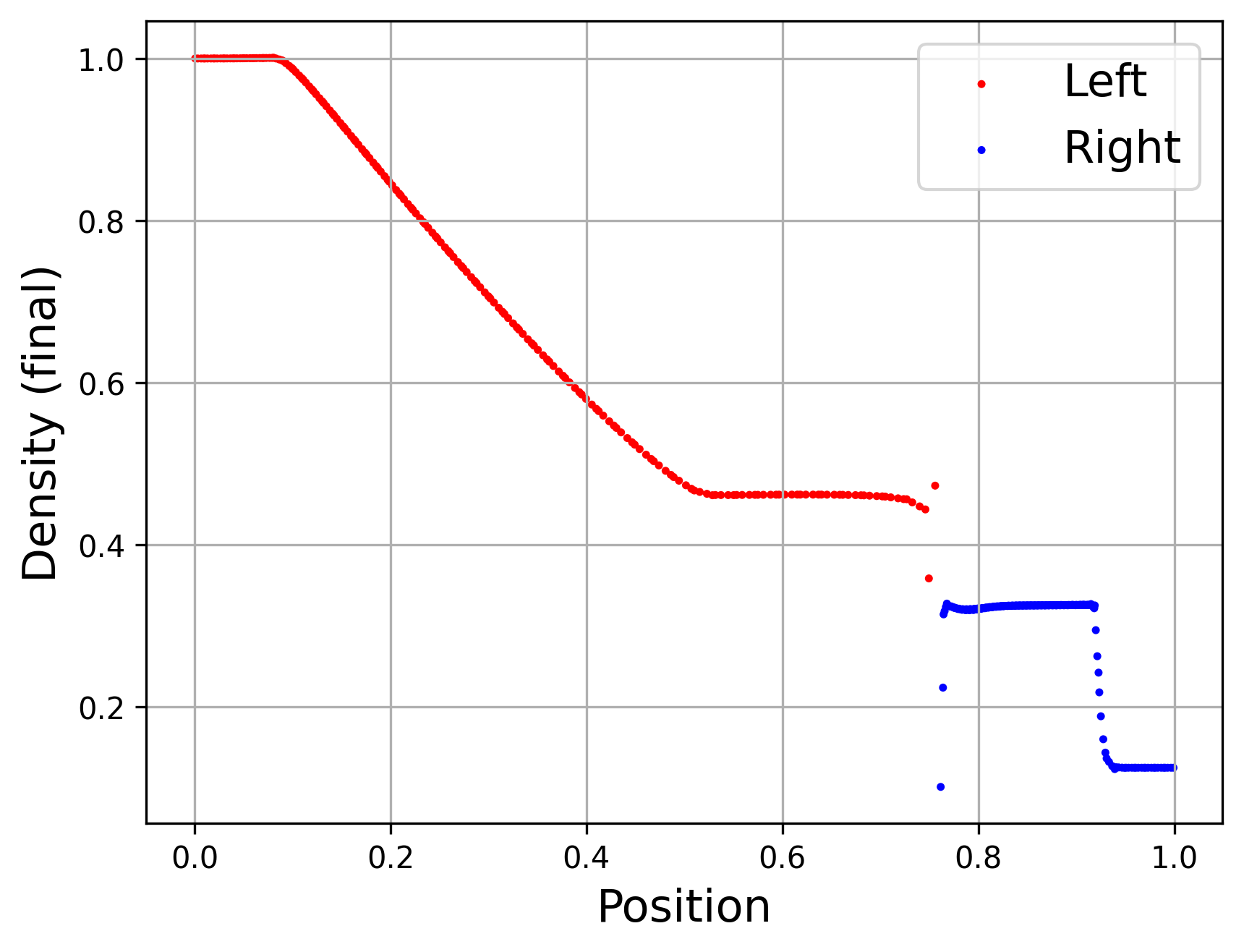} \hfil
  \includegraphics[width=0.3\textwidth]{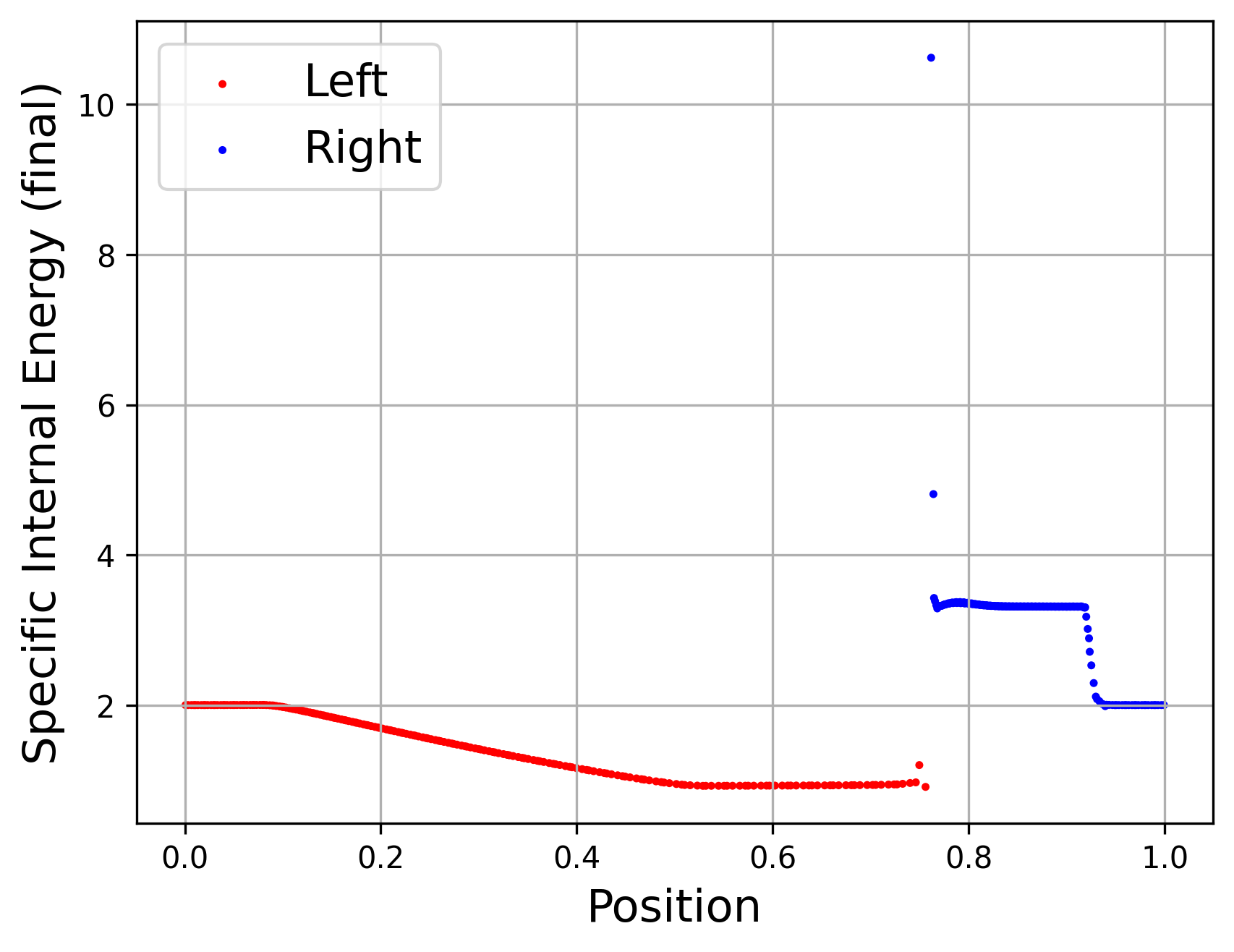}
}
\caption{Final velocity, material densities, and material specific internal
         energies for the one-dimensional Sod tube test.}
\label{fig_sod_vre}
\end{figure}

The problem is run on 100 elements with a $Q_2-Q_1$ discretization.
To stress the effects of the shifted face integrals, we compare two sets of
results, namely, (i) the results obtained with the presented method, and
(ii) results obtained by the same simulation without including the shifted
face integrals.
Results from case (i) are in Figure \ref{fig_sod_shift}, showing
the final material pressures, their evolution in time at the interface, and
the time evolution of the material volume fractions.
Visualized point values (e.g., on the left panel of Figure \ref{fig_sod_shift})
correspond to the quadrature points that are used
to compute all volumetric integrals (4 points per element in this case).
For each material, a quantity is visualized at a given point $x$ whenever
$\eta_0(x_0) \geq 0$, where $x_0$ is the initial position of $x$.
We observe that the pressures equilibrate similarly to the fitted case, while
the volume of the left material expands in the mixed element, as expected.
The corresponding results of case (ii) are in Figure \ref{fig_sod_mix}.
We observe that the pressures do not equilibrate, and the materials do not expand
(left material) and compress (right material) correctly.
The final velocity, material densities, and material internal energies
for case (i) are shown in Figure \ref{fig_sod_vre}.
Although the pressures are in equilibrium, oscillations in the material
densities and specific internal energies are present in the mixed element.
This is a known issue for Lagrangian methods that do not add artificial
diffusion around the contact discontinuity.
Note that the shifted interface terms are also unable to detect such
oscillations, as the terms use information only about the gradients of
pressure and velocity, which are zero at the equilibrated contact.

\subsection{1D Water-Air Shock Tube}
\label{sec_wa}

The next test is a more severe one dimensional Riemann problem.
It represents a water-air interaction with a large initial pressure jump.
This problem tests the robustness of the Shifted Interface Method under
extreme conditions.

The domain is $[0,1]$ with $v(0) = v(1) = 0$ boundary conditions.
The high-pressure water (left side) material is represented by the
stiffened gas equation of state $p = (\gamma-1) \rho e - \gamma 6 \times 10^8$.
The low-pressure air (right side) material represented by ideal
gas with equation of state $p = (\gamma-1) \rho e$.
The interface is at $x_c = 0.7$ for the fitted simulations and is moved
to $x_c = 0.7 + 0.5 \Delta x$ for the shifted interface simulations.
The two initial states are:
\[
(v, \rho, p, \gamma)=
  \begin{cases}
     (0, 10^3, 10^9, 4.4)  & \text{if } x < x_c ~~ \text{(Left material)}, \\
     (0, 50,   10^5, 1.4)  & \text{if } x < x_c ~~ \text{(Right material)}. \\
  \end{cases}
\]

\begin{figure}[h!]
\centerline
{
  \includegraphics[width=0.3\textwidth]{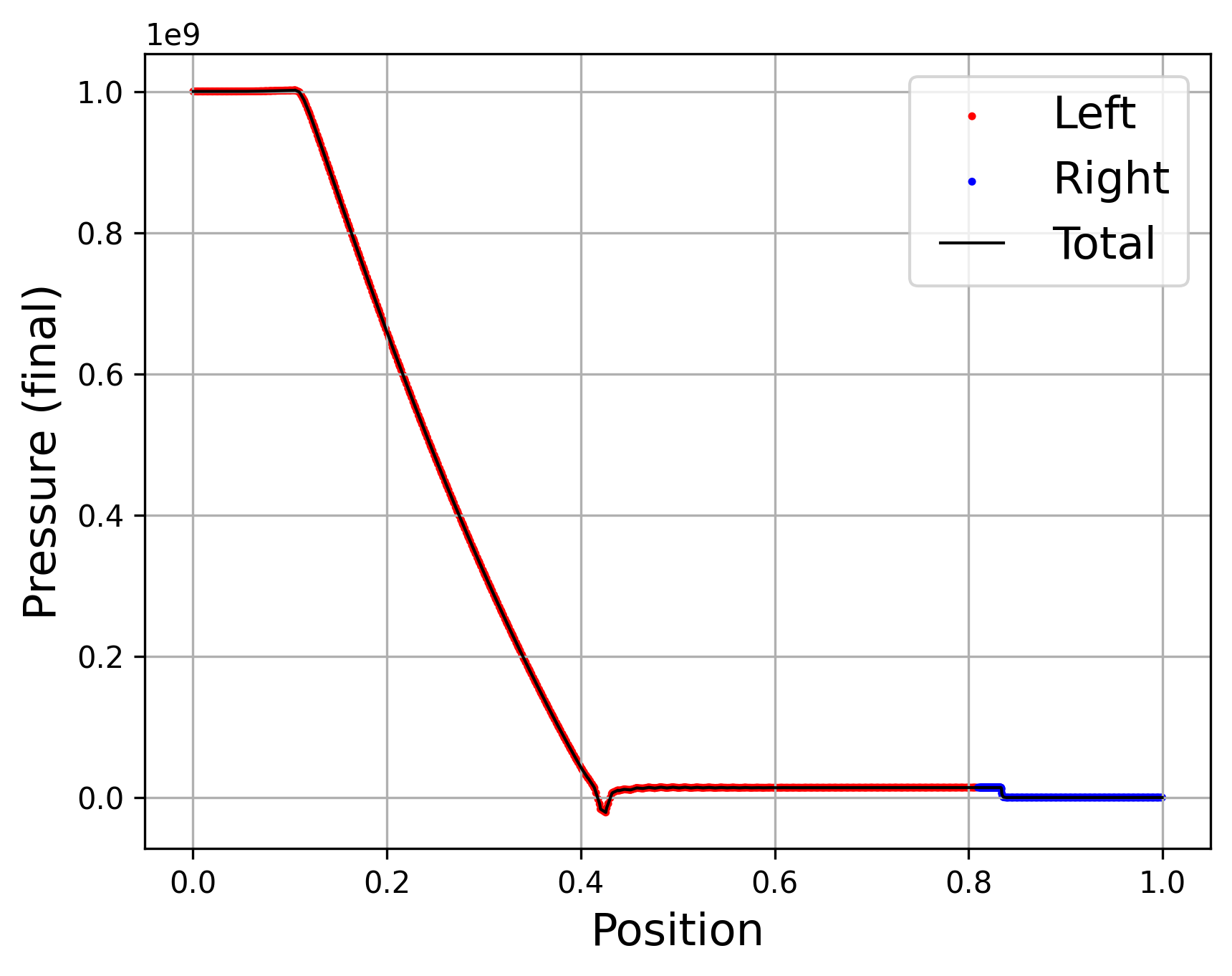} \hfil
  \includegraphics[width=0.3\textwidth]{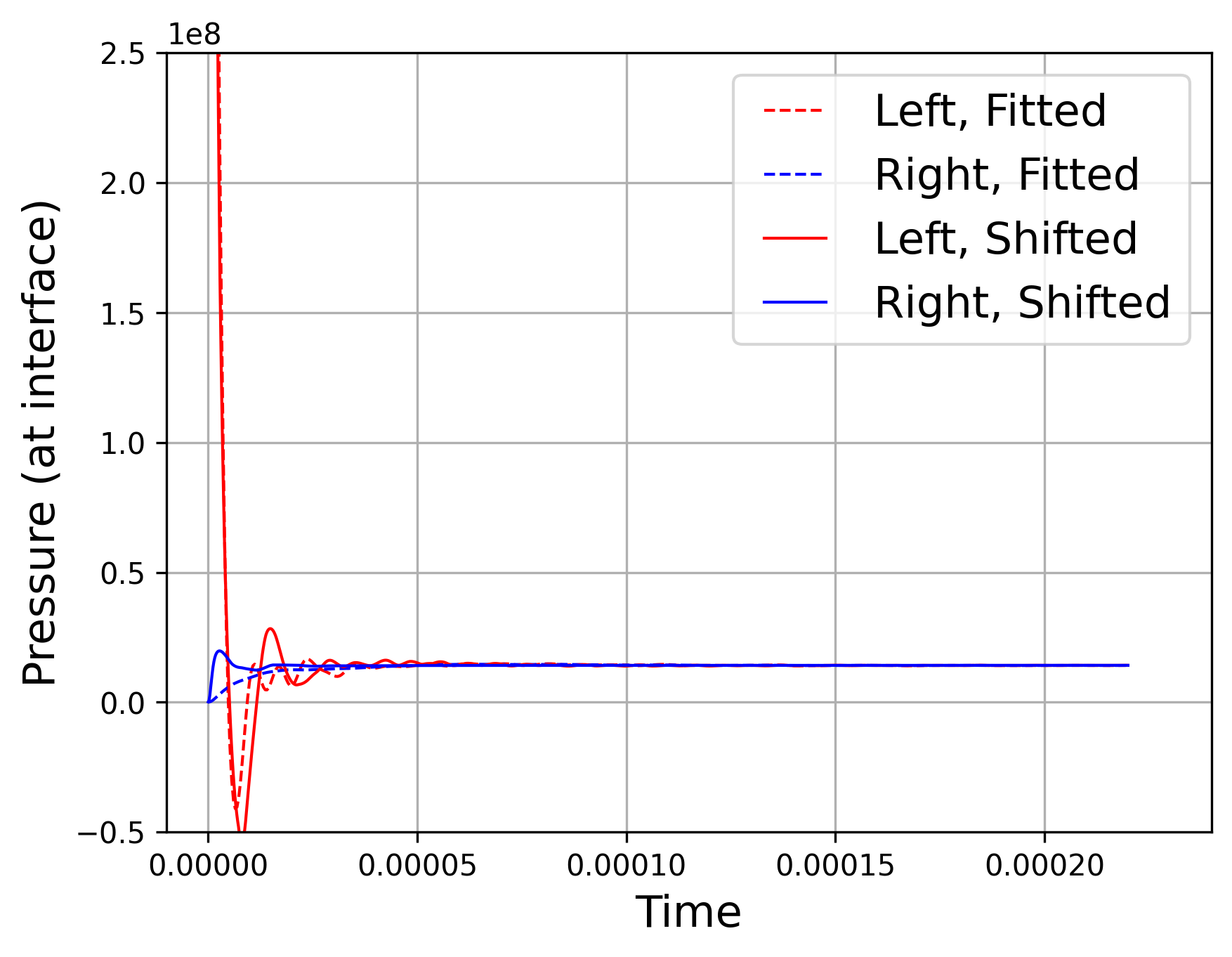}  \hfil
  \includegraphics[width=0.3\textwidth]{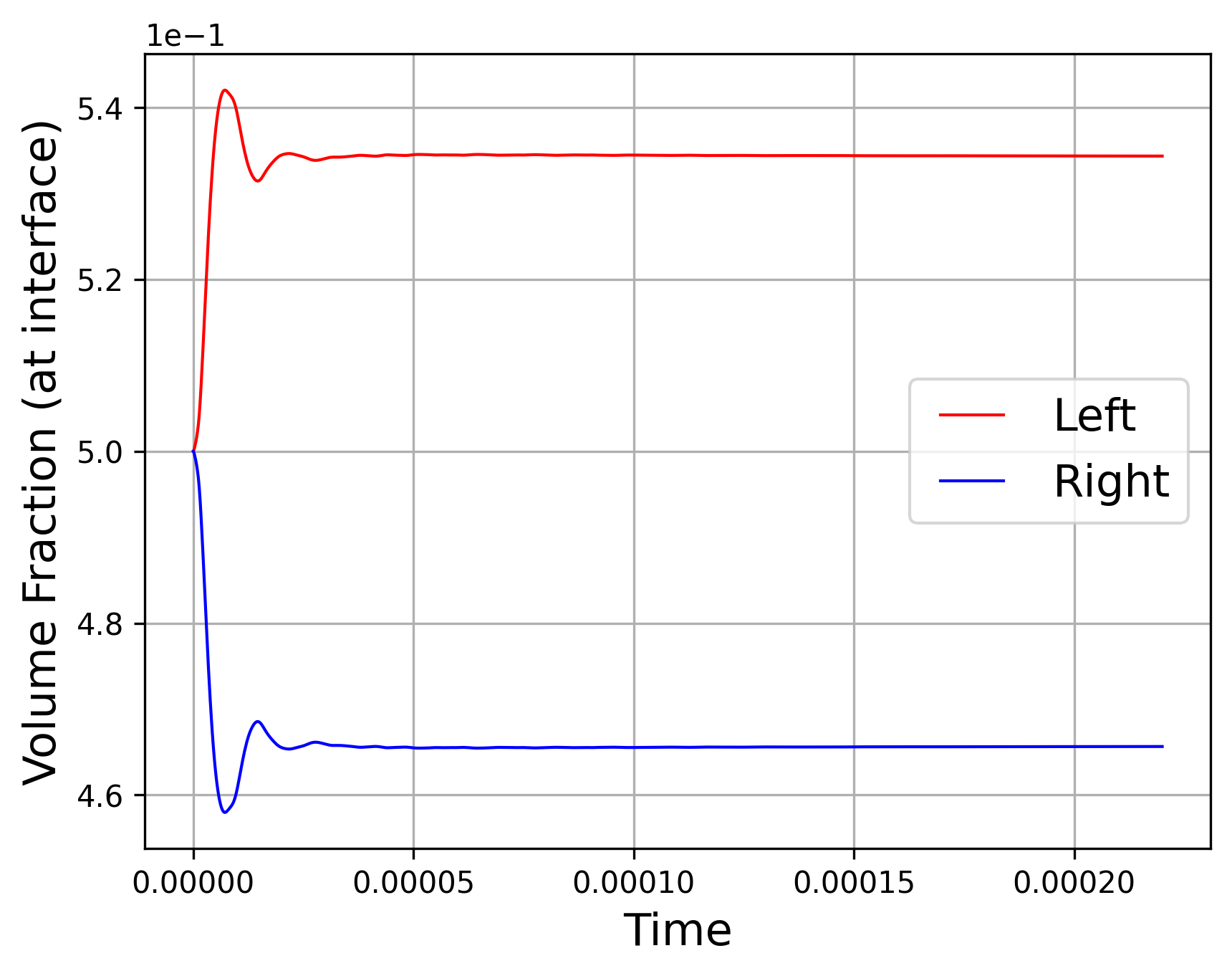}
}
\caption{Simulation with \textit{active} shifted interface terms:
         final material pressures (left),
         pressures' time-history at the interface (middle),
         volume fractions' time-history at the interface (right),
         for the one-dimensional Water-Air shock tube test.}
\label{fig_wa_shift}
\end{figure}

\begin{figure}[h!]
\centerline
{
  \includegraphics[width=0.3\textwidth]{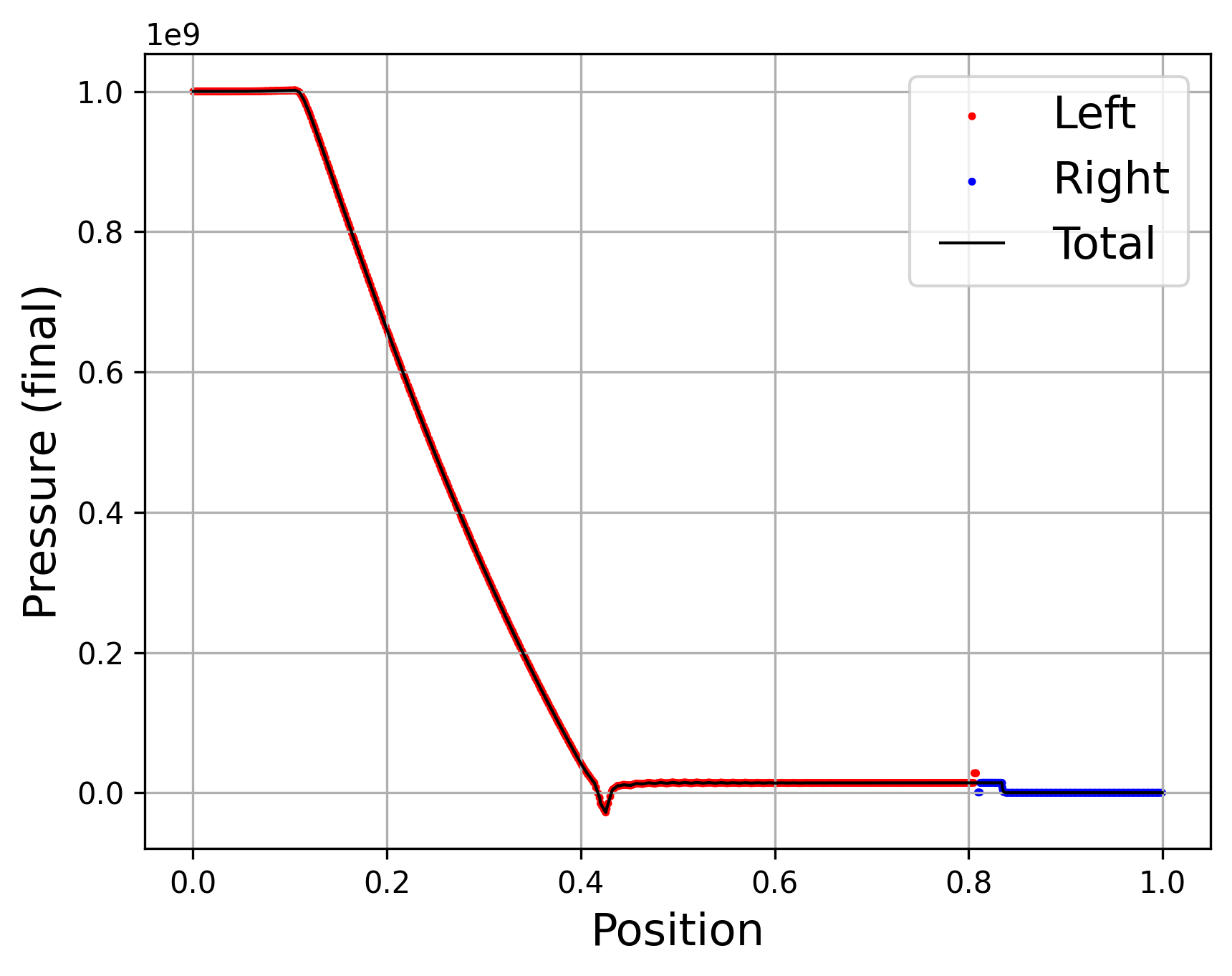} \hfil
  \includegraphics[width=0.3\textwidth]{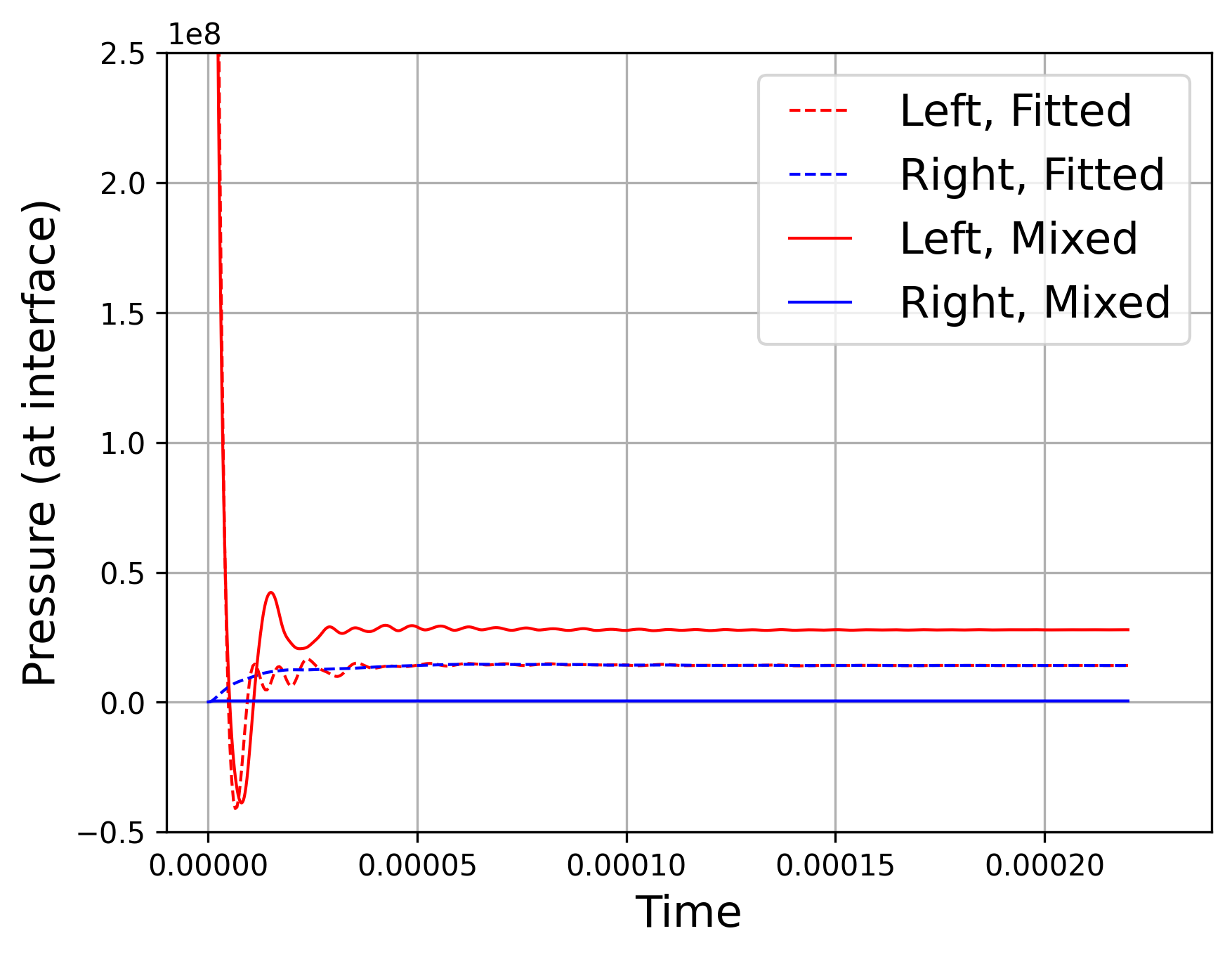}  \hfil
  \includegraphics[width=0.3\textwidth]{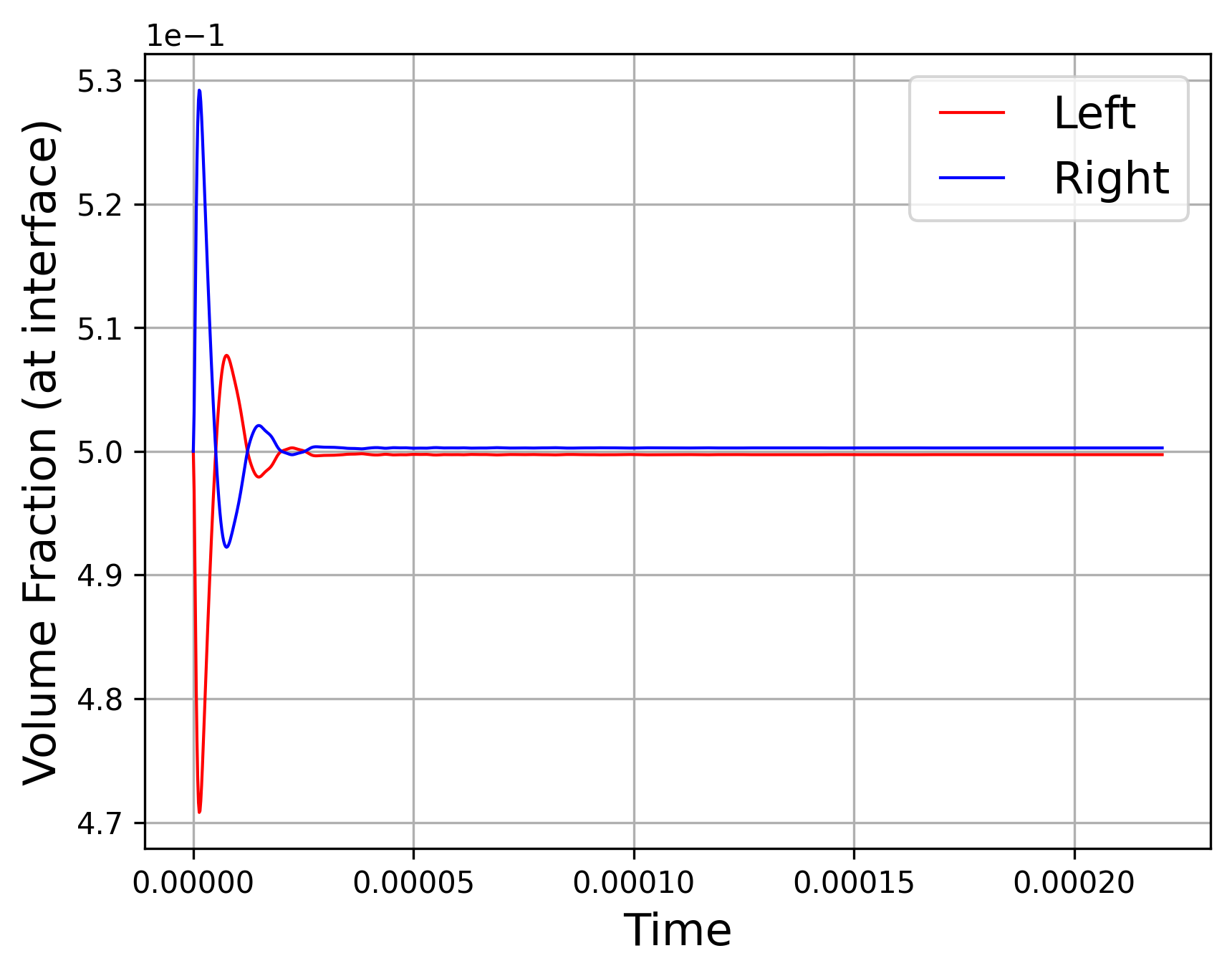}
}
\caption{Simulation with \textit{inactive} shifted interface terms:
         final material pressures (left),
         pressures' time-history at the interface (middle),
         volume fractions' time-history at the interface (right),
         for the one-dimensional Water-Air shock tube test.}
\label{fig_wa_mix}
\end{figure}

\begin{figure}[h!]
\centerline
{
  \includegraphics[width=0.3\textwidth]{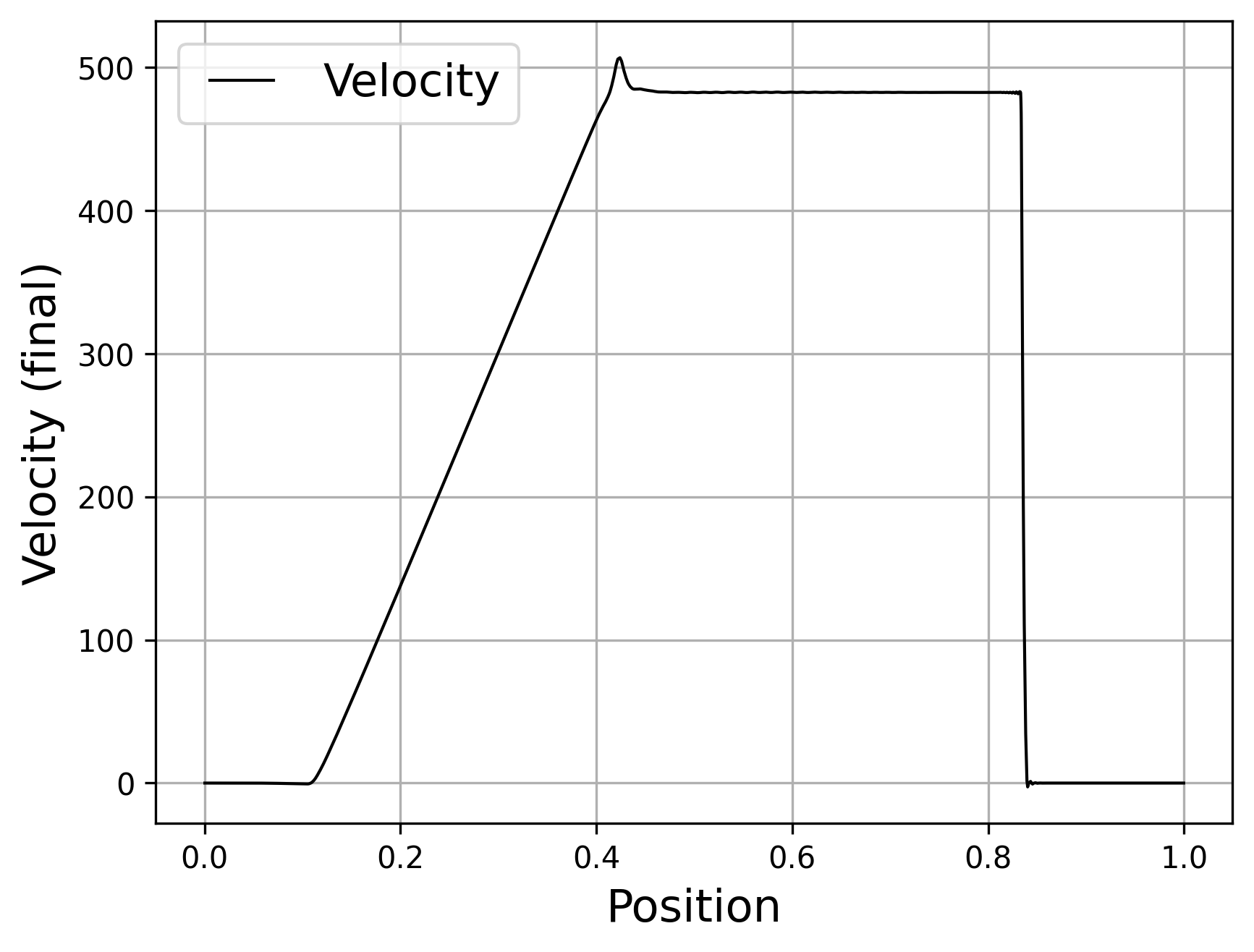}   \hfil
  \includegraphics[width=0.3\textwidth]{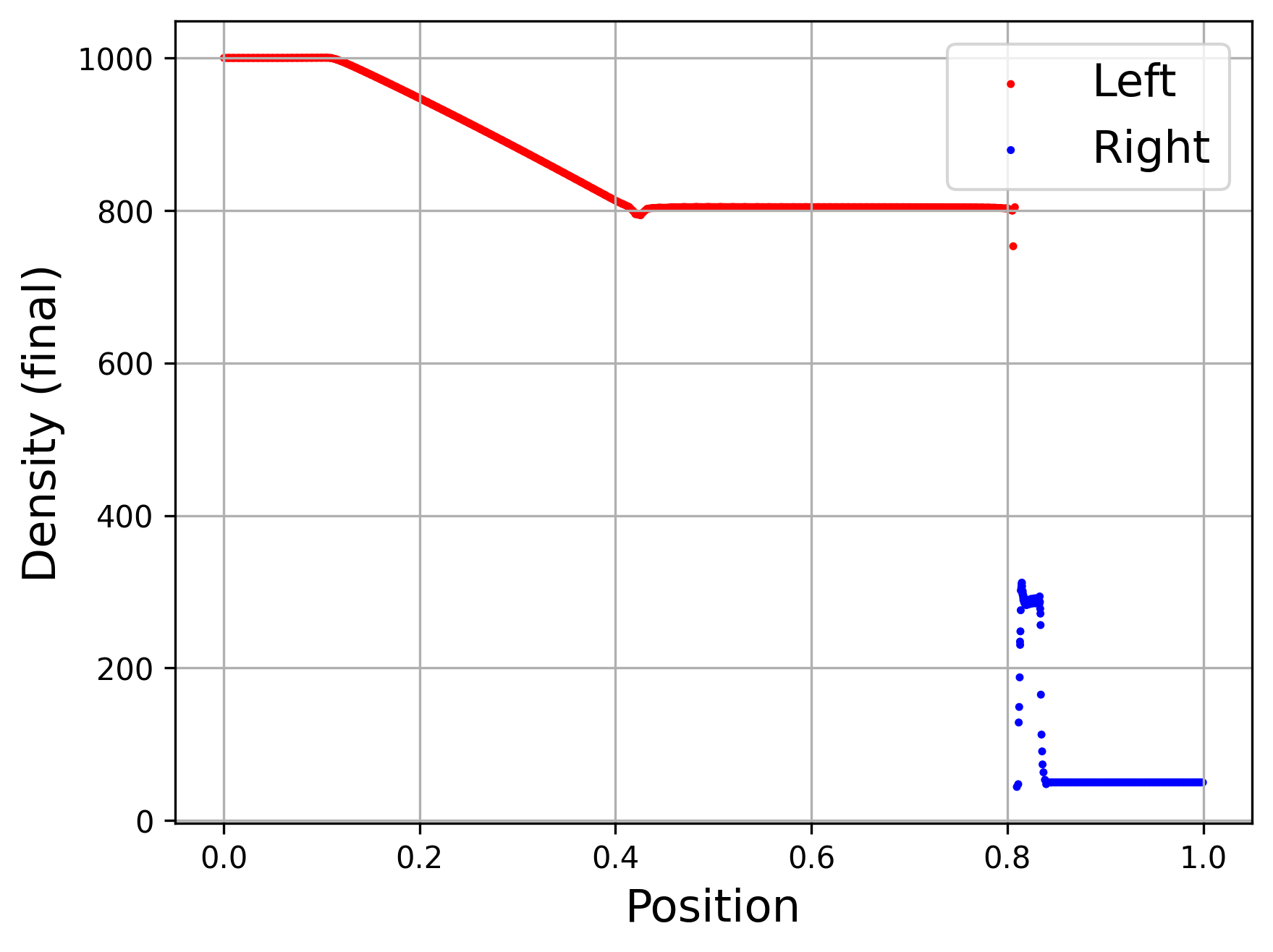} \hfil
  \includegraphics[width=0.3\textwidth]{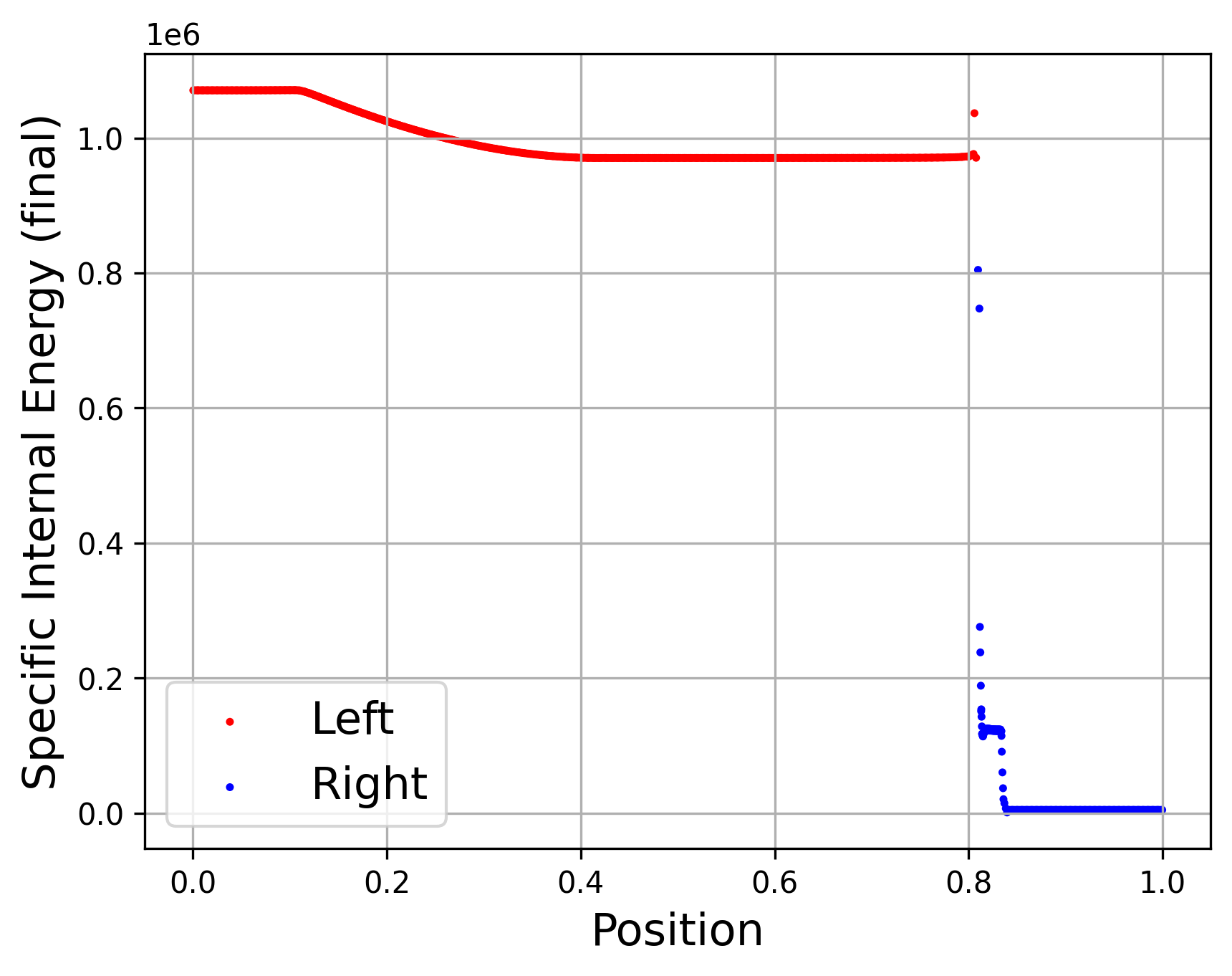}
}
\caption{Final velocity, material densities, and material specific internal
         energies for the Water-Air shock tube test.}
\label{fig_wa_vre}
\end{figure}

The problem is run on 200 elements with a $Q_2-Q_1$ discretization.
Again we compare (i) the results obtained with the presented method, and
(ii) results obtained by the same simulation without including the shifted
face integrals.
Results from case (i) are in Figure \ref{fig_wa_shift}, showing
the final material pressures, their evolution in time at the interface, and
the time evolution of the material volume fractions
(Section \ref{sec_sod} explains the visualization of the point values).
The corresponding results of case (ii) are in Figure \ref{fig_wa_mix}.
We make similar observations as in the Sod test, namely, the shifted terms
are able to equilibrate the material pressures, and expand the left material
near the interface.
The final velocity, material densities, and material internal energies
for case (i) are shown in Figure \ref{fig_wa_vre}.
Again we observe some oscillations of the material-specific
densities and specific internal energies around the cut element, due to the
lack of artificial viscosity around the contact region.
The more concerning behavior, however, occurs at the end of the rarefaction
region of the left material, where the final pressure is negative and the
corresponding velocity has a pronounced oscillation.
We do not fully understand this behavior and
we are exploring ways to improve it.

\subsection{Two-dimensional and three-dimensional Triple Point Interaction.}
\label{sec_3p}

This example demonstrates the main advantage of the method, namely, the
ability to maintain sharp curved material interfaces inside curved zones
for complex two-dimensional and three-dimensional problems.
The test is a modification of the standard Triple Point problem from
\cite{Dobrev2012}, Section 8.5.

\begin{figure}[h!]
\centerline
{
   \includegraphics[width=0.3\textwidth]{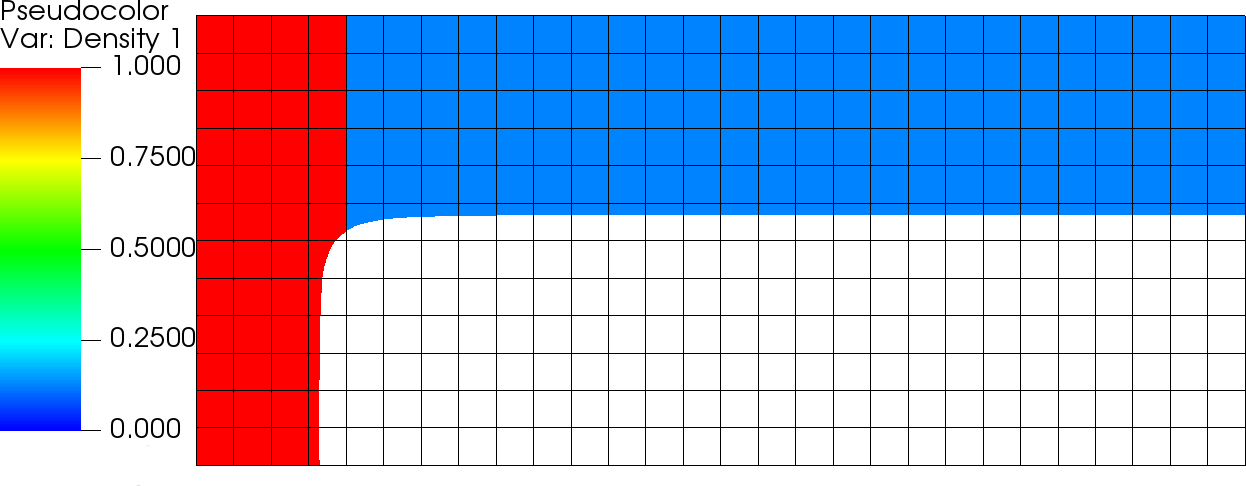}\hfil
   \includegraphics[width=0.3\textwidth]{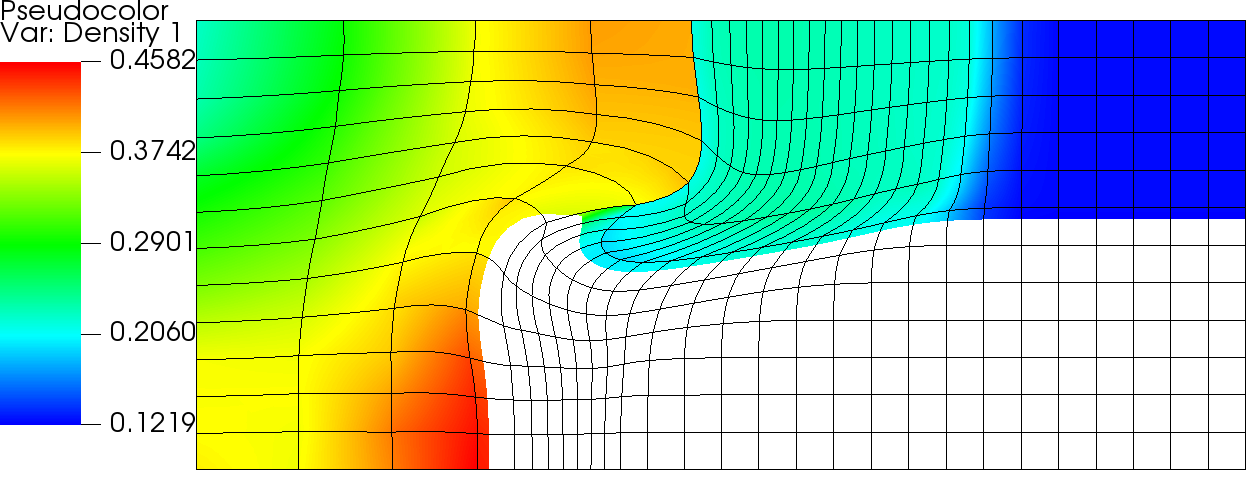}\hfil
   \includegraphics[width=0.3\textwidth]{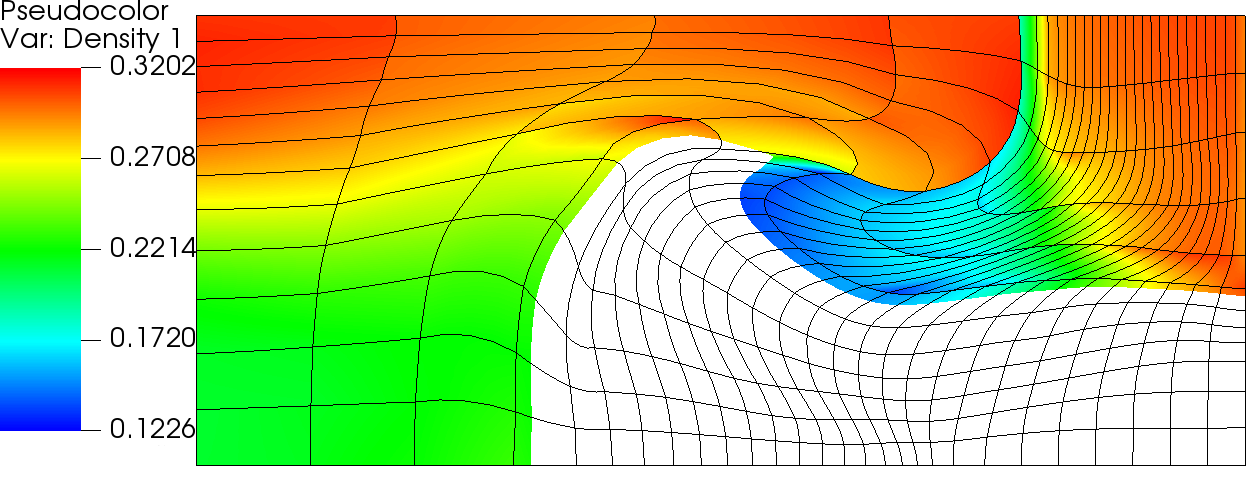}
}
\centerline
{
   \includegraphics[width=0.3\textwidth]{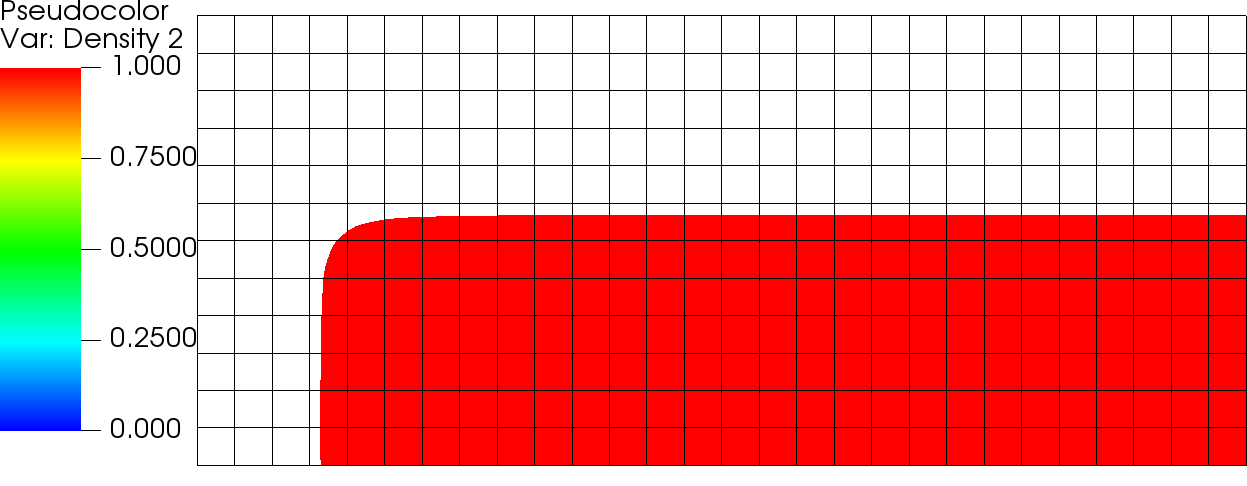}\hfil
   \includegraphics[width=0.3\textwidth]{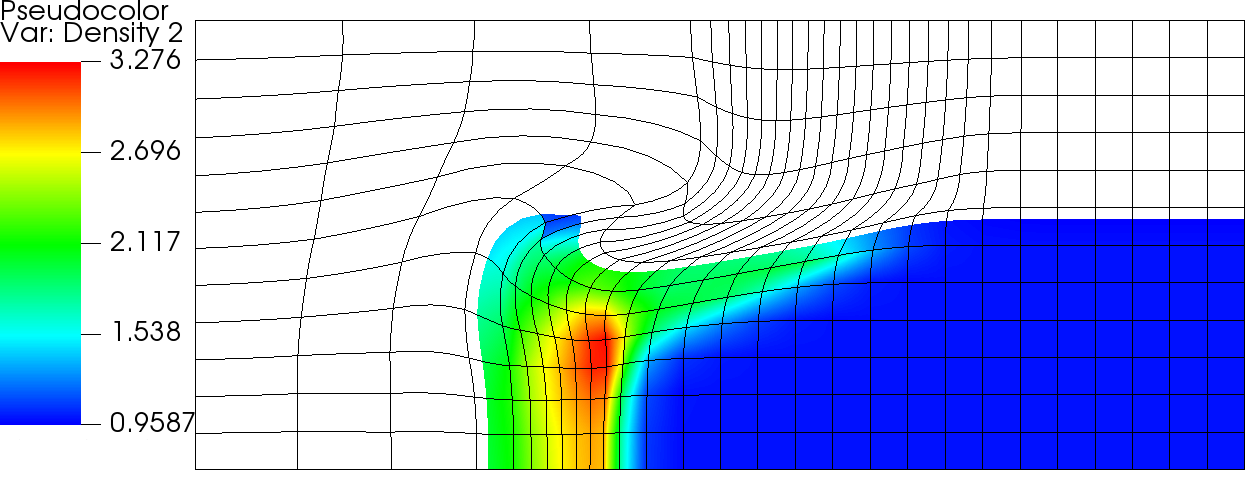}\hfil
   \includegraphics[width=0.3\textwidth]{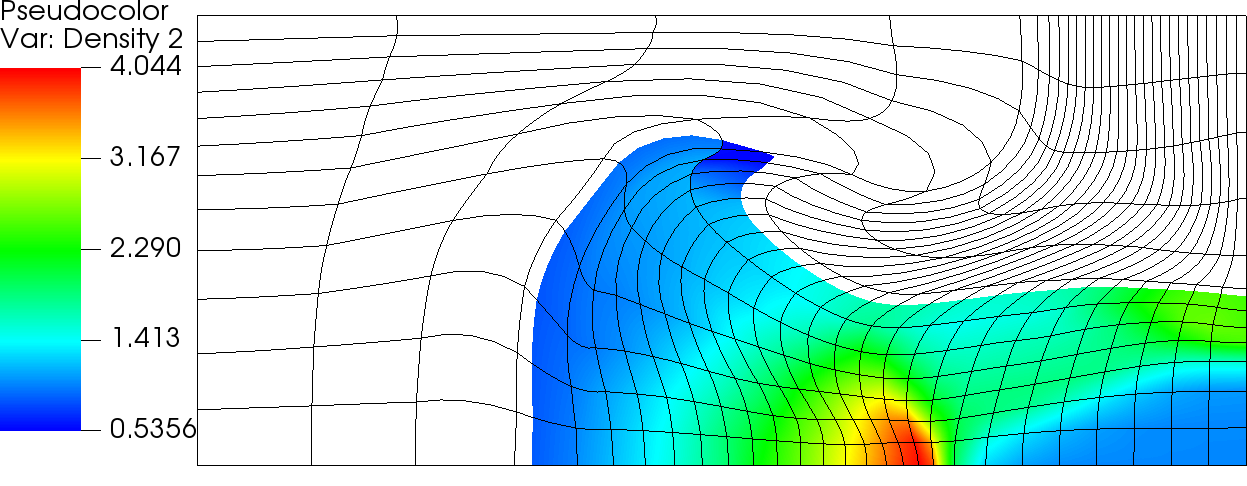}
}
\caption{Material densities (\#1 on top, \#2 on bottom)
         at times $t=0$, $t=2.5$, and $t=5$
         for the two-dimensional Triple Point interaction problem.}
\label{fig_3p_2D}
\end{figure}

\begin{figure}[h!]
\centerline
{
   \includegraphics[width=0.3\textwidth]{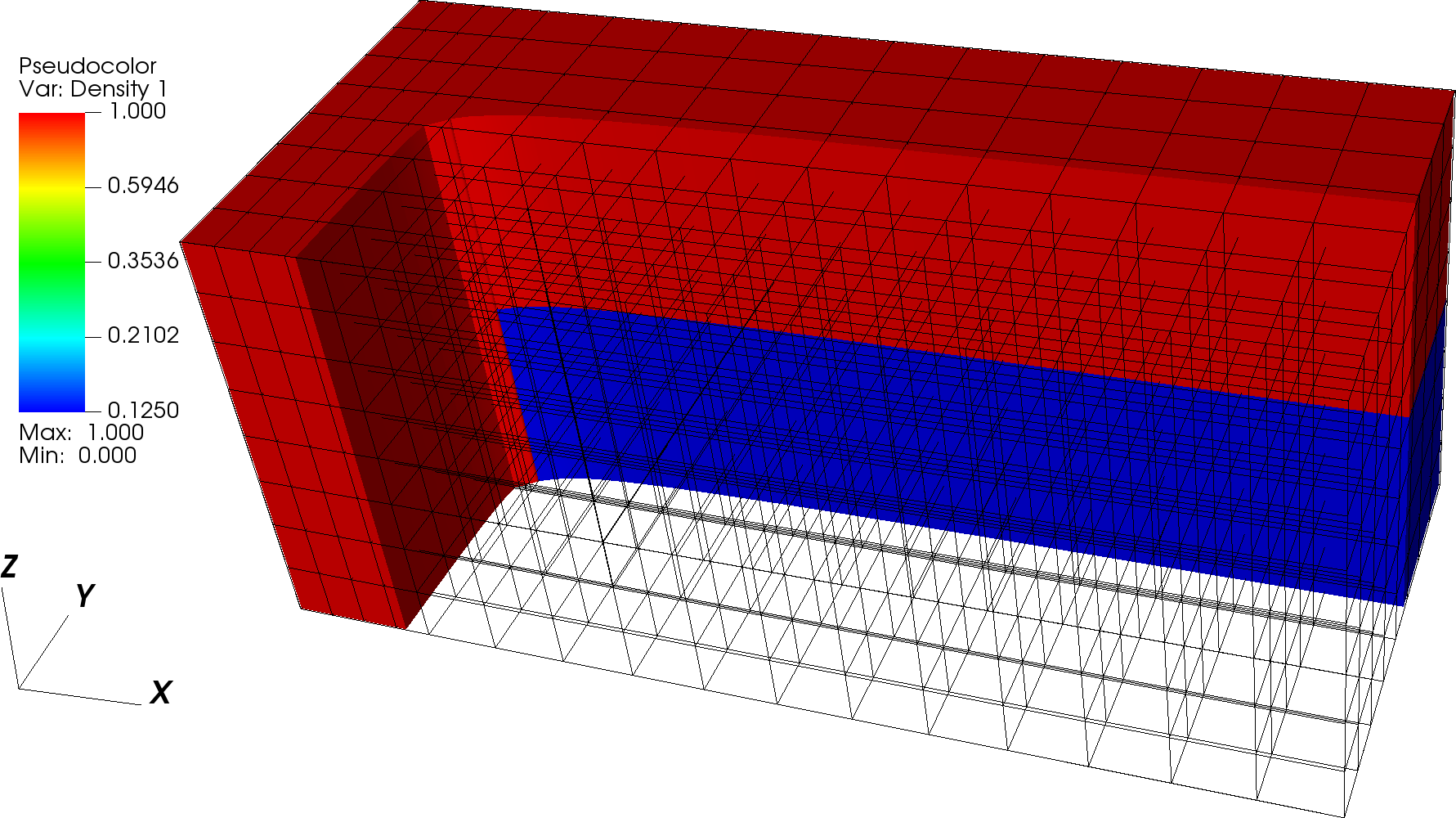}\hfil
   \includegraphics[width=0.3\textwidth]{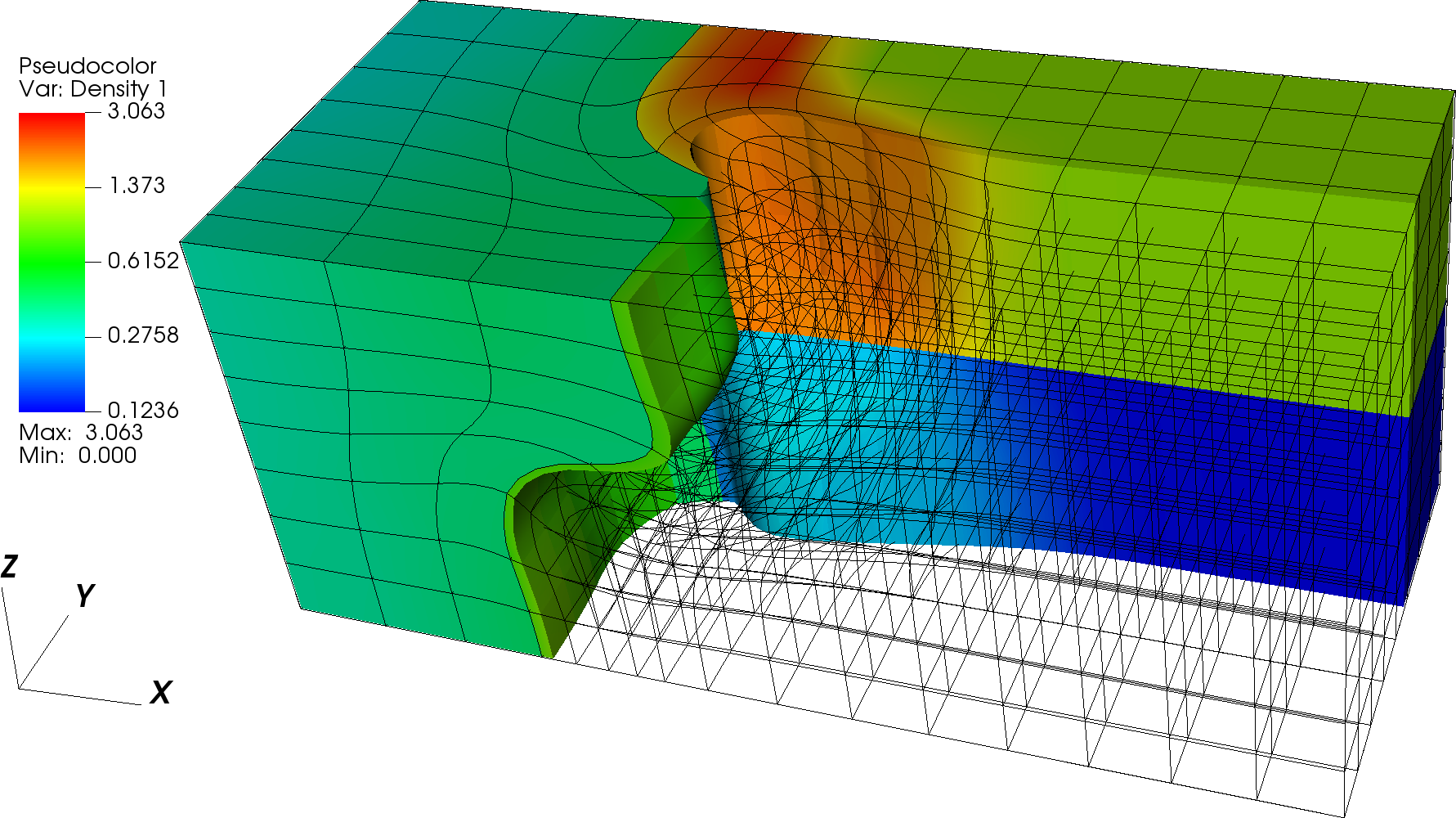}\hfil
   \includegraphics[width=0.3\textwidth]{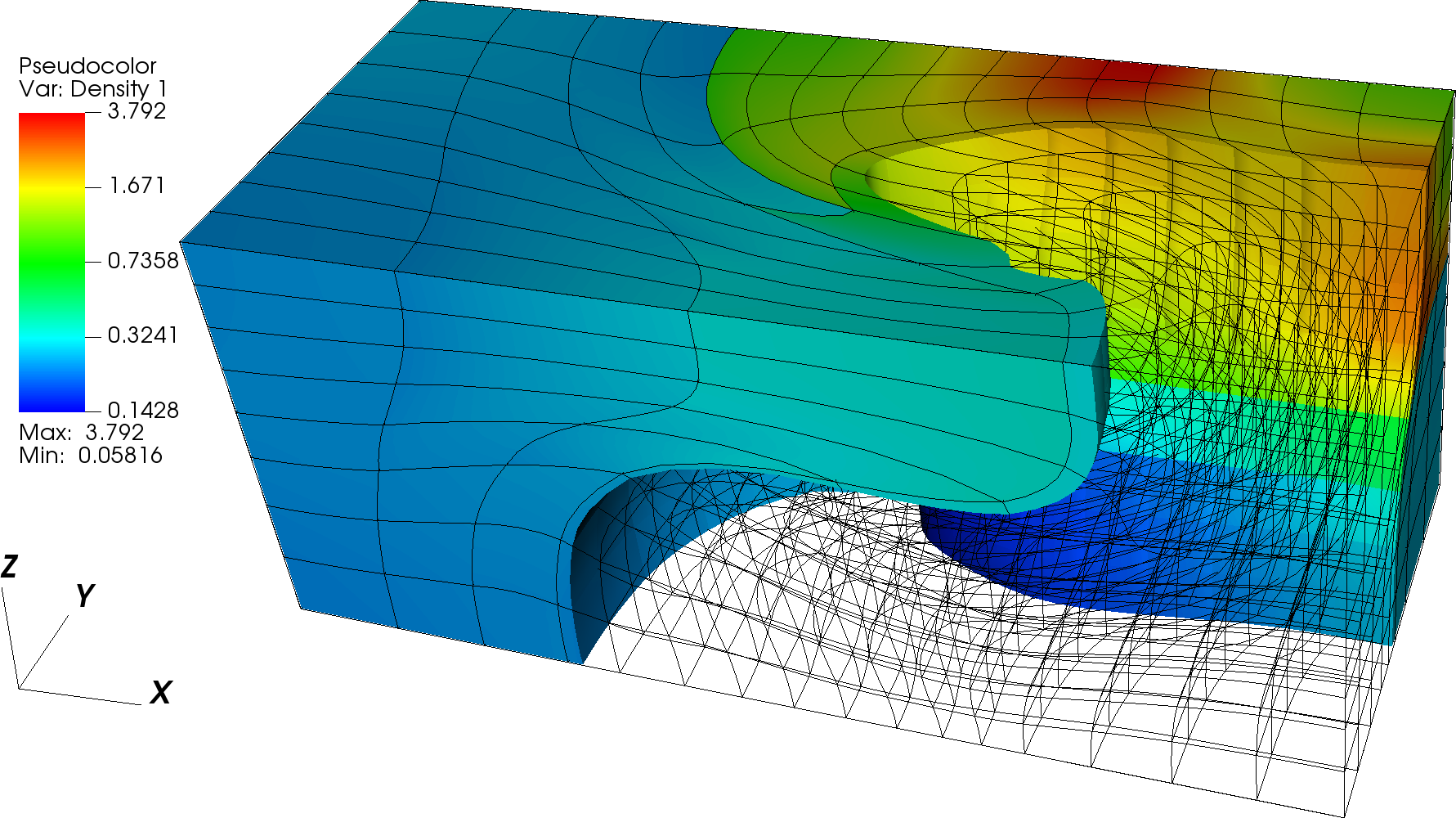}
}
\centerline
{
   \includegraphics[width=0.3\textwidth]{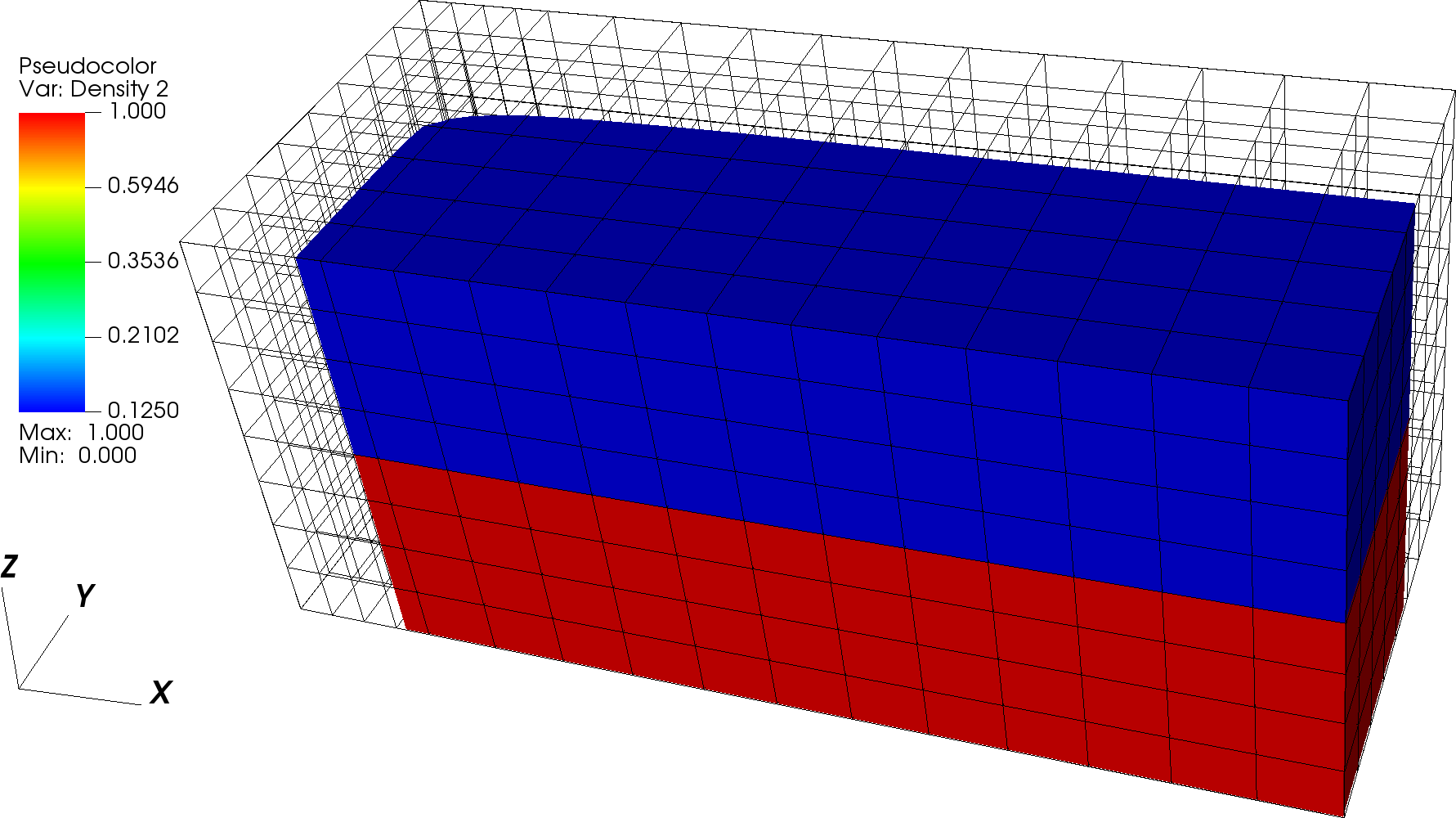}\hfil
   \includegraphics[width=0.3\textwidth]{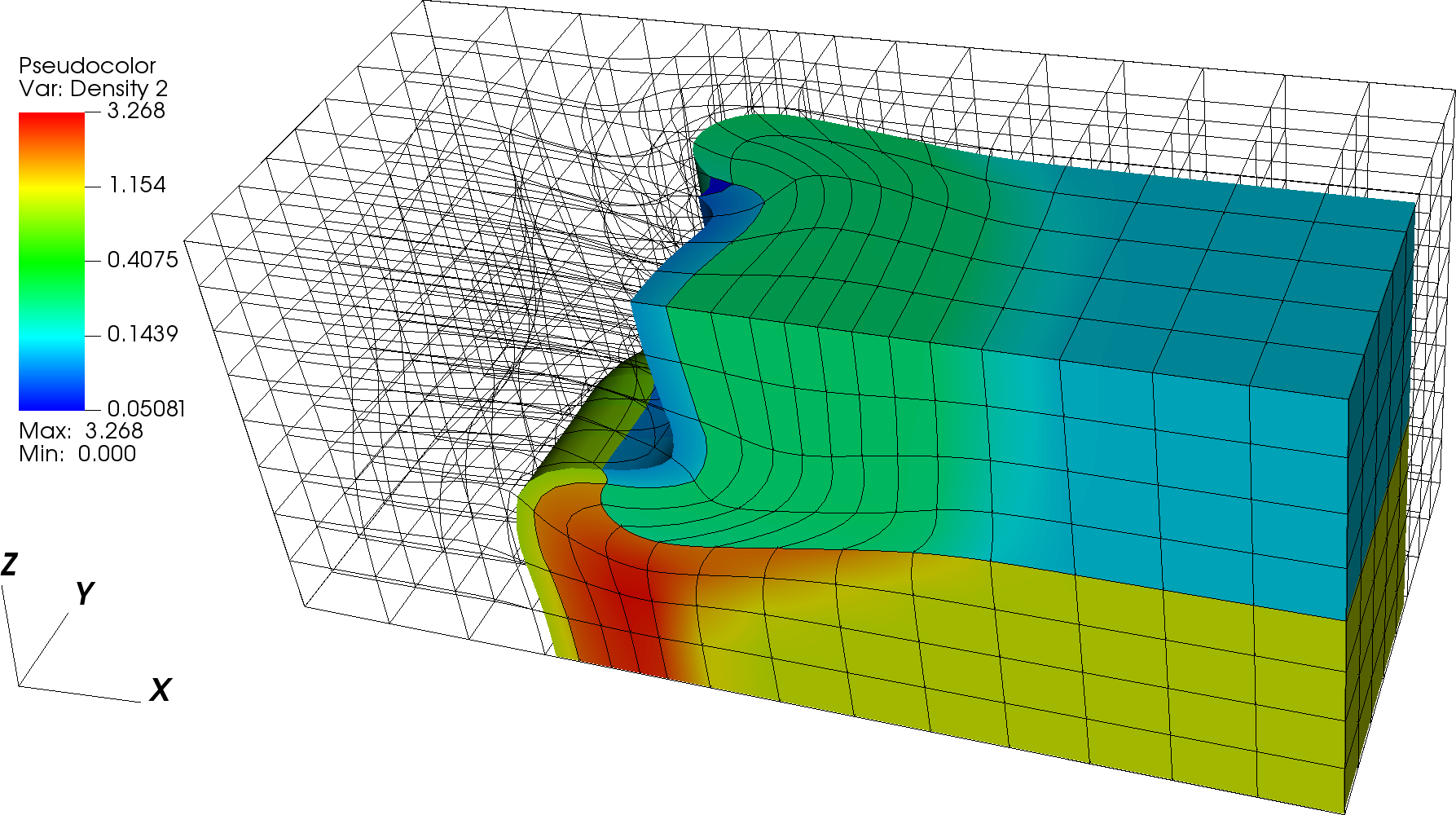}\hfil
   \includegraphics[width=0.3\textwidth]{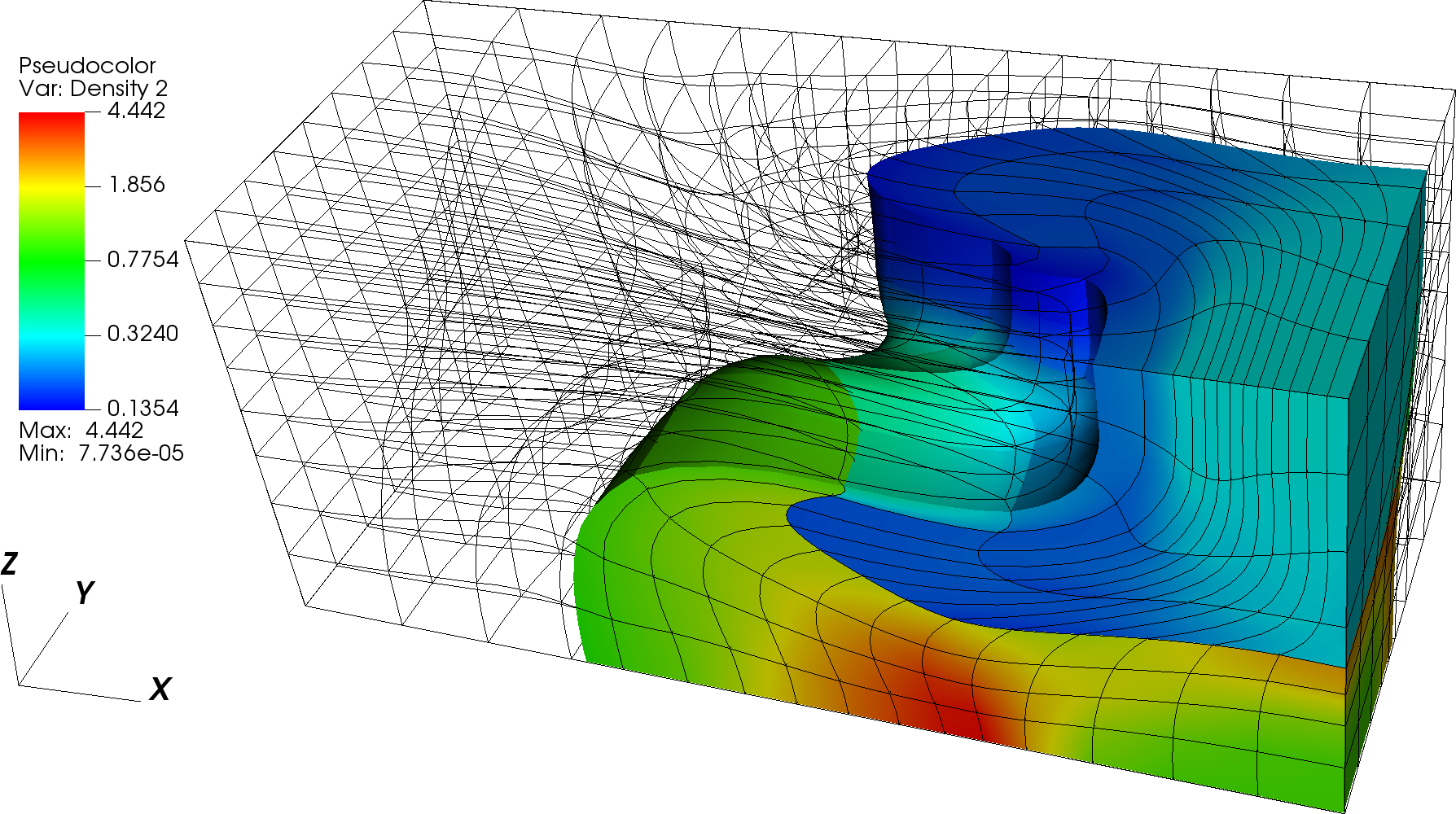}
}
\caption{Material densities (\#1 on top, \#2 on bottom)
         at times $t=0$, $t=2.5$, and $t=5$
         for the three-dimensional Triple Point interaction problem.}
\label{fig_3p_3D}
\end{figure}

The domain is $x \in [0,7] \times y \in [0,3]$ in two dimensions, and
$x \in [0,7] \times y \in [0,3] \times z \in [0,3]$ in three dimensions with
$\bs{v} \cdot \bs{n} = 0$ boundary conditions.
At $t=0$ we set a non-fitted material interface as shown in the left panels
of Figures \ref{fig_3p_2D} and Figures \ref{fig_3p_3D}.
We denote the left/top material with \#1, while the bottom/right region is
occupied by material \#2.
The three initial states are (where $z = 0$ in the two-dimensional setup):
\[
(v, \rho, p, \gamma)=
  \begin{cases}
     (0, 1, 1, 1.5)        & \text{for Material \#1}, ~ x < 1, \\
     (0, 0.125, 0.1, 1.5)  & \text{for Material \#1}, ~ x > 1, z < 1.5, \\
     (0, 1, 0.1, 1.5)      & \text{for Material \#1}, ~ x > 1, z > 1.5, \\
     (0, 1, 0.1, 1.4)      & \text{for Material \#2}, ~ z < 1.5, \\
     (0, 0.125, 0.1, 1.4)  & \text{for Material \#2}, ~ z > 1.5.
  \end{cases}
\]

Although the jumps in material properties are not large,
the resulting material interfaces are rather complex.
Density plots for both materials at times $t = 2.5$ and $t = 5$ are shown
in Figures \ref{fig_3p_2D} and \ref{fig_3p_3D}
for the two- and three-dimensional simulations, respectively.
We observe that the Shifted Interface Method is able to produce reasonable
results in both two and three dimensions, while maintaining the sharp interface
representation, as designed in the method's construction.
This test is also a good example of the method's generality with respect to dimension,
as it is based on finite elements and does not involve explicit
geometric operations.

%% file: appendix.tex
\appendix

\section{Conservation properties}
\label{sec_conservation}

To start the discussion, let us consider a body-fitted variational formulation that corresponds to the limit case of the WSIM when $\bs{d} \rightarrow \bs{0}$, so that $\ti{\G} \rightarrow \G_I$, $\ti{\bs{n}} \rightarrow \bs{n}$, $\ti{\bs{n}} \cdot \bs{n} = 1$ and $\ti{\bs{n}} \cdot \bs{\tau}_i = 0$. In this limit, \eqref{eq_weak_v}, \eqref{eq_weak_e_1}, and
\eqref{eq_weak_e_2} reduce to:
\begin{subequations}
	\label{eq:SH_Cons1A}
	\begin{align}
		\label{eq:SH_Cons1A1}
	(\bs{\psi} \, , \, \rho \dot{\bs{v}})_{\Om}
	- (\nabla \cdot \bs{\psi} \, , \, p )_{\Om}
  + (\mu \nabla^s \bs{v}, \nabla \bs{\psi})_{\Om}
	&=\;
	0 \; ,
	\\
	(\phi \, , \, \rho \dot{e} )_{\Om}
	+( \phi \, , \,  p \, \nabla \cdot \bs{v} )_{\Om}
  -	 ( \mu \nabla^s \bs{v}, \nabla \bs{v} )_{\Om}
	&=\;
	0 \; .
	\end{align}
\end{subequations}
Let us assume, as is customary when proving conservation statements, that only homogenous Neumann (zero-traction) boundary conditions are imposed. 
Setting $\bs{\psi}=\bs{e}_i$ (where $\bs{e}_i$ is the unit vector along the $i$th coordinate axis):
\begin{subequations}
	\label{eq:SH_Cons2_genA}
	\begin{align}
	\label{eq:SH_Cons2A}
	\partial_t (1 \, , \, \rho \bs{v} \cdot \bs{e}_i)_{\Om}
	=
	0
	\; , \qquad i=1 \dots d \, ,
	\end{align}
which states that global integral of the momentum is conserved over time for the semi-discrete system of equations.
\end{subequations}
The global statement of conservation of total energy is obtained setting $\phi=1$ in the energy equation, $\bs{\psi}=\bs{v}$ in the momentum equation (to obtain a kinetic energy balance), and summing the results:
	\begin{align}
	\label{eq:SH_Cons3a}
	\partial_t (1 \, , \, \rho ( \bs{v} \cdot \bs{v} /2 + e) )_{\Om}
	&=\;
	(\nabla \cdot \bs{v} \, , \, p )_{\Om}
	-( 1 \, , \,  p \, \nabla \cdot \bs{v} )_{\Om}
	= 0 \, .
	\end{align}
The conservation of momentum statement for the WSIM is achieved
analogously to the body-fitted variational form. 
Namely, setting $\bs{\psi}=\bs{e}_j$ in~\eqref{eq_weak_v} yields
\begin{subequations}
\begin{equation}
\begin{split}
\label{eq_weak_v_conserve_p1_A1}
			(\rho_{1}  , \, \dot{\bs{v}} \cdot \bs{e}_j   )_{\tO_{1}}
	  + (\rho_{2} , \, \dot{\bs{v}} \cdot \bs{e}_j   )_{\tO_{2}}
  & + (\rho_{i}  \,\alpha_{i}  , \, \dot{\bs{v}} \cdot \bs{e}_j   )_{\Om_{c}}
	=
\\ &
\langle \jump{ \sum_{k = 1}^{m-1}  \frac{\mathcal{D}^{k}_{\bs{d}} \, p^{-} }{k!}  } \cdot \bs{n}_{1},  \alpha_{2}^{-} \, \bs{e}_j  \cdot \ti{\bs{n}}^{+} \rangle _{\ti{\G}_{1}}
+\langle \jump{ \sum_{k = 1}^{m-1}  \frac{\mathcal{D}^{k}_{\bs{d}} \, p^{+} }{k!}  } \cdot \bs{n}_{2} ,  \alpha_{1}^{+} \, \bs{e}_j \cdot \ti{\bs{n}}^{-}  \rangle _{\ti{\G}_{2}}
\\ &  
+\langle \frac{\alpha^{+}_{1} - \alpha^{-}_{1}}{2}\left( \jump{ \sum_{k = 1}^{m-1}  \frac{\mathcal{D}^{k}_{\bs{d}} \, p^{-} }{k!}  } \cdot \bs{n}_{1}   +\jump{ \sum_{k = 1}^{m-1}  \frac{\mathcal{D}^{k}_{\bs{d}} \, p^{+} }{k!}  } \cdot \bs{n}_{1} \right)   , \bs{e}_j   \cdot \ti{\bs{n}}^{+} \rangle_{\mathcal{E}^{o}_{c}}\, ,
\end{split}
\end{equation}
	Adding 
$	\langle \jump{p^{-}} \cdot \bs{n}_{1},  \alpha^{-}_{2} \bs{e}_{j} \cdot \ti{\bs{n}}^{+} \rangle _{\ti{\G}_{1}}
	+ \langle \jump{p^{+}} \cdot \bs{n}_{2} ,  \alpha^{+}_{1}  \bs{e}_{j} \cdot \ti{\bs{n}}^{-} \rangle _{\ti{\G}_{2}}
	+ \avg{\frac{\alpha^{+}_{1} - \alpha^{-}_{1}}{2}\left(  \jump{p^{-}}  + \jump{p^{+}} \right) \cdot \bs{n}_{1}  \, , \, \bs{e}_j  \cdot \ti{\bs{n}}^{+}}_{\mathcal{E}^{o}_{c}} $
	to both sides of~\eqref{eq_weak_v_conserve_p1_A1} yields,
\begin{equation}
\begin{split}
\label{eq_weak_v_conserve}
	(\rho_{1}    \, , \, \dot{\bs{v}} \cdot \bs{e}_j   )_{\tO_{1}}
		& +
	(\rho_{2}   \, , \, \dot{\bs{v}} \cdot \bs{e}_j   )_{\tO_{2}}
		+
(\rho_{i}  \,\alpha_{i}  \, , \, \dot{\bs{v}} \cdot \bs{e}_j   )_{\Om_{c}}
+
			\langle \jump{p^{-}} \cdot \bs{n}_{1} \, , \,  \alpha^{-}_{2} \bs{e}_{j} \cdot \ti{\bs{n}}^{+} \rangle _{\ti{\G}_{1}}
\\  & 
		+
		\langle \jump{p^{+}} \cdot \bs{n}_{2} \, , \,   \alpha^{+}_{1}  \bs{e}_{j} \cdot \ti{\bs{n}}^{-} \rangle _{\ti{\G}_{2}}
	+ \avg{\frac{\alpha^{+}_{1} - \alpha^{-}_{1}}{2}\left(  \jump{p^{-}}  + \jump{p^{+}} \right) \cdot \bs{n}_{1}  \, , \, \bs{e}_j  \cdot \ti{\bs{n}}^{+}}_{\mathcal{E}^{o}_{c}} 
=
		e_{mom}(\bs{e}_{j}) \, ,
\end{split}
\end{equation}
where
\begin{equation}
\begin{split}
		e_{mom}(\bs{e}_{j})
		=&\; 
\langle \jump{ \tS^{m-1} p^{-}  } \cdot \bs{n}_{1},  \alpha_{2}^{-} \, \bs{e}_j  \cdot \ti{\bs{n}}^{+} \rangle _{\ti{\G}_{1}}
		+\langle \jump{ \tS^{m-1} p^{+} } \cdot \bs{n}_{2} ,  \alpha_{1}^{+} \, \bs{e}_j \cdot \ti{\bs{n}}^{-}  \rangle _{\ti{\G}_{2}} 
  \\ &  \;
	+\langle \frac{\alpha^{+}_{1} - \alpha^{-}_{1}}{2}  \jump{ \tS^{m-1} (p^{-} + p^{+})  } \cdot \bs{n}_{1}   , \bs{e}_j   \cdot \ti{\bs{n}}^{+} \rangle_{\mathcal{E}^{o}_{c}} \; .
\end{split}
\end{equation}
\end{subequations}
The face terms in the left-hand side of \eqref{eq_weak_v_conserve_p1_A} represent the fluxes entering and exiting $\tO_{1}$, $\tO_{2}$ and $\Om_{c}$, respectively, due to pressure jump between material $\#1$ and $\#2$.
Therefore the left-hand side of \eqref{eq_weak_v_conserve_p1_A} represents the momentum balance in the whole domain.
The presence of the material pressure jump terms across $\mathcal{E}^{o}_{c}$ is due to the fact that the material pressures are not restricted to just $\tO_{1}$ and $\tO_{2}$, but are extended to $\Om_{c}$.
Finally, from~\cite{TheoreticalPoissonAtallahCanutoScovazzi2020},
$e_{mom} \rightarrow 0$ in the asymptotic limit at a rate of $O(h^{m})$.
In other words, the formulation admits a momentum error that is caused by the Taylor expansion approximation; this error converges to zero with the order of the Taylor expansion.

A statement of global total energy conservation can be achieved summing the energy equation, tested with $\phi=1$, and the momentum equation, tested with $\bs{\psi}=\bs{v}$:
\begin{subequations}
\begin{equation}
\begin{split}
\label{eq:SH_Cons4b_v}
	&	\partial_t (1 \, , \, \rho_{1} ( \bs{v} \cdot \bs{v} /2 + e_{1}) )_{\tO_{1}}
		+ 	\partial_t (1 \, , \, \rho_{2} ( \bs{v} \cdot \bs{v} /2 + e_{2}) )_{\tO_{2}}
			+ 	\partial_t (1 \, , \, \alpha_{1} \, \rho_{1} ( \bs{v} \cdot \bs{v} /2 + e_{1}) )_{\Om_{c}}
    \\ &	
		+ 	\partial_t (1 \, , \, \alpha_{2} \, \rho_{2} ( \bs{v} \cdot \bs{v} /2 + e_{2}) )_{\Om_{c}}
	= 
		\avg{  \jump { \bs{S}_{h}^{m} \bs{v}} \cdot \bs{n}_{1} (\ti{\bs{n}}^{+} \cdot \bs{n}_{1})   ,  \{\alpha_{1}  p_{1}\}_{\gamma_{1}}+\{\alpha_{2}  p_{2}\}_{\gamma_{2}} }_{\ti{\G}_{1} \cup \ti{\G}_{2} \cup \mathcal{E}^{o}_{c}}
		\\ &	
		+ 
		\langle \jump{ \sum_{k = 1}^{m-1}  \frac{\mathcal{D}^{k}_{\bs{d}} \, p^{-} }{k!}  } \cdot \bs{n}_{1},  \alpha^{-}_{2} \, \bs{v} \cdot \ti{\bs{n}}^{+} \rangle _{\ti{\G}_{1}}
		+ \langle \jump{ \sum_{k = 1}^{m-1}  \frac{\mathcal{D}^{k}_{\bs{d}} \, p^{+} }{k!}  } \cdot \bs{n}_{2},  \alpha^{+}_{1} \, \bs{v}  \cdot \ti{\bs{n}}^{-}\rangle _{\ti{\G}_{2}}
		\\ &
		+\langle \frac{\alpha^{+}_{1} - \alpha^{-}_{1}}{2}\left( \jump{ \sum_{k = 1}^{m-1}  \frac{\mathcal{D}^{k}_{\bs{d}} \, p^{-} }{k!}  } \cdot \bs{n}_{1}   +\jump{ \sum_{k = 1}^{m-1}  \frac{\mathcal{D}^{k}_{\bs{d}} \, p^{+} }{k!}  } \cdot \bs{n}_{1} \right)    , \bs{v}   \cdot \ti{\bs{n}}^{+} \rangle_{\mathcal{E}^{o}_{c}}
		\, ,
\end{split}
\end{equation}
Adding 
$	\langle \jump{p^{-}} \cdot \bs{n}_{1},  \alpha^{-}_{2} \bs{v} \cdot \ti{\bs{n}}^{+} \rangle _{\ti{\G}_{1}}
+ \langle \jump{p^{+}} \cdot \bs{n}_{2} ,  \alpha^{+}_{1}  \bs{v}\cdot \ti{\bs{n}}^{-} \rangle _{\ti{\G}_{2}}
+ \avg{\frac{\alpha^{+}_{1} - \alpha^{-}_{1}}{2}\left(  \jump{p^{-}}  + \jump{p^{+}} \right) \cdot \bs{n}_{1}  \, , \, \bs{v} \cdot \ti{\bs{n}}^{+}}_{\mathcal{E}^{o}_{c}} $
to both sides of~\eqref{eq:SH_Cons4b_v} yields,
\begin{equation}
\begin{split}
\label{eq:SH_Cons4b-A}
	&	\partial_t (1 \, , \, \rho_{1} ( \bs{v} \cdot \bs{v} /2 + e_{1}) )_{\tO_{1}}
+ 	\partial_t (1 \, , \, \rho_{2} ( \bs{v} \cdot \bs{v} /2 + e_{2}) )_{\tO_{2}}
+ 	\partial_t (1 \, , \, \alpha_{1} \, \rho_{1} ( \bs{v} \cdot \bs{v} /2 + e_{1}) )_{\Om_{c}}
\\ &	
+ 	\partial_t (1 \, , \, \alpha_{2} \, \rho_{2} ( \bs{v} \cdot \bs{v} /2 + e_{2}) )_{\Om_{c}}
+ \langle \jump{p^{-}} \cdot \bs{n}_{1} \, , \,  \alpha^{-}_{2} \bs{v} \cdot \ti{\bs{n}}^{+} \rangle _{\ti{\G}_{1}}
+ \langle \jump{p^{+}} \cdot \bs{n}_{2} \, , \,   \alpha^{+}_{1}  \bs{v} \cdot \ti{\bs{n}}^{-} \rangle _{\ti{\G}_{2}}
\\ &	
	+ \avg{\frac{\alpha^{+}_{1} - \alpha^{-}_{1}}{2}\left(  \jump{p^{-}}  + \jump{p^{+}} \right) \cdot \bs{n}_{1}  \, , \, \bs{v}  \cdot \ti{\bs{n}}^{+}}_{\mathcal{E}^{o}_{c}} 
\\ &	
= 
\avg{  \jump { \bs{S}_{h}^{m} \bs{v}} \cdot \bs{n}_{1} (\ti{\bs{n}}^{+} \cdot \bs{n}_{1})   ,  \{\alpha_{1}  p_{1}\}_{\gamma_{1}}+\{\alpha_{2}  p_{2}\}_{\gamma_{2}} }_{\ti{\G}_{1} \cup \ti{\G}_{2} \cup \mathcal{E}^{o}_{c}}
+ e_{mom}(\bs{v})
	\, .
\end{split}
\end{equation}
Similarly to the momentum derivation, the left-hand side of \eqref{eq:SH_Cons4b-A}
represents the total energy balance in the whole domain. 
In the asymptotic limit, from~\cite{TheoreticalPoissonAtallahCanutoScovazzi2020}, $\jump { \bs{S}_{h}^{m} \bs{v}} \cdot \bs{n}_{1} \rightarrow 0$ at a rate of $O(h^{m+1})$ and $e_{mom} \rightarrow 0$ at a rate of $O(h^{m})$.
Therefore the formulation admits a total energy error that goes to zero with the order of the Taylor expansion.

The form of the momentum \eqref{eq_weak_v_conserve} and the
total energy \eqref{eq:SH_Cons4b-A} statements are not traditional.
A different way to interpret these statements is to use the Gauss divergence theorem on $(\alpha_{i} \nabla p_{i} , \bs{e}_{j})_{\Om_c}$. Applying~\eqref{eq:continuity_p1}-~\eqref{eq:continuity_p2} and following the same steps as Section~\ref{sec_weak_form} yields
\begin{equation}
\begin{split}
\label{eq:Gauss}
(\alpha_{i} \nabla p_{i} , \bs{e}_{j})_{\Om_c}
& = 
\avg{ \jump{p^{-}} \cdot \bs{n}_{1}  \, , \, \alpha_{2}^{-}\, \bs{e}_j \cdot \ti{\bs{n}}^{+}}_{\ti{\G}_{1}}
+ \langle p^{-}_{1} , \bs{e}_j  \cdot \ti{\bs{n}}^{-}  \rangle _{\ti{\G}_{1}}
	\\ &
+\avg{ \jump{p^{+}} \cdot \bs{n}_{2}  \, , \, \alpha_{1}^{+}\, \bs{e}_j  \cdot \ti{\bs{n}}^{-}}_{\ti{\G}_{2}} 
+ \langle p^{+}_{2} , \bs{e}_j  \cdot \ti{\bs{n}}^{+}  \rangle _{\ti{\G}_{2}} 
+ \avg{\frac{\alpha^{+}_{1} - \alpha^{-}_{1}}{2}\left(  \jump{p^{-}}  + \jump{p^{+}} \right)\cdot \bs{n}_{1}  \, , \, \bs{e}_j  \cdot \ti{\bs{n}}^{+}}_{\mathcal{E}^{o}_{c}}
\end{split}
\end{equation}
Thus, plugging~\eqref{eq:Gauss} in~\eqref{eq_weak_v_conserve} yields
	\begin{align}
\label{eq_weak_v_conserve_p1_A}
(\rho_{1}   \, , \, \dot{\bs{v}} \cdot \bs{e}_j   )_{\tO_{1}}
+ (\rho_{2}   \, , \, \dot{\bs{v}} \cdot \bs{e}_j   )_{\tO_{2}}
- \langle p^{-}_{1} , \bs{e}_j  \cdot \ti{\bs{n}}^{-}  \rangle _{\ti{\G}_{1}}
- \langle p^{+}_{2} , \bs{e}_j  \cdot \ti{\bs{n}}^{+}  \rangle _{\ti{\G}_{2}} 
& = \; 
e_{mom}(\bs{e}_{j}) 	- \underbrace{(\rho_{i}  \,\dot{\bs{v}}  +  \nabla p_{i}\, , \, \alpha_{i} \,  \bs{e}_j   )_{\Om_{c}}}_{O(h^{m})}\, ,
\end{align}
where
\begin{equation}
\begin{split}
		e_{mom}(\bs{e}_{j})
		& = \; 
 \langle \jump{ \tS^{m-1} p^{-} } \cdot \bs{n}_{1},  \alpha_{2}^{-} \, \bs{e}_j  \cdot \ti{\bs{n}}^{+} \rangle _{\ti{\G}_{1}}s
		+\langle \jump{ \tS^{m-1} p^{+}  } \cdot \bs{n}_{2} ,  \alpha_{1}^{+} \, \bs{e}_j \cdot \ti{\bs{n}}^{-}  \rangle _{\ti{\G}_{2}} 
\\ &  \;
	+\langle \frac{\alpha^{+}_{1} - \alpha^{-}_{1}}{2}  \jump{ \tS^{m-1} (p^{-} + p^{+})  } \cdot \bs{n}_{1}   , \bs{e}_j   \cdot \ti{\bs{n}}^{+} \rangle_{\mathcal{E}^{o}_{c}}
\\ &  \;
		= O(h^{m})
		\; .
\end{split}
\end{equation}
In regards to the statement of global total energy conservation~\eqref{eq:SH_Cons4b-A}, we first rewrite it as
\begin{equation}
\begin{split}
\label{eq:SH_Cons4b}
\partial_t (1 \, , \, \rho ( \bs{v} \cdot \bs{v} /2 + e_{1}) )_{\tO_{1}}
 & + \partial_t (1 \, , \, \rho ( \bs{v} \cdot \bs{v} /2 + e_{2}) )_{\tO_{2}}
= 
-(\rho_{1} \dot{e}_{1} + p_{1} \nabla  \cdot \bs{v} , \alpha_{1}  )_{\Om_{c}}
-	(\rho_{2} \dot{e}_{2} + p_{2} \nabla  \cdot \bs{v} , \alpha_{2}  )_{\Om_{c}}
\\	&	\; 
+ \avg{  \jump { \bs{S}_{h}^{m} \bs{v}} \cdot \bs{n}_{1} (\ti{\bs{n}}^{+} \cdot \bs{n}_{1})   ,  \{\alpha_{1}  p_{1}\}_{\gamma_{1}}+\{\alpha_{2}  p_{2}\}_{\gamma_{2}} }_{\ti{\G}_{1} \cup \ti{\G}_{2} \cup \mathcal{E}^{o}_{c}}
\\	&	\; 
-	( \alpha_{i} \, \rho_{i}  , \dot {\bs{v}} \cdot \bs{v} )_{\Om_{c}} 
+	( \alpha_{i} \, p_{i} \, , \, \nabla \cdot \bs{v}  )_{\Om_{c}} 
\\ &	\; 
+ \langle \jump{ \sum_{k = 1}^{m-1}  \frac{\mathcal{D}^{k}_{\bs{d}} \, p^{-} }{k!}  } \cdot \bs{n}_{1},  \alpha^{-}_{2} \, \bs{v} \cdot \ti{\bs{n}}^{+} \rangle _{\ti{\G}_{1}}
+ \langle \jump{ \sum_{k = 1}^{m-1} \frac{\mathcal{D}^{k}_{\bs{d}} \, p^{+} }{k!}  } \cdot \bs{n}_{2},  \alpha^{+}_{1} \, \bs{v}  \cdot \ti{\bs{n}}^{-}\rangle _{\ti{\G}_{2}}
\\ &	\; 
+\langle \frac{\alpha^{+}_{1} - \alpha^{-}_{1}}{2}\left( \jump{ \sum_{k = 1}^{m-1}  \frac{\mathcal{D}^{k}_{\bs{d}} \, p^{-} }{k!}  } \cdot \bs{n}_{1}   +\jump{ \sum_{k = 1}^{m-1}  \frac{\mathcal{D}^{k}_{\bs{d}} \, p^{+} }{k!}  } \cdot \bs{n}_{1} \right)    , \bs{v}   \cdot \ti{\bs{n}}^{+} \rangle_{\mathcal{E}^{o}_{c}}
\, ,
\end{split}
\end{equation}
The first three terms on the right-hand side approach $0$ at a rate of $O(h^{m+1})$.
At this point, we employ Gauss divergence theorem for $( \alpha_{i} \, p_{i} \, , \, \nabla \cdot \bs{v}  )_{\Om_{c}} $,  recall~\eqref{eq:continuity_p1}-~\eqref{eq:continuity_p2} and follow the same steps as Section~\ref{sec_weak_form}, yielding
\begin{equation}
\begin{split}
\label{eq:SH_Cons4b-2}
&	\partial_t (1 \, , \, \rho_{1} ( \bs{v} \cdot \bs{v} /2 + e_{1}) )_{\tO_{1}}
+ 	\partial_t (1 \, , \, \rho_{2} ( \bs{v} \cdot \bs{v} /2 + e_{2}) )_{\tO_{2}}
- \langle p^{-}_{1} , \bs{v}  \cdot \ti{\bs{n}}^{-}  \rangle _{\ti{\G}_{1}}
- \langle p^{+}_{2} , \bs{v} \cdot \ti{\bs{n}}^{+}  \rangle _{\ti{\G}_{2}} 
 = 
	\\ & 	\; 
  O(h^{m+1})
	- \underbrace{(\rho_{i}  \,\dot{\bs{v}}  +  \nabla p_{i}\, , \, \alpha_{i} \,  \bs{v} )_{\Om_{c}}}_{O(h^{m})}
+e_{mom}(\bs{v})
 = 
O(h^{m})
\, .
\end{split}
\end{equation}
\end{subequations}